\documentclass{amsart}
\usepackage[colorlinks, linkcolor=blue,anchorcolor=Periwinkle,
    citecolor=blue,urlcolor=Emerald]{hyperref}
\usepackage{amsmath}
\usepackage{amssymb}
\usepackage{mathrsfs}
\usepackage{amscd}
\usepackage{amsthm,color}
\usepackage[all]{xypic}
\SelectTips{cm}{}
\usepackage{tikz}
\usepackage{dirtytalk}
\usepackage{mathtools}

\let\margin\marginpar
\newcommand\myMargin[1]{\margin{\raggedright\scriptsize #1}}
\renewcommand{\marginpar}[1]{\myMargin{#1}}

\newtheorem{lemma}{Lemma}[section]
\newtheorem{theorem}[lemma]{Theorem}
\newtheorem{corollary}[lemma]{Corollary}
\newtheorem{conjecture}[lemma]{Conjecture}
\newtheorem{prop}[lemma]{Proposition}
\theoremstyle{definition}
\newtheorem{definition}[lemma]{Definition}
\newtheorem{example}[lemma]{Example}

\newtheorem{remark}[lemma]{Remark}

\theoremstyle{remark}
\newtheorem*{proof*}{Proof}
\numberwithin{equation}{section}

\newcommand{\iso}{\stackrel{_\sim}{\rightarrow}}
\newcommand{\liso}{\xleftarrow{_\sim}}
\newcommand{\ko}{\;,}

\newcommand{\colim}{\mbox{colim}}

\def\Spec{{\bf {Spec}}}
\def\RHom{{\mathbf{R}\mathrm{Hom}}}

\def\D{\mathrm{D}}
\def\ev{\mathrm{ev}}
\def\Ev{\mathrm{Ev}}
\def\iD{\mathcal{D}}

\def\Ext{{\mathrm{Ext}}}

\def\Hom{{\mathrm{Hom}}}

\def\End{{\mathrm{End}}}

\def\per{{\mathrm{per}}}

\def\twoper{{\mbox{2-per}}}
\def\Spec{{\mathrm{Spec\ }}}
\def\deg{{\mathrm{deg}}}

\def\res{{\mathrm{res}}}
\def\tria{{\mathrm{tria}}}
\def\Tot{{\mathrm{Tot}}}

\def\2silt{{2\mathrm{silt}}}
\def\del{{\partial}}

\def\PP{{\mathbb P}}

\def\ZZ{{\mathbb Z}}
\def\CC{{\mathbb C}}

\def\HH{{\mathbb H}}
\def\DD{{\mathbf D}}

\def\HoH{{\mathrm{HH}}}
\def\HC{{\mathrm{HC}}}

\def\hY{{\widehat{Y}}}
\def\hX{{\widehat{X}}}
\def\hf{{\widehat{f}}}
\def\bR{{\bf{R}}}

\def\Ga{{\Gamma}}
\def\La{{\Lambda}}

\def\cF{{\mathcal{F}}}
\def\cE{{\mathcal{E}}}
\def\cO{{\mathcal{O}}}
\def\cC{{\mathcal{C}}}
\def\cL{{\mathcal{L}}}
\def\cM{{\mathcal{M}}}

\def\cN{{\mathcal{N}}}
\def\cA{{\mathcal{A}}}
\def\cB{{\mathcal{B}}}
\def\cH{{\mathcal{H}}}

\def\cT{{\mathcal{T}}}
\def\cU{{\mathcal{U}}}

\def\cR{{\mathcal{R}}}
\def\cQ{{\mathcal{Q}}}

\def\cX{{\mathcal{X}}}
\def\cY{{\mathcal{Y}}}

\def\fm{{\mathfrak{m}}}

\def\lgg{{\langle\langle}}
\def\rgg{{\rangle\rangle}}
\def\ol{\overline}
\def\ul{\underline}
\def\ot{{\otimes}}
\def\bt{{\boxtimes}}
\def\ten{{\otimes}}
\def\lten{{\overset{\mathbb{L}}{\ten}}}

\def\op{{\oplus}}
\def\wh{\widehat}
\def\wt{\widetilde}

\def\sC{{\mathscr{C}}}

\def\bD{{\mathbb{D}}}

\def\deg{{\mathrm{deg}}}
\def\MF{{\mathrm{MF}}}

\def\Lm{{\Lambda}}

\def\p{{\prime}}
\def\pp{{\prime\prime}}
\def\1{{\bf{1}}}

\def\coh{{\mathrm{coh}}}
\def\Qcoh{{\mathrm{Qcoh}}}
\def\Mod{{\mathrm{Mod}~}}
\def\mod{{\mathrm{mod}~}}
\def\cy{{\mathrm{cyc}}}

\def\id{{\mathrm{Id}}}

\def\fD{{\mathfrak{D}}}

\def\proj{{\mathrm{proj}}}

\def\add{{\mathrm{add}}}

\def\CM{{\mathrm{CM}}}
\def\thick{{\mathbf{thick}}}

\def\Sg{\Sigma}
\def\Def{{\mathrm{Def}}}

\def\DEF{{\mathrm{DEF}}}

\def\Art{{\mathbf{Art}}}
\def\Set{{\mathbf{Set}}}
\def\cArt{{\mathbf{cArt}}}
\def\dgArt{{\mathbf{dgArt}}}

\def\Ex{{\mathrm{Ex}}}
\def\rad{{\mathrm{rad}}}
\def\dd{{\mathrm{d}}}
\def\kiu{{k[u^{-1}]}}
\def\kuiu{{k[u,u^{-1}]}}
\def\Gpd{{\bf{Gpd}}}
\def\tp{{2^\p}}
\def\ut{{u^{-1}\text{-}tor}}

\def\PCAlgc{{\mathrm{PCAlgc}}}
\def\IndQ{{\mathrm{IndCoh}}}

\newcommand{\dgcat}{\mathrm{dgcat}}
\newcommand{\St}{\mathcal{S}t}
\renewcommand{\Pr}{\mathcal{P}r}
\newcommand{\Tria}{\mathrm{Tria}}

\newcommand{\eps}{\varepsilon}

\begin{document}
\title{Cluster categories and rational curves}

\author{Zheng Hua}
\address{Department of Mathematics\\
	The University of Hong Kong\\
	Hong Kong SAR\\
	China
}
\email{huazheng@maths.hku.hk}
\urladdr{http://hkumath.hku.hk/~huazheng/}

\author{Bernhard Keller}
\address{Universit\' Paris Cit\'e and Sorbonne Universit\'e, CNRS, IMJ-PRG, F-75013 Paris, France}
\email{bernhard.keller@imj-prg.fr}
\urladdr{https://webusers.imj-prg.fr/~bernhard.keller/}

\begin{abstract}
We study rational curves on smooth complex Calabi--Yau threefolds via noncommutative algebra. By the general theory of derived noncommutative deformations due to Efimov, Lunts and Orlov, 
the structure sheaf of a rational curve in a smooth CY 3-fold $Y$ is pro-represented by a nonpositively graded dg algebra $\Ga$. The curve is called \emph{nc rigid} if $H^0\Ga$ is finite dimensional. When $C$ is contractible, $H^0\Ga$ is isomorphic to the contraction algebra defined by Donovan and Wemyss. 
More generally, one can show that there exists a $\Ga$ pro-representing the (derived) multi-pointed deformation (defined by Kawamata) of a collection of rational curves $C_1,\ldots,C_t$ with $\dim(\Hom_Y(\cO_{C_i},\cO_{C_j}))=\delta_{ij}$. The collection is called \emph{nc rigid} if $H^0\Ga$ is finite dimensional.
We prove that $\Ga$ is a homologically smooth bimodule 3CY algebra. As a consequence, we define a (2CY) cluster category $\cC_\Ga$ for such a collection of rational curves in $Y$. It has finite-dimensional morphism spaces 
iff the collection is nc rigid. When $\bigcup_{i=1}^tC_i$ is (formally) contractible by a morphism $\hY\to \hX$, 
then $\cC_\Ga$ is equivalent to the singularity category of $\hX$
and thus categorifies the contraction algebra of 
Donovan and Wemyss.  The Calabi-Yau structure on $Y$ determines a canonical class $[w]$
(defined up to right equivalence) in the zeroth Hochschild homology of $H^0\Ga$. 
Using our previous work on the noncommutative Mather--Yau theorem and singular Hochschild cohomology, we prove that the singularities underlying a 3-dimensional smooth flopping contraction are classified by the derived equivalence
class of the pair $(H^0\Ga, [w])$. We also give a new necessary condition for contractibility of rational curves in terms of $\Ga$.
\end{abstract}



\maketitle

\section{Introduction}
\newtheorem{maintheorem}{\bf{Theorem}}
\renewcommand{\themaintheorem}{\Alph{maintheorem}}
\newtheorem{mainconjecture}[maintheorem]{\bf{Conjecture}}
\renewcommand{\themainconjecture}{}

The study of rational curves in algebraic varieties lies at the core of birational geometry. A smooth rational curve $C$ in a quasi-projective variety $Y$ is called \emph{rigid} if the component of the Hilbert scheme of curves containing $C$ is a finite scheme. Note that this is weaker than the notion of \emph{infinitesimally rigid}, which says that $\Ext^1_Y(\cO_C,\cO_C)=0$.  If a curve is not rigid then we call it \emph{movable}. When $Y$ is a smooth projective surface, a smooth rational curve $C\subset Y$ is rigid if and only if its normal bundle $N_{C/Y}$ is negative. And if $C$ is rigid then it is \emph{contractible}, i.e. for the formal completion 
$\hY$ of $Y$ along $C$ there exists a birational morphism $f:\hY\to \hX$ to a normal surface $\hX$ that contracts $C$. The definition of contractibility in general can be found in Definition \ref{def:contractible}.

In this article, we will focus on the case when $Y$ is a smooth complex Calabi--Yau threefold, i.e. $\omega_Y$ is trivial. The situation is much more complicated than the surface case. We call a rational curve $C\subset Y$ of type $(a,b)$ if it has normal bundle $\cO(a)\op \cO(b)$. By the adjunction formula, we have $a+b=-2$. 
A $(-1,-1)$-curve is contractible. The underlying singular variety $\hX$ is 
equivalent to the singular hypersurface $x^2+y^2+u^2+v^2=0$. There exists a different resolution $\hY^+\to \hX$ and the birational map $\hY\dasharrow \hY^+$ is called the \emph{Atiyah flop}. In \cite{Reid}, Reid proves that a $(0,-2)$-curve is either contractible or movable. The contractible case corresponds to the \emph{Pagoda flops}. Laufer proves that a contractible curve is of the types $(-1,-1)$, $(0,-2)$ or $(1,-3)$ (c.f. \cite{Pink83}).  Katz and Morrison show that any simple flopping contraction (see definition in Section \ref{sec:flop}) can be constructed as base change of a universal contraction \cite{KM92}.
 In general, it is not true that all rigid curves are contractible. A counter example was constructed by Clemens \cite{Clem89}.

We study the contractibility of rational curves in Calabi-Yau 3-folds via 
noncommutative methods. In general, given a rational curve $C\subset Y$ the problem is two-fold:
\begin{enumerate}
\item[$(1)$] Find infinitesimal criteria for the contractibility of $C$.
\item[$(2)$] If $C$ is contractible, determine the underlying singularity of the contraction.
\end{enumerate}
Our research is motivated by a remarkable paper of Donovan and Wemyss. In \cite{DW13}, Donovan and Wemyss considered the algebra $\Lm$ that represents the noncommutative deformation functor of $\cO_C$ for a contractible rational curve $C\subset Y$. They prove that $\Lm$ is finite dimensional and call it the \emph{ contraction algebra}. Indeed, the contraction algebra can be defined in a more general context where $Y$ may be neither CY nor smooth, and the birational morphism may contract a divisor containing $C$. However we will focus on the special case when $Y$ is a smooth CY 3-fold and the contraction is not divisorial. Donovan and Wemyss conjectured that the 3-dimensional simple flops are classified by the isomorphism types of the contraction algebras 
(cf. Conjecture \ref{DWconj}). 

In order to deal with the case of general flops where the exceptional fiber can have multiple irreducible components, Kawamata proposes to study the multi-pointed noncommutative deformation of a \emph{semi-simple collection} (see definition in Section \ref{sec:pre})
 \footnote{Kawamata calls it a \emph{simple collection}.}
of sheaves $\cE_1,\ldots, \cE_t$.  A case of special interest is when the collection is $\cO_{C_1},\ldots,\cO_{C_t}$ where $C_i$ are irreducible components of the reduced exceptional fiber of a contraction (see Example 6.5 of \cite{KaComp}).
We consider the derived noncommutative deformation theory of $\cE:=\bigoplus_{i=1}^t\cE_i$ for a semi-simple collection of sheaves $\{\cE_i\}_{i=1}^t$ on a smooth CY 3-fold $Y$. By a result of Efimov, Lunts and Orlov (cf. Theorem \ref{rep-ELO}), such deformation functor is pro-represented by a nonpositively graded dg algebra $\Ga$. We call $\Ga$ the \emph{derived deformation algebra} of the semi-simple collection $\{\cE_i\}_{i=1}^t$. We call a semi-simple collection $\{\cE_i\}_{i=1}^t$  
\emph{nc rigid} (``nc'' stands for noncommutative)  if $H^0\Ga$ is finite dimensional. Given a collection of smooth rational curve $C_1,\ldots,C_t$ such that $\{\cO_{C_i}\}_{i=1}^t$ form a semi-simple collection, we call $\{C_i\}_{i=1}^t$ a \emph{nc-rigid} collection of rational curves if $H^0\Ga$ is finite dimesional. If $t=1$ and $C=C_1$ is  nc rigid rational curve then the abelianization of $H^0\Ga$ represents the commutative deformation functor of $\cO_C$. Therefore, a nc rigid curve is in particular rigid. Our first result is:
\begin{maintheorem}\label{thmA}(Corollary \ref{GaCY})
Let $C_1,\ldots,C_t$ be a collection of rational curves   in a smooth quasi-projective Calabi-Yau 3-fold $Y$ such that $\{\cO_{C_i}\}_{i=1}^t$ form a semi-simple collection. The derived deformation algebra $\Ga$ of $\bigoplus_{i=1}^t\cO_{C_i}$ is a non positive pseudo-compact dg algebra that is 
\begin{enumerate}
\item[$(1)$] homologically smooth;
\item[$(2)$] bimodule 3CY.
\end{enumerate}
Moreover, $\Ga$ is exact 3CY in either one of the following cases
\begin{enumerate}
\item[$(a)$] $Y$ is projective;
\item[$(b)$] there is a (formal) contraction $\hf:\hY\to \hX$ such that  $\Ex(\hf)=\bigcup_{i=1}^t C_i$, where $\Ex(\hf)$ stands for the reduced exceptional fiber of $\hf$.
\end{enumerate}
\end{maintheorem}
This theorem establishes a link between birational geometry and the theory of cluster categories. We consider the triangle quotient $\cC_\Ga:=\per(\Ga)/\D_{fd}(\Ga)$ (see Section \ref{subsec:cluster category}).  It is $\Hom$-finite if and only if $\{C_i\}_{i=1}^t$ is nc rigid. 
By a result of Amiot \cite{Am}, it is then a 2CY category. When $C$ is contractible
by a morphism $\hY \to \hX$, then $H^0\Ga$ is isomorphic to the contraction algebra $\Lm$ defined in \cite{DW13} and $\cC_\Ga$ is equivalent to the singularity
category of $\hX$ and thus categorifies the contraction algebra of Donovan and Wemyss. If $\Ga$ is exact 3CY, Van den Bergh proved that it is quasi-isomorphic to a (complete) Ginzburg algebra $\fD(Q,w)$ for some finite quiver $Q$ and a potential $w$. If we fix the CY structure on $Y$, then there is a canonical class $[w]$, 
defined up to right equivalence,  in the zeroth Hochschild homology of $H^0\Ga$ (see Proposition \ref{canonicalclass}).  The class $[w]$ can be viewed as the ``classical shadow'' of the Calabi--Yau structure on $Y$. Our second result is: 

\begin{maintheorem}\label{thmB}(Theorem \ref{weakDW})
Let $\hf:\hY\to \hX$ and $\hf^\p: \hY^\p\to \hX^\p$ be two formal flopping contractions with reduced exceptional fibers $\Ex(\hf)=\bigcup_{i=1}^t C_i$ and $\Ex(\hf^\p)=\bigcup_{i=1}^s C^\p_i$. Denote respectively by $\Ga$ and $\Ga^\p$
the derived deformation algebras of $\bigoplus_{i=1}^t \cO_{C_i}$ and $\bigoplus_{i=1}^s \cO_{C^\p_i}$, and by $[w]\in\HoH_0(H^0\Ga)$ and $[w^\p]\in\HoH_0(H^0\Ga^\p)$ the canonical classes. Suppose there is a triangle equivalence
\[
?\lten_{H^0 \Ga} Z: \D(H^0\Ga)\iso \D(H^0\Ga^\p)
\] 
given by a dg bimodule $Z$ such that $\HoH_0(Z)$ (as defined
in \cite{Keller98}) takes $[w]$ to $[w^\p]$
in $\HoH_0(H^0\Ga^\p)=\HoH_0(\Ga^\p)$. Then $\hX$ is isomorphic to $\hX^\p$. In particular, $s$ is equal to $t$.
\end{maintheorem}

For general (non-simple) flopping contractions, there exist derived equivalent algebras $H^0\Ga$ and $H^0\Ga^\p$ that are non-isomorphic. August proves that the isomorphism classes of such algebras in a fixed derived equivalence class of $H^0\Ga$ are 
precisely the contraction algebras for the iterated flops of $\hY$ (see Theorem 1.4 of \cite{August18a}). Different contraction algebras in the same derived equivalence class are related by the iterated mutations of the tilting objects. The mutations are the homological counterpart of flops between different minimal models.  We refer to \cite{Wemyss18, Wemyss21} for the general framework of the homological minimal model program.

Theorem B says that the underlying singularity type of the smooth minimal models is determined by the derived equivalence class of the pair $(H^0\Ga,[w])$. We sketch the idea of the proof. From 3-dimensional birational geometry
we know that  the underlying (isolated) singularity of a smooth flopping contraction is a hypersurface 
(see Section \ref{sec:flop}). It is a classical theorem of Mather and Yau that 
up to isomorphism, a germ of isolated hypersurface singularity is determined by its Tyurina algebra (see \cite{MY} for the analytic case   and \cite{GP} for the formal case). Next we prove that the derived equivalence class of $H^0\Ga$ together with the canonical class $[w]$
recovers the Tyurina algebra of the singularity. We solve this problem in two steps. First, we prove that the Tyurina algebra, therefore the isomorphism class of the hypersurface singularity, can be recovered from the ($\ZZ$-graded dg enhanced)
cluster category $\cC_\Ga$. This result, proved in Section 5,  should have independent interest.  Secondly, we show that the isomorphism class of the Ginzburg algebra 
$\fD(Q,w)$ that is quasi-isomorphic to 
$\Ga$ can be recovered from the data $(H^0\Ga,[w])$. The proof uses a result of the first author and Gui-song Zhou in noncommutative differential calculus of potentials with finite-dimensional Jacobi-algebras \cite{HuaZhou}. Finally, we prove (in Section \ref{subsec:silting} and \ref{subsec:cyclic}) that any derived Morita equivalence 
$\D(H^0\Ga)\simeq\D(H^0\Ga^\p)$ preserving the canonical class yields a derived Morita 
equivalence $\per(\Ga)\simeq \per(\Ga^\p)$.

Note that in \cite{HT16}, the first author and Toda gave an alternative definition of the contraction algebra associated to a flopping contraction using the category of matrix factorizations. In this definition, the contraction algebra carries an additional (compared with the definition in \cite{DW13}) $\ZZ/2$-graded $A_\infty$-structure. In \cite{Hua17}, the first author proved that the Tyurina algebra 
of the singularity can be recovered from the $\ZZ/2$-graded $A_\infty$-structure. Our proof of Theorem \ref{thmB} shows that the $\ZZ/2$-graded $A_\infty$-structure on the contraction algebra can be recovered from the class $[w]$. Theorem B without
the condition on the preservation of the canonical class is precisely 
the generalization of the conjecture by Donovan and Wemyss 
stated by August in Conjecture 1.3 of \cite{August18a}. See Conjecture \ref{DWconj} for the original conjecture of Donovan and Wemyss, which is for simple flopping contractions. The generalized form of the conjecture has recently been proved
in Appendix~A of \cite{JassoMuro22},  cf. also \cite{JassoKellerMuro23}, by combining the derived Auslander--Iyama 
correspondence of [loc.~cit.] with our Theorem~\ref{HHsg=HHDsg}.
The geometric meaning of the class $[w]$ remains to be understood. It is believed that the vanishing of $[w]$ is closely related to the condition that $\hX$ is quasi-homogeneous.

Our third result is a necessary condition on the contractibility of a nc rigid rational curve in a smooth CY 3-fold. Let $u$ be a variable of degree $2$.
\begin{maintheorem}\label{thmC}(Proposition \ref{contractible=>kiu})
Let $C$ be a nc rigid rational curve in a smooth CY 3-fold $Y$. If $C$ is contractible, 
then its derived deformation algebra $\Ga$ is $\kiu$-enhanced (see definition in Section \ref{sec:contractibility}). 
Moreover, $H^0\Ga$ is a symmetric Frobenius algebra.
\end{maintheorem}
We conjecture that a nc rigid rational curve $C$ is contractible if and only if $\Ga$ is $\kiu$-enhanced (see Conjecture \ref{conj-contractibility}). 

The paper is organized as follows. In Section \ref{sec:pre} we review basics on derived noncommutative deformation theory,  noncommutative crepant resolutions and flopping contractions. Concerning derived deformations, we complement the results of
\cite{ELO2} by explaining the link to classical deformations in abelian categories
in subsection~\ref{ss:classical-defo}. In Section \ref{sec:CYcluster}, we discuss various notions of Calabi-Yau structures in geometry and algebra and prove Theorem \ref{thmA}. The notion of cluster category is introduced in Section \ref{subsec:cluster category}. When the curve is contractible, there are two cluster categories associated to it: one via derived deformation and the other via the NCCR. We prove that these two cluster categories are algebraically equivalent. 
In Section \ref{sec:Ginz}, we recall the definition of Ginzburg algebras and several results in noncommutative differential calculus including the noncommutative Mather-Yau theorem and noncommutative Saito theorem. We further show that for a contractible curve in a CY 3-fold, there exists a Ginzburg algebra weakly equivalent to the derived deformation algebra whose potential is canonically defined up to right equivalence. Then we establish a relation between the silting theory of a non positive dg algebra and the silting theory of its zeroth homology. In Section \ref{sec:sing}, we study the relation between the cluster category associated 
to the contractible curves and their underlying singularities via Hochschild cohomology. 
In particular, Theorem \ref{thmB} is proved. In Section \ref{sec:contractibility}, we introduce the notion of $\kiu$-enhancement of dg algebras.  For derived deformation algebra, we establish a link between the existence of $\kiu$-enhancement and contractibility of rational curve and prove Theorem \ref{thmC}.

\vspace{.5cm}
\paragraph{\bf{Acknowledgments.}}   
The first author wishes to thank Yujiro Kawamata, Mikhail Kapranov, Sheldon Katz and Gui-song Zhou for inspiring discussions, and Aron Heleodoro for the help on understanding \cite{Gaits11}. 
The second author thanks Akishi Ikeda for kindly explaining the details of Section~3.5 of \cite{IkedaQiu18} and Gustavo
Jasso for help with section~7.2. He is grateful to 
Zhengfang Wang for inspiring discussions on singular Hochschild cohomology and
for providing reference \cite{GGRV}. Both authors want to thank Michael Wemyss for many valuable comments and suggestions, in particular for drawing our attention to the work of August \cite{August18a,August18b}. They are greatly indebted to an anonymous referee for reading previous versions of the manuscript
with great care, pointing out numerous local errors and inaccuracies and helping to improve the readability of the paper.
The research of the first author was supported by RGC General Research Fund no. 17330316, no. 17308017.

\section{Preliminaries}\label{sec:pre}

\subsection{Notation and conventions} Throughout $k$ will be a ground field unless stated otherwise. Unadorned tensor products are over $k$. Let $V$ be a $k$-vector space. We denote its dual vector space by $DV$. When $V$ is graded, $DV$ is understood as the the dual in the category of graded vector spaces. For a subspace $V^\p$ of a complete topological vector space $V$, we denote the closure of $V^\p$ in $V$ by $(V^\p)^c$.
By definition a \emph{pseudo-compact $k$-vector space} is a linear topological vector space which is complete and whose topology is generated by subspaces of finite codimension.  Following \cite{VdB15} we will denote the corresponding category 
by $PC(k)$. We have inverse dualities
\[
\bD: Mod(k)\to PC(k)^{op} :V  \mapsto \Hom_k(V,k)
\]
\[
\bD : PC(k) \to Mod(k)^{op} : W  \mapsto \Hom_{PC(k)}(W, k)
\]
where we recall that for $V \in Mod(k)$ the topology on $\bD V$ is generated by the kernels
of $\bD V \to \bD V^\p$ where $V^\p$ runs through the finite dimensional subspaces of $V$.
Similarly, if $V$ is graded then $\bD$ is understood in the graded sense. For the definition of $\Hom$-space and tensor product in $PC(k)$, we refer to Section 4 of \cite{VdB15}. Using the tensor product in $PC(k)$, we define the pseudo-compact dg algebras, modules and bimodules to be the corresponding objects in the category of
graded objects of $PC(k)$. Let $A$ be pseudo-compact dg $k$-algebra. Denote by $PC(A^e)$ the category of pseudo-compact $A$-bimodules.
We will sometimes take a finite dimensional separable $k$-algebra $l$ to be the ground ring. The definition of the duality functor $\bD$ on $PC(l^e)$ requires some extra care due to the noncommutativity of $l$. We refer to Section 5 of \cite{VdB15} for the detailed discussion.

Denote by $\PCAlgc(l)$ the category of augmented pseudo-compact
dg algebras $A$ whose underlying graded algebras have their augmentation ideal equal
to their Jacobson radical (cf. Proposition 4.3 and section 6 of \cite{VdB15}).
Our main interest is in the case when $l\cong ke_1\times ke_2\times \ldots\times ke_n$ for central orthogonal idempotents $(e_i)_i$. For an object $A\in \PCAlgc(l)$, we use 
$\Hom$-spaces and tensor products in $PC(l^e)$ to define the Hochschild and cyclic (co)homology. 
For details, we refer to Section 7 and Appendix A of \cite{VdB15}. 
If $A$ is an $l$-algebra in $\PCAlgc(l)$,  we use $\HoH_*(A), \HoH^*(A), \HC_*(A)$ to denote the continuous Hochschild homology, cohomology and cyclic homology of $A$. Because for a pseudo-compact dg algebra, we will only consider continuous Hochschild homology, cohomology and cyclic homology, there is no risk of confusion. 
By an abuse of notations, for $A\in \PCAlgc(l)$ we denote by $\D(A)$  
the pseudo-compact derived category of $A$. Its subcategories $\per(A)$ and 
$\D_{fd}(A)$ are defined as the thick subcategory generated by the free
$A$-module $A$ and as the full subcategory of all objects with
homology of finite total dimension. Similar to the algebraic case, the notion of 
homological smoothness can be defined in the pseudo-compact setting. 
We refer to the Appendix of \cite{KY11} for a careful treatment.
For the bar-cobar formalism and Koszul duality of pseudo-compact dg algebras, we refer to Appendix A and D of \cite{VdB08}.

\subsection{Derived deformation theory}
We briefly recall the setup of derived noncommutative deformation theory of 
Efimov, Lunts and Orlov.  In this section, we fix a field $k$.
We refer to \cite{Keller06} for foundational material
on dg categories. For a dg category $\cA$, we denote by $\D(\cA)$ the
derived category of right dg $\cA$-modules. Fix a positive integer $n$ and
let $l$ be the separable $k$-algebra $ke_1\times\ldots \times ke_n$.
An \emph{$l$-algebra} $A$ is a $k$-algebra together with a morphism
of $k$-algebras $l\to A$ (note that $l$ is {\em not} necessarily central in $A$).
An equivalent datum is that of the $k$-category with $n$ objects $1$, \ldots, $n$, 
whose morphism space from $i$ to $j$ is given by $e_j A e_i$.
An \emph{$l$-augmented} (dg) algebra is a (dg) $l$-algebra  $\cR$ together with an 
$l$-algebra morphism $\cR\to l$ such that the composition $l\to \cR\to l$ is the identity morphism. Its \emph{augmentation ideal} is the kernel of the augmentation
morphism $\cR \to l$. 
An \emph{artinian} $l$-algebra is an augmented $l$-algebra
whose augmentation ideal is finite-dimensional and nilpotent.
A dg $l$-algebra is \emph{artinian} if it is an augmented dg $l$-algebra whose
augmentation ideal is finite-dimensional and nilpotent.
Denote by $\Art_l$ and $\cArt_l$ the categories of artinian $l$-algebras and of
commutative artinian $l$-algebras. Denote by $\dgArt_l$ the category of artinian dg algebras and by $\dgArt_l^-$ the subcategory of 
$\dgArt_l$ consisting of dg algebras concentrated in nonpositive degrees.

Fix a dg category $\cA$ and a dg $\cA$-module $E$ with a decomposition
$E=E_1\oplus \cdots \oplus E_n$. We view $E$ as an $l^{op}\ten\cA$-module
in the natural way. The {\em dg endomorphism $l$-algebra} of $E$ is
the dg endomorphism algebra over $\cA$ of the sum $E$ viewed as
an $l$-algebra in the natural way.
We are going to define a pseudo-functor $\Def(E)$ from  
$\dgArt_l$ to the category $\Gpd$ of groupoids. This pseudo-functor assigns to an artinian dg $l$-algebra $\cR$ the groupoid $\Def_\cR(E)$ of 
$\cR$-deformations of $E$ in the derived category $\D(\cA)$. We will mostly
follow the notations of \cite{ELO1} and identify $\cR$ with the dg category with 
$n$ objects $1$, \ldots, $n$, where the morphism complex from $i$ to $j$ is given
by $e_j \cR e_i$. Denote the dg category $\cR^{op}\ten\cA$ by $\cA_\cR$. 
The augmentation $\eps:\cR\to l$ yields the
functor of extension of scalars $\eps^*$ taking a dg $\cR$-module $S$ to 
the dg $l^{op}\ten\cA$-module
\[
\eps^*(S)=l\lten_\cR S.
\]

\begin{definition}\cite[Definition 10.1]{ELO1} \label{def-Def}
Fix an artinian augmented
dg $l$-algebra $\cR$. An object of the groupoid $\Def_\cR(E)$ is a pair 
$(S, \sigma)$, where $S$ is an object of $\D(\cA_\cR)$ and
\[
\sigma : \eps^*(S) \to E
\]
is an isomorphism in $\D(l^{op}\ten\cA)$.
A morphism $f : (S,\sigma) \to (T,\tau)$ between two $\cR$-deformations of $E$ is an isomorphism $f: S\to T$  in $\D(\cA_\cR)$ such that
\[
\tau \circ \eps^*(f)=\sigma.
\]
This defines the groupoid $\Def_\cR(E)$. A homomorphism of augmented artinian dg 
$l$-algebras $\phi : \cR \to \cQ$ induces the functor
\[
\phi^* : \Def_\cR(E) \to \Def_\cQ(E)
\]
given by $\cQ\lten_\cR ?$. Thus we obtain a pseudo-functor
\[
\Def(E) : \dgArt_l \to \Gpd.
\]
We call $\Def(E)$ the {\em pseudo-functor of derived deformations of $E$}.  
We denote by $\Def_-(E)$  the restrictions of the
pseudo-functor $\Def(E)$ to the subcategory $\dgArt_l^-$.
\end{definition}

The category of augmented dg $l$-algebras can be naturally enhanced to a weak 
2-category. We refer to Definition 11.1 of \cite{ELO2} for the precise definition of the 
2-category structure.  In particular, we denote the corresponding 2-categorical enhancements of $\dgArt_l$, $\dgArt_l^-$ and $\Art_l$ by $2$-$\dgArt_l$, 
$2$-$\dgArt_l^-$ and $2$-$\Art_l$ (in \cite{ELO2}, they are denoted by
$\tp$-$\dgArt_l$ etc.).
By Proposition 11.4 of \cite{ELO2}, there exists a pseudo-functor $\DEF(E)$ from 
$2$-$\dgArt_l$ to $\Gpd$ and which is an extension to $2$-$\dgArt_l$ of the 
pseudo-functor $\Def(E)$. Similarly, there exists a pseudo-functor
$\DEF_-(E)$  extending $\Def_-(E)$.

The main theorem of \cite{ELO2} is:
\begin{theorem} \cite[Theorem 15.1, 15.2]{ELO2}\label{rep-ELO}
Let $E_1$, \ldots, $E_n$ be a collection of objects in $\D(\cA)$.
Let $E$ be the direct sum of the $E_i$ and $C$ the extension
algebra $\Ext^*_\cA(E,E)$ considered as a graded $l$-algebra.
Assume that we have 
\begin{enumerate}
\item[$(a)$] $C^p=0$ for all $p<0$;
\item[$(b)$] $C^0=l$;
\item[$(c)$] $\dim_k C^p<\infty$ for all $p$ and $\dim_k C^p=0$ 
for all $p\gg 0$.
\end{enumerate}
Denote by $\cC$ the dg endomorphism $l$-algebra of $E$.
Let $A$ be a strictly unital minimal model of $\cC$. 
Then the pseudo-functor $\DEF_-(E)$ is pro-representable by the dg $l$-algebra 
$\Ga = \bD B{A}$, where $B$ denotes the bar construction. 
That is, there exists an equivalence of pseudo-functors 
$\DEF_-(E) \simeq h_{\Ga}$ from $2$-$\dgArt_l^-$ to $\Gpd$, 
where $h_{\Ga}$ denotes the groupoid of $1$-morphisms $\mbox{1-}\Hom(\Ga, ?)$.
\end{theorem}

In the case where the dg category $\cA$ is given by an algebra $A$ concentrated
in degree $0$ and $E$ is a one-dimensional $A$-module, Booth 
\cite[Theorem 3.5.9]{Booth18} obtains an analogous prorepresentability result for
the set-valued framed deformation functor $\mathrm{Def}^{\mathrm{fr},\leq 0}_A(E)$
without having to impose the finiteness condition (c).

Let $Y$ be a smooth algebraic variety. A collection of compactly supported coherent sheaves  
$\cE_1, \ldots, \cE_t$ on $Y$ is called \emph{semi-simple} if 
$\Hom_Y(\cE_i,\cE_i)\cong k$ for all $i$ and $\Hom_Y(\cE_i,\cE_j)=0$ for
all $i\neq j$. 
The finiteness assumption in Theorem \ref{rep-ELO} is satisfied by any 
semi-simple collection. 
Let $\cE$ be the direct sum of such a collection of coherent sheaves on $Y$. 
We may denote the completion of $Y$  along the support of $\cE$ by $\hY$. 
\begin{corollary}\label{localgenerator}
Given a semi-simple collection $\{\cE_i\}_{i=1}^t$ in $\D^b(\coh \hY)$, denote by $\cC$ the dg endomorphism $l$-algebra of  $\cE:=\bigoplus_{i=1}^t\cE_i$.
Let $A$ be a strictly unital minimal model of $\cC$ and $\Ga= \bD B{A}$.  Then there
is an equivalence
\[
\DEF_-(\cE)\simeq  h_\Ga.
\]
\end{corollary}
We call $\Ga$ the \emph{derived deformation algebra} of the collection $\{\cE_i\}_{i=1}^t$ in $Y$.
When we want to emphasize the dependence on $Y$ and $\{\cE_i\}_{i=1}^t$, $\Ga$ is replaced by $\Gamma^Y_{\cE}$.  The semi-simple collection $\{\cE_i\}_{i=1}^t$ is called \emph{nc rigid} if 
$\dim_kH^0(\Ga^Y_\cE)<\infty$.
In this paper, we are mainly interested in the case when $\{\cE_i\}_{i=1}^t$ is (the structure sheaves) a collection of smooth rational curves $C_1,\ldots, C_t$ that satisfies the condition that $\Hom_Y(\cO_{C_i},\cO_{C_j})=0$ for $i\neq j$. For  such a collection of rational curves $C:=\{C_i\}_{i=1}^n$ where we write $\Ga_{C}^Y$ for $\Ga^Y_{\{\cO_{C_i}\}_{i=1}^t}$.

In the context of classical noncommutative deformation theory, the representability of noncommutative deformations of contractible rational curves was proved by Donovan and Wemyss.
\begin{theorem}(Proposition 3.1, Corollary 3.3 \cite{DW13})
Let $f:Y\to X$  be a simple flopping contraction of 3-folds (see definition in Section \ref{sec:flop}) and let $C$ be the reduced exceptional fiber of $f$. 
The functor
\[
\pi_0(\Def^{cl}(\cO_C)): \Art_k\to \Set
\]
is representable. The artinian algebra $\Lm$ representing it is called
the {\em contraction algebra} associated to $f: Y\to X$.
\end{theorem}
The definition of the classical deformation functor $\Def^{cl}$ is recalled in the next
section (cf. Section 2 of \cite{DW13}). If $\Ga$ is the derived deformation
algebra of $C$ (with $t=1$), it follows from the above Theorem and 
Theorem~\ref{thm:classical-defo} below that the contraction algebra
$\Lm$ is isomorphic to $H^0\Ga$. Indeed, they both represent the
same deformation functor and this determines them up to
(non unique) isomorphism (cf. the proof of Theorem~2.14 in
\cite{Segal08}).

\subsection{Link to classical deformations} 
\label{ss:classical-defo} Let $\cA$ be a dg category.
Let $\cH\subset\D(\cA)$ be the heart of a $t$-structure on $\D(\cA)$.
We assume that $\cH$ is {\em faithful}, i.e. the higher extension
groups computed in $\cH$ are canonically isomorphic to those
computed in $\D(\cA)$.

Let $\cR$ be an augmented artinian $l$-algebra. By an {\em $\cR$-module in
$\cH$}, we mean an object $M$ of $\cH$ endowed with an algebra
homomorphism $\cR \to \End(M)$. Given such an $\cR$-module,
we denote by $?\ten_\cR M$ the unique right exact functor
$\mod \cR \to \cH$ extending the obvious additive functor $\proj ~\cR \to \cH$
taking $\cR$ to $M$. Here we denote by $\proj~ \cR$ the category
of finitely generated projective (right) $\cR$-modules and by $\mod\cR$ the
category of finitely generated $\cR$-modules. Notice that $?\ten_\cR$ is
a left Kan extension and thus unique up to unique isomorphism and
functorial in $M$. It can be computed using projective resolutions. It is obvious how
to define morphisms of $\cR$-modules in $\cH$.

Let $E$ be the direct sum of a collection of $n$ objects $E_1$, \ldots, $E_n$  of $\cH$.
We view $E$ as an $l$-module in $\cH$ in the natural way.
For an augmented artinian $l$-algebra $\cR$, we define the groupoid
$\Def^{cl}_\cR(E)$ of classical deformations of $E$ as follows:
Its objects are pairs $(M,\mu)$ where $M$ is an $\cR$-module in
$\cH$ such that the functor $?\ten_\cR M$ is exact and $\mu: l\ten_\cR M \iso E$
is an isomorphism of $l$-modules in $\cH$. A morphism $(L,\lambda) \to (M,\mu)$
is an isomorphism $f: L \to M$ of $\cR$-modules in $\cH$ such that
$\mu \circ (l\ten_\cR f) =\lambda$. 

For an augmented $l$-algebra $A$ and an augmented artinian
$l$-algebra $\cR$, we define $G(A, \cR)$ to be the groupoid
whose objects are the morphisms $A \to \cR$ of augmented
$l$-algebras and whose morphisms $\phi_1 \to \phi_2$ are the invertible
elements $r$ of $\cR$ such that $\phi_2(a) = r \phi_1(a) r^{-1}$ for
all $a$ in $A$.

\begin{theorem} \label{thm:classical-defo}
Suppose that in addition to the above assumptions,
$E$ satisfies the hypotheses of Theorem 
\ref{rep-ELO}.  
Let $\Ga$ be the pseudo-compact dg $l$-algebra defined there.
Let $\cR$ be an augmented artinian $l$-algebra.
Then $H^0\Gamma$ represents the classical deformations of $E$ in the
sense that there is an equivalence of groupoids
\[
\Def^{cl}_\cR(E) \iso G(H^0\Gamma, \cR).
\]
\end{theorem} 

\begin{proof} By Theorem \ref{rep-ELO}, we have an equivalence of groupoids
\[
\DEF_-(E)(\cR) \iso \mbox{1-}\Hom(\Ga, \cR),
\]
where $\mbox{1-}\Hom$ denotes the groupoid of $1$-morphisms in $2$-$\dgArt$.
We will show that $\DEF_-(E)(\cR)$ is equivalent to $\Def^{cl}_\cR(E)$ and
$\mbox{1-}\Hom(\Ga, \cR)$ is equivalent to $G(H^0\Gamma,\cR)$. We start with the
second equivalence. By Definition~11.1 of \cite{ELO2}, an object
of $\mbox{1-}\Hom(\Ga,\cR)$ is a pair $(M,\theta)$ consisting of
\begin{itemize}
\item a dg bimodule $M$ in $\D(\Ga^{op}\ten \cR)$
such that the restriction to $\cR$ of $M$ is isomorphic to $\cR$ in $\D(\cR)$ and
\item an isomorphism $\theta : M\lten_\cR l \to l$ in $\D(\Ga^{op})$.
\end{itemize}
A $2$-morphism $f:(M_1,\theta_1) \to (M_2,\theta_2)$ is an isomorphism
$f:M_1\to M_2$ in $\D(\Ga^{op}\ten\cR)$ such that 
$\theta_2\circ (f\lten_\cR l) =\theta_1$. We define a functor
$F: G(H^0\Ga,\cR) \to \mbox{1-}\Hom(\Ga,\cR)$ as follows: Let $\phi: H^0(\Ga) \to\cR$
be a morphism of augmented $l$-algebras. Since $\Ga$ is concentrated in
degrees $\leq 0$, we have a canonical algebra morphism $\Ga \to H^0(\Ga)$.
By composing it with $\phi$ we get a morphism of augmented dg $l$-algebras
$\Ga \to \cR$. It defines a structure of dg bimodule $M$ on $\cR$. We put
$F\phi=(M,\theta)$, where $\theta: \cR\lten_\cR l \to l$ is the canonical
isomorphism.   Now let $\phi_1$ and $\phi_2$ be two morphisms of augmented
algebras $H^0(\Ga)\to \cR$. Put $(M_i,\theta_i)=F\phi_i$, $i=1,2$.  Let $r$ an 
be an invertible element of $\cR$ such that
$\phi_2(a)=r\phi_1(a)r^{-1}$ for all $a$ in $H^0(\Ga)$. Then it is clear
that the left multiplication with $r$ defines an isomorphism of bimodules
$M_1 \to M_2$ compatible with the $\theta_i$. Recall that $\Ga^{op}\ten\cR$
is concentrated in degrees $\leq 0$ so that its derived category has a 
canonical $t$-structure.
Since the $M_i$ live
in the heart of this $t$-structure on $\D(\Ga^{op}\ten\cR)$, it
is also clear that $F$ is fully faithful.  It remains to be checked that $F$
is essentially surjective. So let $(M,\theta)$ be given. Since $M$ is
quasi-isomorphic to $\cR$ when restricted to $\cR$, its homology is
concentrated in degree $0$. We can therefore replace $M$ with
$H^0(M)$, which is an ordinary 
$H^0(\Ga)$-$\cR$-bimodule isomorphic to $\cR$ as a right $\cR$-module
(we also consider it as a left $\Ga$-module via the canonical morphism
$\Ga\to H^0(\Ga)$).
In particular, $M$ is right projective and so $M\lten_\cR  ?= M\ten_\cR ?$.
We choose an isomorphism $f: M \iso \cR$ of right $\cR$-modules.
After multiplying $f$ with an invertible element of $l$, we may assume that
$f\ten_\cR l=\theta$. The left $\Ga$-module structure on $M$ yields
an algebra morphism 
\[
\phi : H^0(\Gamma) \to \End_\cR(M) \iso \End_\cR(\cR) =\cR.
\]
It is clear that $f$ yields an isomorphism between $(M,\theta)$ and $F\phi$.

We now construct an equivalence from $\DEF_-(E)(\cR)$ to  $\Def^{cl}_\cR(E)$.
Recall from Proposition~11.4 of \cite{ELO2} that the groupoid 
$\DEF_-(E)(\cR)$ equals the groupoid $\Def_\cR(E)$ of
Definition~\ref{def-Def} (but $\DEF_-$ has enhanced $2$-functoriality).
Let $P \to E$ be a cofibrant resolution of $E$. Since the graded
algebra $\Ext^*(E,E)$ has vanishing components in degree $-1$ and
in all sufficiently high degrees, we can apply Theorem~11.8 of \cite{ELO2}
to conclude that the groupoid $\Def_\cR(E)$ is equivalent to the groupoid 
$\Def_\cR^h(P)$ of homotopy deformations of Definition~4.1 of \cite{ELO1}. 
We now construct
an equivalence $F$ from $\Def_\cR^h(P)$ to $\Def_\cR^{cl}(E)$.
Let $(S,\sigma)$ be an object of $\Def_\cR^h(P)$. We may
assume that $S=\cR\ten_l P$ as a graded bimodule and that $\sigma$
is the canonical isomorphism $l\ten_\cR(\cR\ten P) \iso P$.
Let $I$ denote the augmentation ideal of $\cR$. Then $S$ has
a finite filtration by the dg submodules $I^p S$, $p\geq 0$, and
each subquotient is isomorphic to a summand of a finite
sum of copies of $l\ten_\cR S = P$.
Thus, the underlying dg $\cA$-module $M$ of $S$ is isomorphic in $\D(\cA)$
to a finite iterated extension of objects of $\add(E)$, the subcategory of direct
factors of finite direct sums of copies of $E$. Therefore, $M$ still lies in
the heart $\cH$. Note that as shown in the
proof of Theorem~11.8 of \cite{ELO2}, $S$ is cofibrant over $\cR^{op}\ten \cA$.
Therefore, $M$ is cofibrant over $\cA$. The left $\cR$-module structure
on $S$ yields an algebra homomorphism $\cR \to \End(M)$. Since each
object of $\mod \cR$ is a finite iterated extension of one-dimensional $l$-modules,
the functor $?\ten_\cR S :\D(\cR) \to \D(\cA)$ takes $\mod \cR$ to $\cH$.
Since $?\ten_\cR S$ is a triangle functor, the induced functor $\mod\cR \to \cH$
is exact. Clearly it restricts to the natural functor $\proj ~\cR \to \cH$ and
is therefore isomorphic to $?\ten_\cR M: \mod \cR \to \cH$.
Finally, the isomorphism $l\ten_\cR S \iso E$ yields an
isomorphism $\mu: l\ten_\cR M \iso E$. In this way, to an object $(S,\sigma)$
of $\Def_\cR^h(P)$, we have associated an object $F(S,\sigma)=(M,\mu)$
of $\Def^{cl}_\cR(E)$. Notice that by what we have just shown, 
we may also describe $M\in\cH$ as the zeroth homology $H^0_\cH(S)$ with respect to 
the $t$-structure associated with $\cH$, that we have an isomorphism
$l \ten_\cR M \iso H^0_\cH(l \lten_\cR S)$ and that the isomorphism 
$\mu: l \ten_\cR  M \iso E$ is induced by $\sigma: l \lten_\cR S \iso E$.
Recall that a morphism $(S_1, \sigma_1) \to (S_2,\sigma_2)$
of $\Def^h_\cR(P)$ is a class of isomorphisms 
$S_1 \to S_2$ of dg $\cR^{op}\ten\cA$-modules
compatible with the $\sigma_i$ modulo homotopies compatible with the $\sigma_i$.
Since the functor $\Def_\cR^h(P) \to \Def_\cR(E)$ is an 
equivalence, these
morphisms are in bijection with the isomorphisms $S_1 \to S_2$ of
$\D(\cR^{op}\ten\cA)$ compatible with the $\sigma_i$. Clearly each
such morphism induces an  isomorphism $(M_1, \lambda_1) \to (M_2,\lambda_2)$, where $(M_i,\lambda_i)=F(S_i,\sigma_i)$, $i=1,2$. It follows from 
Lemma~\ref{lemma:Toda-condition} below 
that this assignment is a bijection. It remains to
be shown that $F: \Def^h_{\cR}(P)\to \Def^{cl}_\cR(E)$ is essentially surjective.
Since we have an equivalence $\Def^h_\cR(P) \iso \Def_\cR(E)$, it
suffices to lift a given object $(M,\mu)$ of $\Def_\cR^{cl}(E)$
to an object $(S,\sigma)$ of $\Def_\cR(E)$. Let $A$ denote the
dg endomorphism $l$-algebra $\RHom_\cA(M,M)$. Then $M$ becomes
canonically an object of $\D(A^{op}\ten\cA)$. Now since $M$ is in the heart
of a $t$-structure, its negative self-extension groups vanish and we have
a quasi-isomorphism $\tau_{\leq 0}A \iso \End_\cH(M)$. Thus, in the
homotopy category of dg $l$-algebras, we have a morphism
\[
\cR \to \End_\cH(M) \iso \tau_{\leq 0} A \to A.
\]
By tensoring with $\cA$ we obtain a morphism $\cR^{op}\ten\cA \to A^{op}\ten \cA$
in the homotopy category of dg categories. The associated restriction
functor $\D(A^{op}\ten\cA) \to \D(\cR^{op}\ten\cA)$ sends $M$ to an object
$S$ of $\D(\cR^{op}\ten\cA)$. By construction, the restriction of $S$ to
$\cA$ is isomorphic to $M$ in $\D(\cA)$ and the left action of $\cR$
on $S$ induces the given algebra morphism $\cR \to \End_\cH(M)$.
Since $\cH$ is the heart of a $t$-structure, we have a canonical
realization functor $\D^b(\cH) \to \D(\cA)$ extending the inclusion
$\cH \to \D(\cA)$, cf. Section~3.1.10 of \cite{BeilinsonBernsteinDeligne82}
or Section~3.2 of \cite{KellerVossieck87}. Moreover, since
$\cH$ is faithful, the realization functor is fully faithful. Since we only know
how to compare tensor functors, we use a different construction
to extend the inclusion $\cH \to \D(\cA)$ to a triangle functor $\D^b(\cH) \to \D(\cA)$.
Let $\cH_{dg}$ be the full subcategory of the dg category of right $\cA$-modules
formed by cofibrant resolutions of the objects of $\cH$. We have an
equivalence of $k$-categories $\cH \iso H^0(\cH_{dg})$. Since $\cH$ is
the heart of a $t$-structure, the homology of the dg category $\cH_{dg}$
is concentrated in degrees $\geq 0$. Thus, we have quasi-equivalences
$\tau_{\leq 0} \cH_{dg} \iso H^0(\cH_{dg}) \iso \cH$. Therefore, in the
homotopy category of dg categories, we obtain a morphism
\[
\cH \to \tau_{\leq 0}(\cH_{dg}) \to \cH_{dg} \to \D_{dg}(\cA)
\]
where $\D_{dg}(\cA)$ denotes the dg category of cofibrant dg $\cA$-modules.
It gives rise to an $\cH$-$\cA$-bimodule $R$. Let $\cH^b(\cH)$ denote
the category modulo homotopy of bounded complexes of objects of $\cH$.
Using the fact that short
exact sequences of $\cH$ give rise to triangles in $\D(\cA)$, one checks
that the induced functor $?\ten_\cH R: \cH^b(\cH) \to \D(\cA)$ vanishes
on the bounded acyclic complexes and therefore induces a triangle
functor $\D^b(\cH) \to \D(\cA)$ still denoted by $?\ten_\cH R$. Since
$\cH$ is a faithful heart, one obtains that this triangle functor is
fully faithful. We claim that we have
a square of triangle functors commutative up to isomorphism
\[
\xymatrix{
\D^b(\mod\cR) \ar[rr]^-{?\ten_\cR M} \ar[d] & &\D^b(\cH) \ar[d]^{?\ten_\cH R} \\
\D(\cR) \ar[rr]_{?\lten_\cA S} & &\D(\cA).}
\]
To check this, one has to check that the bimodules $S$ and $M\ten_\cH R$ 
are isomorphic in $\D(\cR^{op}\ten\cA)$. This is easy using 
Lemma~\ref{lemma:Toda-condition}
below. Since $?\ten_\cR M: \mod \cR \to \cH$ is exact, the
given isomorphism $l\ten_\cR M \iso E$ yields an isomorphism
$l\lten_\cR M \iso E$ in $\D^b(\cH)$ and thus an isomorphism
$\sigma : l\lten_\cR S \iso E$ in $\D(\cA)$. It is now clear
that we can recover $M$ as $H^0_\cR(S)$ and that
$\mu: l \ten_\cR M \iso E$ is the morphism induced by 
$\sigma : l\lten_\cR S \iso E$ in $H^0_\cR$. By the description
of $F$ given above, this shows that we do have $F(S,\sigma) \iso (M,\mu)$.
\end{proof}

Let $\cA$ be a dg category and $B$ an ordinary $k$-algebra. Let
$X$ and $Y$ be objects of $\D(B^{op}\ten\cA)$ and let $\res(X)$
be the restriction of $X$ to $\cA$. The left action of $B$ on
$X$ defines an algebra morphism
\[
\alpha_X: B \to \End_{\D(\cA)}(\res(X)).
\]
Let $\cM(X,Y)$ be the space of all morphisms $f:\res(X) \to \res(Y)$ in
$\D(\cA)$ such that
\[
f\circ \alpha_X(b) = \alpha_Y(b)\circ f
\]
for all $b\in B$. The restriction functor induces a natural map
\[
\Phi: \Hom_{\D(B^{op}\ten\cA)}(X,Y) \to \cM(X,Y).
\]
\begin{lemma} \label{lemma:Toda-condition} If we have
\[
\Hom_{\D(\cA)}(\res(X), \Sigma^{-n}\res(Y))=0
\]
for all $n>0$, the map $\Phi$ is bijective.
\end{lemma}
\begin{proof} We adapt the argument of Section~5 of \cite{Keller00}.
We may suppose that $X$ is cofibrant over $B^{op}\ten\cA$
and in particular cofibrant over $\cA$. Then the sum total dg module
of the bar resolution
\[
\xymatrix{
\ldots \ar[r] & B \ten B^{\ten p} \ten X \ar[r] & \ldots \ar[r] & B\ten B \ten X \ar[r] & 
B\ten X \ar[r] & 0}
\]
is still cofibrant over $B^{op}\ten\cA$ and quasi-isomorpic to $X$. We use it
to compute $\Hom_{\D(B^{op}\ten \cA)}(X,Y)$. By applying $\Hom_{B^{op}\ten \cA}(?,Y)$ to the bar resolution, we get a double complex $D$ of the form
\[
\xymatrix{
\Hom_{\cA}(X,Y) \ar[r] & \Hom_\cA(B\ten X, Y) \ar[r] & \ldots \ar[r] &
\Hom_\cA(B^{\ten p}\ten X,Y) \ar[r] & \ldots}
\]
We have to compute $H^0$ of the product total complex $\Tot^\Pi D$. 
Let $D_{\geq 0}$ be the double complex obtained by applying the
intelligent truncation functor $\tau_{\geq 0}$ to each column of $D$.
Let $D_{<0}$ be the kernel of $D \to D_{\geq 0}$. We claim that the
product total complex of $D_{<0}$ is acyclic. Indeed, the homology
of the $p$th column of $D_{<0}$ in degree $-q$ is isomorphic to
\[
\Hom_{\D(\cA)}(B^{\ten p}\ten X, \Sigma^{-q} Y).
\]
It vanishes for $-q<0$ by our assumption. To show that the product total
complex of $D_{<0}$ is acyclic, we consider the column filtration 
$F_p D_{<0}$. Then $D_{<0}$ is the inverse limit of the $F_p D_{<0}$.
By induction on $p$, each $F_p D_{<0}$ has an acyclic total complex.
Moreover, the transition maps $F_{p+1} D_{<0} \to F_p D_{<0}$ induce
componentwise surjections in the total complexes. It follows that
the inverse limit of the total complexes of the $F_p D_{<0}$ is
still acyclic and this inverse limit is the product total complex of $D_{<0}$.
So it is enough to compute $H^0$ of $\Tot^{\Pi} D_{\geq 0}$. 
For this, let us denote by 
\[
\ol{\Hom}^0_\cA(X,Y)
\]
the quotient of $\Hom^0_\cA(X,Y)$ by the nullhomotopic morphisms
and similarly for $\ol{\Hom}^0_\cA(B\ten X,Y)$. The space we have
to compute is the homology in degree $0$ of the total complex
of the double complex
\[
\xymatrix{ 
\Hom^1_\cA(X,Y) \ar[r] & \Hom^1_\cA(B\ten X,Y) \\
\ol{\Hom}^0_\cA(X,Y) \ar[u]\ar[r] & \ol{\Hom}^0_\cA(B\ten X,Y) \ar[u]
} \ko
\]
where the lower left corner is in bidegree $(0,0)$. This equals
the intersection of the kernels of the maps from bidegree $(0,0)$ to bidegrees
$(0,1)$ and $(1,0)$. The kernels of the vertical maps
are respectively $\Hom_{\D(\cA)}(X,Y)$ and $\Hom_{\D(\cA)}(B\ten X, Y)$
and the space we have to compute thus identifies with the
kernel of the map
\[
\Hom_{\D(\cA)}(X,Y) \to \Hom_{\D(\cA)}(B\ten X, Y)
\]
which is easily seen to take $f: X \to Y$ to 
$f \circ \alpha_X - \alpha_Y\circ f$. Thus, the homology to
be computed is isomorphic to $\cM(X,Y)$.
\end{proof}

\subsection{Noncommutative crepant resolutions}
\label{ss:NCCR}

\begin{definition}
Let $(R,\fm)$ be a complete commutative Noetherian local Gorenstein $k$-algebra of Krull dimension $n$ with isolated singularity and with residue field $k$. Denote the category of maximal Cohen--Macaulay (MCM) modules by $\CM_R$ and its stable category by $\ul{\CM}_R$. Let $N_0=R,N_1,N_2,\ldots,N_t$ be pairwise non-isomorphic indecomposables in $\CM_R$ and $A:=\End_R(\bigoplus_{i=0}^t N_i)$. We call $A$ a \emph{noncommutative resolution} (NCR) of $R$ if it has finite global dimension. A NCR is called a \emph{noncommutative crepant resolution} (NCCR) if $A$ further satisfies that
\begin{enumerate}
\item[$(a)$] $A\in \CM_R$;
\item[$(b)$] ${\rm{gldim}}(A)=n$.
\end{enumerate}
\end{definition}
If $A$ is a NCCR, we call $\op_{i=0}^t N_i$ a \emph{tilting module}.
Under the above conditions, Iyama shows that $\op_{i=1}^t N_i$ is a \emph{cluster tilting object} (see definition in Section \ref{subsec:cluster category}) in $\ul{\CM}_R$. Denote $\ol{l}$ for $A/\rad A $ and $e_0$ for the idempotent given by the projection $R\op\bigoplus_{i=1}^t N_i\to R$. Let $S_0,S_1,\ldots,S_t$ be the simple $A$-modules  with $S_0$ corresponding to the summand $R$ of $R\op\bigoplus_{i=1}^t N_i$. De Thanhoffer de V\"olcsey and Van den Bergh prove that $\ul{\CM}_R$ admits an explicit 
dg model in this case.
\begin{theorem}(\cite[Theorem 1.1]{ThV})\label{VdBTh}
There exists a finite dimensional graded $\ol{l}$-bimodule $V$ and a minimal model 
$(\wh{T}_{\ol{l}}V,d)\iso A$ for $A$, where $\wh{T}_{\ol{l}}V$ is the graded completion of the ordinary tensor algebra $T_{\ol{l}} V$ with respect to the two sided ideal generated by $V$. Put $\Gamma=\wh{T}_{\ol{l}}V/\wh{T}_{\ol{l}}V e_0 \wh{T}_{\ol{l}} V$. 
Then one has 
\[
\ul{\CM}_R\cong \per(\Ga)/\thick(S_1,\ldots, S_t)
\] and furthermore $\Gamma$ has the following properties
\begin{enumerate}
\item[$(1)$] $\Gamma$ has finite dimensional cohomology in each degree;
\item[$(2)$] $\Gamma$ is concentrated in negative degrees;
\item[$(3)$] $H^0\Gamma=A/Ae_0A$;
\item[$(4)$] As a graded algebra $\Gamma$ is of the form  $\wh{T}_{l} V^0$ for $V^0=(1-e_0)V(1-e_0)$ with $l:=\ol{l}/ke_0$.
\end{enumerate}
\end{theorem}

\subsection{Flopping contraction}\label{sec:flop}
\begin{definition}\label{def:contractible}
A smooth rational curve $C$ in a normal variety $Y$ is called \emph{contractible} if there exists an open subscheme $Y^\circ\subset Y$ containing $C$ and a proper birational morphism $f^\circ: Y^\circ\to X^\circ$ such that 
\begin{enumerate}
\item[$(1)$] $X^\circ$ is normal,
\item[$(2)$] the exceptional locus $\Ex(f^\circ)$  contains $C$,
\item[$(3)$] $f^\circ$ is an isomorphism in codimension one.
\end{enumerate}
\end{definition}
The above definition of contractibility is more restrictive than the standard one since it rules out the divisorial contraction (by the last condition). If $Y$ is a 3-fold (which is our main interest), then $\Ex(f)$ must have dimension one by condition (3). However, it may contain other components besides $C$.
If $C$ is a contractible curve in $Y$, denote by $\hX$ the formal completion of $X^\circ$ along the exceptional subscheme, i.e. where $f^\circ$ is not an isomorphism. Consider the Cartesian diagram
\[
\xymatrix{
\hY\ar[d]^\hf\ar[r] & Y^\circ\ar[d]^{f^\circ}\\
\hX\ar[r] & X^\circ
}
\] where $\hY$ is  the fiber product. We call $\hf:\hY\to \hX$ the \emph{formal contraction} associate to the contraction $f^\circ: Y^\circ\to X^\circ$.

The following definition is a special case of Definition 6.10 of \cite{KoMo}.
\begin{definition}
Let $Y$ be a normal variety of dimension 3. A \emph{flopping contraction} is a proper birational morphism $f:Y\to X$ to a normal variety $Y$ such that $f$ is an isomorphism in codimension one, and $K_Y$ is $f$-trivial. If $Y$ is smooth, then we call  $f$ a smooth flopping contraction.
\end{definition}

In this paper, we only consider smooth flopping contractions unless  stated otherwise.
Given a 3-dimensional flopping contraction $f:Y\to X$, 
let $D$ be a divisor on $Y$ such that $-(K_Y+D)$ is $f$-ample. By Theorem 6.14 of \cite{KoMo}, there exists a $D$-flop $f^+:Y^+\to X$. To be more precise, $f^+$ is a proper birational morphism that is an isomorphism in codimension one,  and $K_{Y^+}+D^+$ is $f^+$-ample where $D^+$ is the birational transform of $D$ on $Y^+$. In particular, $X$ is Gorenstein terminal. Without loss of generality, we may work locally near the exceptional fiber of $f$. By the classification theorem of 3-dimensional terminal singularities, $X$ has an isolated cDV singularity (see Corollary 5.38 of \cite{KoMo}). Recall that a 3-fold singularity $(X,0)$ is called cDV if a generic hypersurface section $0\in H\subset X$ is a Du Val singularity. Because $H$ has embedded dimension 3, $X$ has embedded dimension 4, i.e. $X$ is a hypersurface.

Denote by $\Ex(f)$ the reduced exceptional fiber of $f$. It is well known that $\Ex(f)$ 
is a tree of rational curves
\[
\Ex(f)=\bigcup_{i=1}^t C_i
\] with normal crossings such that $C_i\cong \PP^1$ (c.f. Lemma 3.4.1 of \cite{VdB04}).
We call a 3-dimensional flopping contraction $f:Y\to X$ \emph{simple} if $\Ex(f)\cong \PP^1$. Let $p$ be the singular point of $X$. By the remark above, $R:=\wh{\cO}_{X,p}$ is a complete local ring of the form $k[[x,y,u,v]]/(g)$. Now we set $\hX:=\Spec R$ and $\hY:=\hX\times_X Y$ and denote by $\hf$ the base change of $f$. We call the triple $(\hY,\wh{f},R)$ the \emph{formal contraction} associated to the flopping contraction $f:Y\to X$. Note that $\hY$ is Calabi-Yau.

Let $(\hY,\hf,R)$ be a formal flopping contraction. 
Now we consider the NCCR associated to a three dimensional 
flopping contraction $\wh{f}: \hY\to \hX$ constructed as follows.
For $i=1,\ldots,t$, let $\wh{\cL}_i$ be a line bundle on $\hY$ such that $\deg_{C_j}\wh{\cL}_i=\delta_{ij}$. Define $\wh{\cN}_i$ to be given by the maximal extension 
\begin{equation}\label{extN}
\xymatrix{
0\ar[r] & \wh{\cL}_i^{-1}\ar[r] & \wh{\cN}_i\ar[r] &\cO_{\hY}^{\oplus r_i}\ar[r] &0
}
\end{equation} associated to a minimal set of $r_i$ generators of $H^1(\hY,\wh{\cL}_i^{-1})$. 
Set $N_i:=\bR f_*\wh{\cN}_i=f_*\wh{\cN}_i$ for $i=1,\ldots,t$. 
We set
\begin{align}\label{tilting algebra}
A:=\End_{\wh{Y}}(\cO_{\wh{Y}}\op\wh{\cN}_1\op \ldots \wh{\cN}_t ) \cong \End_R(R \oplus N_1\op\ldots\op N_t ).
\end{align}
The second isomorphism can be proved as follows. It is clear that the natural morphism from $A$ to $\End_R(R \oplus N_1\op\ldots\op N_t )$ is an isomorphism away from the exceptional locus, which has codimension 3. Since both are reflexive $R$-modules, it must be an isomorphism.

For simplicity, we denote by $\wh{\cN}$ the direct sum $\bigoplus_{i=1}^t\wh{\cN}_i$ and denote by $N$ the direct sum $\bigoplus_{i=1}^t N_i$.
\begin{theorem}(\cite[Corollary 3.2.10]{VdB04})\label{VdBeq}
The functor
$\wh{F}:=\RHom_{\wh{Y}}(\cO_{\wh{Y}}\op\wh{\cN} ,-)$ defines a triangle equivalence between $\D^b(\coh\wh{Y})$ and $\D^b(\Mod A)$, with quasi-inverse $F^{-1}:=(-)\ot_{A} (\cO_{\wh{Y}}\op\wh{\cN})$. In addition, $A$ is itself Cohen-Macaulay.
\end{theorem}
It follows that $A$ is a NCCR. 
\begin{corollary}\label{nccr-vs-contraction}
Let $(\hY,\hf,R)$ be a 3-dimensional formal flopping contraction. Then
\begin{enumerate}
\item[$(1)$] The structure sheaves of the irreducible components $\{C_i\}_{i=1}^t$ of $\Ex(\hf)$ form a semi-simple collection in $\coh(\wh{Y})$.
\item[$(2)$] The derived deformation algebra $\Ga^\hY_C$ of the collection $C:=\{C_i\}_{i=1}^t$ is linked by quasi-isomorphisms to $\Ga$ in Theorem \ref{VdBTh}.
\item[$(3)$] For any $i=1,\ldots,t$, $C_i$ is nc rigid.
\end{enumerate}
\end{corollary}
\begin{proof}
To prove that $\Hom_{\wh{Y}}(\cO_{C_i},\cO_{C_j})=0$ for $i\neq j$, we simply need to use the condition that $\Ex(\hf)$ is a tree of rational curves with normal crossings.

For any $i=1,\ldots,t$, it is easy to check that $S_i\cong \wh{F}(\Sg\cO_{C_i}(-1))$. Note that $(\bigoplus_{i=1}^t\cO_{C_i}) \ot (\bigotimes_{i=1}^t \cL_i^{-1})\cong \bigoplus_{i=1}^t\cO_{C_i}(-1)$.
By Theorem \ref{VdBeq}, $\wh{F}$ induces an isomorphism of $A_\infty$ algebras
\begin{equation}\label{isoext}
\Ext_\hY^*(\bigoplus_{i=1}^t \cO_{C_i},\bigoplus_{i=1}^t \cO_{C_i})\cong \Ext_A^*(\bigoplus_{i=1}^t S_i,\bigoplus_{i=1}^t S_i).
\end{equation}
Let $\ol{l}=A/\rad(A)$ and $l=\ol{l}/ ke_0$.  Then there is a natural isomorphism of $A$-modules $l\cong \bigoplus_{i=0}^tS_i$, where $S_0$ is the simple $A$-module that corresponds to the summand $R$ of $R\op \bigoplus_{i=1}^t N_i$. By Lemma 4.1 of \cite{ThV}, the vector space $V$ in Theorem \ref{VdBTh} can be chosen as $D\big(\Sg \Ext^{\geq 1}_A(\ol{l},\ol{l})\big)$. Therefore, $\Ga:=\wh{T}_{\ol{l}}V/\wh{T}_{\ol{l}}Ve_0\wh{T}_{\ol{l}}V\cong \bD B(\Ext^*_A(l,l))$ represents the noncommutative deformations of semi-simple collection $\{S_i\}_{i=1}^t\in \D^b(\mod A)$. Part $(2)$ follows from the isomorphism \ref{isoext}. Denote by $\Ga^\hY_{C_i}$ the derived deformation algebra of $\cO_{C_i}$. We have
\[
\Ga^\hY_{C_i}\cong \Ga/\sum_{j\neq i}\Ga e_j \Ga,~~~\text{and}~~~ H^0\Ga^\hY_{C_i}\cong H^0\Ga/\sum_{j\neq i}H^0\Ga e_j H^0\Ga.
\]
Then part $(3)$  follows from $(1)$ of Theorem \ref{VdBTh}. 
\end{proof}

\section{Calabi-Yau structure and cluster category}

\label{sec:CYcluster}
In this section, we first review several notions of Calabi--Yau property
for triangulated categories, for homologically smooth dg algebras and
for proper dg algebras. Then we recall geometric versions of the Calabi--Yau
property and translate them into algebraic notions for endomorphism
algebras of generators respectively for derived deformation algebras.
Finally, we classify Calabi--Yau structures for $3$-dimensional flopping
contractions and review the cluster category.

\subsection{CY structures}

\subsubsection{CY triangulated categories}
Let $\cT$ be a $\Hom$-finite $k$-linear triangulated category. 
\begin{definition} A {\em right Serre functor} for $\cT$ is a triangle functor
$S\colon \cT \to \cT$ such that there are bifunctorial isomorphisms
\[
\cT(Y,SX) \to D\cT(X,Y)
\]
for all $X,Y\in\cT$. It is a {\em Serre functor} if it is an autoequivalence.
\end{definition}
One can show that a right Serre functor exists if and only if for each object
$X$ of $\cT$, the functor $D\cT(X,?)$ is representable and in this case, the
right Serre functor is unique up to canonical isomorphism of triangle functors
\cite{BondalKapranov89, VdB08}.
Let $d$ be an integer. The triangulated category $\cT$ is \emph{$d$-Calabi--Yau} if
it admits a Serre functor isomorphic to $\Sigma^d$.

\subsubsection{CY smooth dg algebras}
A dg-$k$-algebra $\Ga$  is called 
{\em homologically smooth} if $\Ga$ is perfect in $\D(\Ga^e)$.
Then one checks that $\D_{fd}(\Ga)$, the subcategory of $\D(\Ga)$ consisting
of the modules whose homology is of finite total dimension, is contained
in the perfect derived category $\per(\Ga)$. Put 
\[
\Theta_\Ga = \RHom_{\Ga^e}(\Ga,\Ga^e).
\]
Then we have a canonical isomorphism
\[
\HoH_d(\Ga) \iso \Hom_{\D(\Ga^e)}(\Theta_\Ga,\Sg^{-d}\Ga).
\]
\begin{definition}
The dg algebra $\Ga$ is called \emph{bimodule} dCY if it is homologically
smooth and there is an isomorphism in 
$\D(\Ga^e)$
\[
\eta: \Theta_\Ga\iso \Sg^{-d}\Ga.
\]
A class $\eta\in \HoH_d(\Ga)$ is called a \emph{dCY structure} if the corresponding morphism $\eta: \Theta_\Ga\to\Sg^{-d}\Ga$ is an isomorphism in $\D(\Ga^e)$. A dCY structure $\eta$ is called \emph{exact} if there exists a class $\xi\in \HC_{d-1}(\Ga)$ such that $B\xi=\eta$, where $B$ is the Connes morphism. A choice of $\xi$ will be called an {\em{exact lifting}} of the dCY structure $\eta$. We call a bimodule dCY algebra $\Ga$ an \emph{exact dCY algebra} 
if the dCY structure is exact in addition.
\end{definition}

\begin{definition}
The dg algebra $\Ga$ is said to satisfy the \emph{relative dCY property} 
if for $L\in \D_{fd}(\Ga)$ and $M\in \per(\Ga)$, we have a bifunctorial isomorphism
\[
D\RHom_{\D(\Ga)}(L,M)\simeq \RHom_{\D(\Ga)}(M,\Sg^dL).
\] 
\end{definition} 

\begin{remark}
If $\Ga$ is a topologically homologically smooth pseudo-compact dg algebra in
$\PCAlgc(l)$, we call $\Ga$ a bimodule dCY if $\eta$ is an isomorphism in the pseudo-compact derived category of bimodules. 
The isomorphism $\eta$ represents a class in the continuous Hochschild homology 
$\HoH_d(\Ga)$. Exactness is defined similarly by taking the continuous cyclic homology. We call a bimodule dCY pseudo-compact dg algebra $\Ga$ in $\PCAlgc(l)$
an \emph{exact dCY algebra} if the dCY structure is exact in addition. Similarly, for a pseudo-compact algebra $\Ga$ in $\PCAlgc(l)$, 
the relative dCY property is defined by replacing 
$\D(\Ga), \per(\Ga)$ and $\D_{fd}(\Ga)$ with their pseudo-dompact 
counter-parts.
\end{remark}

Given a homologically smooth dg algebra $\Ga$, it follows from Lemma~3.4
in \cite{KV09} that we have the implications
\begin{align*}
\text{$\Ga$ is bimodule dCY} & \Rightarrow \text{$\Ga$ satisfies the relative dCY property} \\
& \Rightarrow \text{$\D_{fd}(\Ga)$ is a $\Hom$-finite dCY triangulated category}.
\end{align*}
A similar chain of implications holds in the pseudo-compact case.

\subsubsection{CY proper dg algebras} 
Let $A$ be a dg algebra.
Suppose that $A$ is {\em proper}, i.e. that its homology is of finite
total dimension. Then the category $\per(A)$ is $\Hom$-finite.
The proper dg algebra $A$ is called {\em perfectly $d$CY} if there is
an isomorphism
\[
DA \iso \Sigma^d A
\]
in $\D(A^e)$. By the following lemma, the triangulated category $\per(A)$ is then
$d$CY.

\begin{definition}
Let $k$ be a field of characteristic zero. Given a finite-dimensional
$A_\infty$-algebra $A$, a \emph{cyclic $A_\infty$-structure} of degree $d$ on $A$ is a non degenerate symmetric bilinear form
\[
(?, ?) : A \times A \to \Sg^dA
\] 
of degree $d$ such that 
\[
(m_n(a_1, \ldots , a_n), a_{n+1}) = (-1)^n(-1)^{|a_1|(|a_2|+\ldots +|a_{n+1}|)}(m_n(a_2, \ldots , a_{n+1}), a_1)
\]
\end{definition}
In this case, we have in particular an isomorphism $DA \iso \Sigma^d A$ in
the derived category of $A$-bimodules. Thus, if a dg algebra is quasi-isomorphic
to an $A_\infty$-algebra admitting a cyclic $A_\infty$-structure of degree $d$,
then it is perfectly $d$CY.

\begin{lemma} \label{lemma:dCY-Koszul-dual} Let $T\geq 1$ be an integer,
$l$ a product of $t$ copies of $k$ and $A$ a dg algebra augmented over $l$ and such that $A$ 
belongs to the triangulated subcategory of $\D(A \ten A^{op})$ generated by $l\ten l^{op}$
(in particular, the dg algebra $A$ is proper). 
Let $\Ga$ be the pseudo-compact dg algebra $\bD BA$, where $B$ denotes the bar 
construction over $l$. Then $\Ga$ is topologically homologically smooth. Moreover, 
if $A$ is perfectly $d$-Calabi--Yau, then $\Ga$ is bimodule $d$-Calabi--Yau.
\end{lemma}

\begin{proof} We refer to  \cite{KeKos}  for a brief summary of the bar-cobar formalism that
we are going to use and to \cite{Lefevre03} and Appendix A of \cite{VdB15}  for in-depth treatments.
Let $C$ be the augmented dg coalgebra $BA$ and 
$\tau: C \to A$ the canonical twisting cochain. Since $\tau$ is acyclic,
the canonical morphism
\[
C \to C \ten_\tau A \ten_\tau C
\]
is a weak equivalence of dg $C$-bicomodules (cf. \cite{KeKos}). Recall
(for example from Appendix A of \cite{VdB15}) that a dg comodule is
fibrant if its underlying graded comodule is cofree. Clearly, this
holds for $C\ten_\tau A \ten_\tau C$.
Since $A$ belongs to the triangulated subcategory of the derived category of $A\ten A^{op}$
generated by $l\ten l$, the dg bicomodule $C\ten_\tau A \ten_\tau C$ belongs
to the triangulated subcategory of the coderived category of $C\ten C^{op}$
generated by $C\ten C^{op}$. By applying the duality $\bD$ we
obtain that the morphism
\begin{equation} \label{bimodule-resolution}
\Ga \ten_\tau \bD A \ten_\tau \Ga \to \Ga
\end{equation}
is a weak equivalence and that the object on the left is cofibrant.
Moreover, we see that $\Ga$ belongs to the perfect derived category
of $\Ga\ten \Ga^{op}$. This means that $\Ga$ is
topologically homologically smooth. Now suppose that $A$ is perfectly
$d$-Calabi--Yau. Since $A$ is proper, it is weakly equivalent to its
pseudo-compact completion $\wh{A}$. By the Calabi--Yau property,
we have an isomorphism
$\bD A \iso \Sigma^d \wh{A}$ in the pseudo-compact derived category of
dg $\wh{A}$-bimodules. Now we compute the inverse dualizing complex of
$\Ga$ using the resolution \ref{bimodule-resolution}. We have isomorphisms
in the pseudo-compact derived category of $\Ga$-bimodules
\begin{align*}
\RHom_{\Ga^e}(\Ga,\Ga^e) &= \Hom_{\Ga^e}(\Ga\ten_\tau \bD A \ten_\tau \Ga, \Ga^e) \\
&= \Hom^{\tau,\tau}_{psc}(\bD A, \Ga^e) \\
&= \Hom^{\tau,\tau}_{psc}(\Sigma^d \wh{A}, \Ga^e) \\
&= \Sigma^{-d} \Ga \ten_\tau \bD A \ten_\tau \Ga \\
&= \Sigma^{-d} \Ga.
\end{align*}
Here $\Hom^{\tau,\tau}_{psc}$ denotes the space of morphisms in the
category of pseudo-compact vector spaces twisted twice by $\tau$.
This shows that $\Ga$ is topologically homologically bimodule $d$-Calabi--Yau.
\end{proof}

\subsubsection{CY structures in geometry, algebraic consequences}
We let $k=\CC$ be the field of complex numbers unless specified otherwise.
\begin{definition}
Let $Y$ be a $d$-dimensional smooth quasi-projective $\CC$-variety. We call $Y$ a $d$-dimensional Calabi--Yau variety if there is an isomorphism $\omega_Y:=\Omega^d_Y\cong\cO_Y$, i.e. there exists a nowhere vanishing $d$-form. We call a nowhere vanishing section $\eta: \cO_Y\to \omega_Y$ a \emph{dCY structure} on $Y$. We call the dCY structure \emph{exact} if the the $d$-form $\eta$ is exact, i.e. there exists a $(d-1)$-form $\xi\in \Omega^{d-1}_Y$ such that $\dd\xi=\eta$. A choice of $\xi$ will be called an {\em{exact lifting}} of the dCY structure $\eta$. If $\hY$ is a smooth formal scheme, we may define dCY structure in a similar way by considering de Rham complex of formal scheme $\Omega_\hY^*$. 
\end{definition}

Given a smooth quasi-projective variety $Y$ of dimension $d$ and a bounded complex of coherent sheaves $L$ on $Y$ whose cohomology has proper support, the derived endomorphism algebra $A:=\RHom_{\D(\Qcoh~Y)}(L,L)$ is a proper dga. We will show that $A$ is perfectly dCY if $Y$ is d-Calabi-Yau.

\begin{lemma}\label{SD}
Let $Y$ be a quasi-projective smooth CY $d$-fold with a fixed CY structure. Let $L\in \D^b_{c}(\coh Y)$ 
be a bounded complex of coherent sheaves with proper support. 
Then 
\[
A=\RHom_{\D(\Qcoh~Y)}(L,L)
\]
is perfectly d-Calabi-Yau.
\end{lemma}
\begin{proof}
Since we work over a field of characteristic zero, there exists $\ol{Y}$, a smooth projective variety that compactifies $Y$. 
Denote by $i:Y\to\ol{Y}$ the canonical embedding.  Since $i$ is an open immersion and $L$ has proper support on $Y$, $\RHom_{\D(\Qcoh \ol{Y})}(i_*L,i_*L)$ is quasi-isomorphic to $A$ as dg algebras. 

From now on, we simply assume that $A=\RHom_{\D(\Qcoh \ol{Y})}(i_*L,i_*L)$. We adopt the notation of the Appendix to write $\Qcoh(Y)$ for the dg category of (fibrant replacements of) unbounded complexes of quasi-coherent sheaves and $\coh(Y)$ for its full dg subcategory of complexes with coherent cohomology and bounded cohomological amplitude.
By Theorem \ref{Serredg} there is a bifunctorial quasi-isomorphism
\[
D\Hom^{dg}_{\coh(\ol{Y})}(M,N)\simeq \Hom^{dg}_{\coh(\ol{Y})}(N,M\otimes \omega_{\ol{Y}}[d]).
\]
Now let $M=N=i_*L$. Then $M\otimes \omega_{\ol{Y}}\cong i_*(L\ot\omega_Y)\cong i_*L$ using the Calabi-Yau structure on $Y$. By the bifunctoriality this is an isomorphism in $D(A^e)$.
\end{proof}

\begin{prop}\label{GaCY}
Let $\{\cE_i\}_{i=1}^t$ be a semi-simple collection of compactly supported sheaves in a smooth quasi-projective CY d-fold $Y$. Write 
$\Ga:=\Gamma^Y_{\cE}$ for the derived deformation algebra of  $\cE:=\bigoplus_{i=1}^t \cE_i$. Then $\Ga$ 
is topologically homologically smooth and bimodule $d$-Calabi--Yau.  
\end{prop}

\begin{proof} Let $A$ be the derived endomorphism algebra of $\cE$. Since $Y$ is smooth and $\cE$
has compact support, it is proper. Moreover, it can be chosen augmented. Clearly, it has its homology
concentrated in non negative degrees and $H^0(A)$ is isomorphic to a product of $t$ copies of
$k$, which is $l$. By Corollary~4.1 of \cite{KN13}, the dg algebra $A$ belongs to the triangulated
subcategory of its derived category generated by the semi-simple object $l$. 
We know that $\Ga$ is quasi-isomorphic to $\bD BA$. Thus, by Lemma~\ref{lemma:dCY-Koszul-dual},
$\Ga$ is topologically homologically smooth. 
Since $Y$ is smooth and $d$-Calabi--Yau, $A$ is perfectly $d$-Calabi--Yau
by Lemma~\ref{SD}. Hence $\Ga$ is 
bimodule $d$-Calabi--Yau by Lemma~\ref{lemma:dCY-Koszul-dual}.  
\end{proof}

\begin{remark}
If we assume that $Y$ is smooth and projective, then one can show that $\Ext_Y^*(\cE,\cE)$ is a cyclic $A_\infty$-algebra. This can be proved by reducing to the analytic case and applying the holomorphic Chern-Simons theory (see Example 10.2.7 \cite{KSnotes}). An algebraic proof of the case when $Y$ is a projective Gorenstein CY curve can be found in \cite{Pocyclic}.
\end{remark}

\begin{prop}
Let $\{\cE_i\}_{i=1}^t$ be a semi-simple collection of sheaves in a smooth projective CY d-fold $Y$. Write 
$\Ga:=\Gamma^Y_{\cE}$ for the derived deformation algebra of $\cE:=\bigoplus_{i=1}^t\cE_i$. 
Then $\Ga$ is an exact $d$-CY-algebra.
\end{prop}
\begin{proof} As in the proof of Lemma~\ref{GaCY}, we see that $\Ga$ is topologically homologically smooth
and bimodule $d$-CY. The exactness of the CY-structure follows from Theorem 12.1 of \cite{VdB15} and the remark above.
\end{proof}

\begin{theorem}(Iyama-Reiten)\label{IR}
Let $R$ be an equi-codimensional Gorenstein normal domain of dimension $d$ over an algebraically closed field $k$, and let $A$ be an NCCR. Then $A$ satisfies the relative dCY property. Moreover, if $R$ is complete local then $A$ is bimodule dCY. \end{theorem}

\begin{proof} The fact that $A$ satisfies the relative dCY property is proved in
Theorem~4.23 of \cite{Wem12}. We have $A=\End_R(R\oplus N)$ for the
Cohen-Macaulay module $N=\bigoplus_{i=1}^tN_i$ with $N_1,\ldots,N_t$ indecomposable. 
Denote by
$\ol{l}$ the algebra $k^{t+1}$. By Theorem~1.1 of \cite{ThV}, the 
algebra $A$ is quasi-isomorphic to a pseudo-compact dg algebra
$\wt{\Ga}:=(\wh{T}_{\ol{l}}(V),d)$ in $\PCAlgc(l)$ 
for a finite-dimensional graded $\ol{l}$-bimodule $V$ 
concentrated in degrees $\leq 0$ and
a differential $d$ taking $V$ into the square of the kernel of the augmentation
$\wh{T}_{\ol{l}}(V)\to \ol{l}$. Since $A$ is of finite global dimension, $\wt{\Ga}$ is
homologically smooth and by the first part, we have bifunctorial isomorphisms
\[
\Hom_{\D(\wt{\Ga})}(M,P)=D\Hom_{\D(\wt{\Ga})}(P,\Sigma^d M)
\]
for $M$ in $\D_{fd}(\wt{\Ga})$ and $P$ in $\per(\wt{\Ga})$. Let
\[
\Theta=\RHom_{\wt{\Ga}^e}(\wt{\Ga},\wt{\Ga}^e)
\]
be the inverse dualizing complex of $\wt{\Ga}$. By Lemma~4.1 of \cite{KeCY},
we have bifunctorial isomorphisms
\[
\Hom_{\D(\wt{\Ga})}(L\lten_{\wt{\Ga}} \Theta, M) = D\Hom_{\D(\wt{\Ga})}(M,L)
\]
for $M$ in $\D_{fd}(\wt{\Ga})$ and an arbitrary object $L$ of $\D(\wt{\Ga})$.
By combining these with the previous isomorphisms we find
\[
\Hom_{\D(\wt{\Ga})}(P\lten_{\wt{\Ga}} \Theta, M) = \Hom_{D(\wt{\Ga})}(\Sigma^{-d}P,M)
\]
for $P$ in $\per(\wt{\Ga})$ and $M$ in $\D_{fd}(\wt{\Ga})$. Since an object $L$
of $\D(\wt{\Ga})$ is perfect if and only if $\Hom_{\D(\wt{\Ga})}(L,M)$ is
finite-dimensional for each $M$ in $\D_{fd}(\wt{\Ga})$, we see that 
$P\lten_{\wt{\Ga}} \Theta$ is perfect for each perfect $P$. Now by taking
$P=\wt{\Ga}$ and $M=\Sigma^n S_i$, where $n\in\ZZ$ and $S_i$ is one of the
$t+1$ simple $\wt{\Ga}$-modules, we see that as a right
$\wt{\Ga}$-module, $\Theta$ is quasi-isomorphic to $\Sigma^{-d}\wt{\Ga}$.
For the rest of the argument, let us replace $\wt{\Ga}$ by the
quasi-isomorphic pseudocompact algebra $G=H^0(\wt{\Ga})$, which is isomorphic
to the original algebra $A$.
Since $G$ and $\wt{\Ga}$ have canonically equivalent derived
categories and derived categories of bimodules, we can view
$\Theta$ as a $G$-bimodule complex concentrated in degree
$d$. After replacing it with a quasi-isomorphic bimodule complex, we
may assume that it is an actual $G$-bimodule concentrated
in degree $d$. Moreover, we know that as a right module,
it is isomorphic to $G$. Thus, 
there is an automorphism $\sigma$ of $G$ such that
$\Sigma^d \Theta$ is isomorphic to $_\sigma G$ as a bimodule.
Since each object $L$ of $\D_{fd}(G)$ is perfect, we have 
\[
\Hom_{\D(G}(L_\sigma, M) = \Hom_{\D(G)}(L,M)
\]
for all $L$ and $M$ in $\D_{fd}(G)$, which shows that there is a functorial
isomorphism $L \iso L_\sigma$ for each $L$ in $\D_{fd}(G)$.
In particular, for $L$, we can take the finite-dimensional quotients of $G$.
We deduce that in each finite-dimensional quotient of $G$, the
automorphism $\sigma$ induces an inner automorphism. Thus,
$\sigma$ itself is inner and $_\sigma G$ is isomorphic to $G$
as a bimodule. This shows that $\Theta$ is quasi-isomorphic to $\Sigma^{-d}\wt{\Ga}$
as a bimodule.
\end{proof}

\begin{corollary}
Let $R$ be a local complete equi-codimensional Gorenstein normal domain of dimension $d$ over an algebraically closed field $k$ of characteristic zero, and let $A$ be an NCCR. Let 
$\Gamma$ be the dg algebra constructed in Theorem \ref{VdBTh}. Then $\Ga$ is topologically homologically smooth and bimodule dCY. 
\end{corollary}

\subsection{Classification of CY structures for 3-dimensional flopping contractions}
If $(\hY,\hf,R)$ is a 3-dimensional formal flopping contraction, then  $R$ is  a hypersurface. The natural isomorphism $\hf_*\omega_{\hY}\cong \omega_R$ identifies a CY structure $\eta$ on $\hY$ with a nonzero section $\hf(\eta)$ of $\omega_R$. By the Gorenstein property, $\hf(\eta)$ defines an isomorphism $R\cong \omega_R$.

\begin{theorem}\label{classificationofCY}
Let $(\hY,\hf,R)$ be a 3-dimensional formal flopping contraction. The space of 3CY structures can be identified with $R^\times$. Moreover, every Calabi-Yau structure on $\hY$ is exact. The space of all exact liftings of a 3CY structure can be identified with the cohomology group $H^1(\hY,\Omega^1_\hY)$.
\end{theorem}

\begin{proof}
Assume that $C:=\Ex(f)$ has $t$ irreducible components $C_1,\ldots,C_t$.
Because $R$ has rational singularities, $H^0(\hY,\Omega^3_\hY)\cong H^0(\hY,\cO_\hY)\cong R$. Then the first claim follows.
The Hodge-to-de Rham spectral sequence with $E_1$ term
\[
E_{1}^{pq}=H^q(\hY,\Omega_{\hY}^p),
\] converges to $H^{p+q}_{DR}(\hY,\CC)$.  We claim that $H^1(\Omega^1_{\hY})\cong\CC^t$. Because the first Chern classes of $\wh{\cL}_i$ for $i=1,\ldots,t$ are linearly independent, $\dim_\CC H^1(\Omega^1_{\hY})\geq t$. We write  $\hX:=\Spec R$. By the Leray spectral sequence
\[
H^p(\wh{X},\bR^q \hf_*\Omega^1)\Rightarrow H^{p+q}(\hY,\Omega^1),
\]
we have an exact sequence
\[
0\to H^0(\wh{X},\bR^1f_*\Omega^1) \to H^1(\hY,\Omega^1)\to H^1(\wh{X},f_*\Omega^1)\to 0
\]
The right most term vanishes since $\wh{X}$ is affine. Since $\hY$ is a small resolution of 3-dimensional Gorenstein singularities,  
the normal bundle of $C_i$ is $\cO(a)\op \cO(b)$ with $(a,b)=\{(-1,-1),(0,-2),(1,-3)\}$ (\cite[Theorem 4]{Pink83}). By the exact sequence
\[
0\to \cO(-a)\op \cO(-b)\to \Omega^1|_{C_i}\to \cO(-2)\to 0
\] we have $H^1(C_i,\Omega^1|_{C_i})\cong H^1(C_i,\cO(-2))=\CC$. By the normal crossing condition, there is a short exact sequence of sheaves
\[
\xymatrix{
0\ar[r] & \Omega^1|_C\to \Omega^1|_{C_1}\op \Omega^1|_{\cup_{i=2}^t C_i}\ar[r] & \Omega^1|p\ar[r] & 0
}
\] where $p=C_1\cap (\cup_{i=2}^t C_i)$. Then we get a surjection 
\[
H^1(C,\Omega^1|_C)\to H^1(C_1,\Omega^1|_{C_1})\op H^1(\cup_{i=2}^t C_i, \Omega^1_{\cup_{i=2}^t C_i}).
\]
By induction, $\dim_\CC H^1(C,\Omega^1|_C) \leq t$.
Then the conclusion follows from the formal function theorem and the Lerray spectral sequence. 

The term $E_2^{30}$ of the
Hodge-to-de Rham spectral sequence is the quotient $H^0(\hY,\Omega^3)/\dd H^0(\hY,\Omega^2)$. Recall that  $H^1(\hY,\Omega^1)$ admits a $\CC$-basis by $c_1$ of the line bundles $\wh{\cL}_i$ with $i=1,\ldots,t$. Because $(1,1)$ classes are $d$-closed, the differential $H^1(\hY,\Omega^1)\to H^0(\hY,\Omega^3)/\dd H^0(\hY,\Omega^2)$ is zero. Moreover since $H^2(\cO_{\hY})=0$, we have $E_r^{30}=E_2^{30}$ for $r\geq 2$. Because $H_{DR}^3(\hY,\CC)\cong H_{DR}^3(C,\CC)=0$, $E_2^{30}$ must vanish. Therefore, all 3-forms are exact.

Denote by $\Omega^{\leq i}_\hY$ the stupid truncation $\sigma_{\leq i}\Omega^*_\hY$ of the de Rham complex. There is a long exact sequence of hypercohomology  
\[\xymatrix{
\ldots\ar[r]^B & H^{i-n}(\hY,\Omega^i_\hY)\ar[r]^I & \HH^{2i-n}(\hY,\Omega_\hY^{\leq i})\ar[r]^S &\HH^{2i-n}(\hY,\Omega_\hY^{\leq i-1})\ar[r]^B & H^{i-n+1}(\hY,\Omega^i_\hY)\ar[r] &\ldots}
\]
Take $i=3$ and $n=4$. Let leftmost term vanishes and $\HH^2(\hY,\Omega_\hY^{\leq 3})=H^2_{DR}(\hY)\cong H^1(\hY,\Omega^1_\hY)$. So the last claim is proved.
\end{proof}

\begin{corollary}
Let $R$ be a complete local equi-codimensional Gorenstein normal domain of dimension $d$ over an algebraically closed field of characteristic zero, and let $A$ be an NCCR. 
Then every dCY structure on $A$ is exact. 
\end{corollary}
\begin{proof}
This is an immediate consequence of Theorem \ref{IR}, Theorem \ref{VdBeq} and Corollary 9.3 of \cite{VdB15}.
\end{proof}

The following proposition follows immediately from the Hochschild-Kostant-Rosenberg theorem.
\begin{prop}
Let $(\hY,\hf,R)$ be a 3-dimensional formal flopping contraction. Let $A=\End_R(R\op N)$ be the corresponding NCCR. Then there is a bijective correspondence between the space of 3CY structures (resp. exact 3CY structures) on $\hY$ and that of $A$.
\end{prop}

\subsection{Cluster category}\label{subsec:cluster category}
Let $\Ga$ be a dg $k$-algebra. Suppose that $\Ga$ has the following properties:
\begin{enumerate}
\item[$(1)$] $\Ga$ is homologically smooth, i.e. $\Ga$ is a perfect $\Ga^e$-module;
\item[$(2)$] for each $p>0$, the space $H^p\Ga$ vanishes;
\item[$(3)$] $H^0\Ga$ is finite dimensional;
\item[$(4)$] $\Ga$ satisfies the relative 3CY property.
\end{enumerate}
By property $(1)$, $\D_{fd}(\Ga)$ is a subcategory of the perfect derived category $\per (\Ga)$. The \emph{generalized cluster category} $\cC_\Ga$ is defined to be the triangle quotient $\per(\Ga)/\D_{fd}(\Ga)$. We denote by $\pi$ the canonical projection functor $\pi:\per(\Ga)\to \cC_\Ga$. For simplicity, we will omit the adjective `generalized' and call 
$\cC_\Ga$ the cluster category associated to $\Ga$. An object $T\in \cC_\Ga$ is called a \emph{cluster-tilting object} if  
\begin{enumerate}
\item[$(1)$] $\Ext_{\cC_\Ga}^1(T,T)=0$;
\item[$(2)$] For any object $X$ such that $\Ext^1_{\cC_\Ga}(T,X)=0$, one has $X\in \add(T)$.
\end{enumerate}
Amiot has proved that $\pi(\Ga)$ is a cluster-tilting object in $\cC_\Ga$ (Theorem 2.1 \cite{Am}).
We call $H^0\Ga$ the \emph{CY tilted algebra} associated to the cluster category 
$\cC_\Ga$, cf. \cite{Reiten10}.

\begin{remark}
If $\Ga$ is a pseudo-compact dg $l$-algebra in $\PCAlgc(l)$, 
we may define a continuous version of cluster category. Condition $(1)$ is replaced by 
\begin{enumerate}
\item[$(1^\p)$] $\Ga$ is topologically homologically smooth,
\end{enumerate}
and the \emph{topological} cluster category $\cC_\Ga$ is defined to be the triangle quotient $\per(\Ga)/\D_{fd}(\Ga)$ where $\per(\Ga)$ and $\D_{fd}(\Ga)$ are considered as subcategories of the pseudo-compact derived category. We refer to the Appendix of \cite{KY11} for the details.
\end{remark}

\begin{remark}
If we drop the assumption that $H^0\Ga$ is finite dimensional, the quotient category 
$\cC_\Ga=\per(\Ga)/\D_{fd}(\Ga)$ is no longer Hom-finite. The Calabi--Yau
property only holds when one restricts to suitable subcategories, cf.
Proposition~2.16 of \cite{Plamondon11}.
\end{remark}

\begin{theorem}\cite[Theorem 2.1]{Am}\footnote{In the original statement of \cite{Am}, the author assumed that $\Ga$ is bimodule 3CY. However, the proof is still valid under the weaker assumption that $\Ga$ satisfies the relative 3CY property.}\label{Amiot}
Let $\Ga$ be a dg $k$-algebra with the above properties. Then the cluster category 
$\cC_\Ga$ is Hom-finite and 2CY as a triangulated category. Moreover, the object 
$\pi(\Ga)$ is a cluster tilting object. Its endomorphism algebra is isomorphic to $H^0\Ga$.
\end{theorem}

\begin{definition}
Let $\{C_i\}_{i=1}^t$ be a collection of smooth rational curves in a smooth quasi-projective CY 3-fold $Y$ with fixed CY-structure $\eta: \cO_Y\iso \omega_Y$, such that $\{\cO_{C_i}\}$ form a semi-simple collection. Denote by $\cC(Y,\{C_i\}_{i=1}^t,\eta)$ for the (topological) cluster category associated to the derived deformation algebra of $\bigoplus_{i=1}^t\cO_{C_i}$. We call $\cC(Y,\{C_i\}_{i=1}^t,\eta)$ \emph{the cluster category associated to the triple $(Y,\eta,\{C_i\}_{i=1}^t)$}. 
\end{definition}

\begin{definition}
Let $R$ be a complete local equicodimensional Gorenstein normal domain of dimension $3$ over an algebraically closed field $k$ of characteristic zero, and let $A$ be the NCCR associated to the collection of indecomposables $R,N_1,\ldots, N_t$. Fix a 3CY structure 
$\eta\in \HoH_3(A,A)$. Denote by $\cC(R,\{N_i\}_{i=1}^t,\eta)$ the cluster category associated to the 
dg algebra $\Ga$ constructed in Theorem \ref{VdBTh}, and call it 
\emph{the cluster category associated to the triple $(R,\{N_i\}_{i=1}^t,\eta)$}. 
\end{definition}
A priori, the dg algebra $\Ga$ constructed in Theorem \ref{VdBTh} is pseudo-compact. However, if $\D(\Ga)$ denotes
the ordinary derived category and $\D_{pc}(\Ga)$ the pseudo-compact derived category,
then the natural functor $\D_{pc}(\Ga) \to \D(\Ga)$ induces equivalences in the perfect
derived categories and in the subcategories of objects with finite-dimensional total
homology. Therefore, the two candidates for the cluster category are equivalent.

The following result is an immediate consequence of Corollary \ref{nccr-vs-contraction}.
\begin{corollary}
Let $(\hY,\hf,R)$ be a 3-dimensional formal flopping contraction, and let $A$ be the NCCR associated to  the collection of indecomposables $R,N_1,\ldots, N_t$ constructed in Section \ref{sec:flop}. Fix a 3CY structure $\eta$ on $\hY$ and denote its counter part on $A$ by the same symbol. Then there is a triangle equivalence bewteen the cluster categories
\[
\cC(Y,\{C_i\}_{i=1}^t,\eta)\simeq \cC(R,\{N_i\}_{i=1}^t,\eta).
\]
\end{corollary}

\section{Ginzburg algebras}\label{sec:Ginz}
In this section we introduce the notion of Ginzburg (dg) algebra and prove several properties of it. The cluster 
category can be defined via the Ginzburg algebra, which provides an effective tool to do computations.

\subsection{Definitions}
Fix  a commutative ring $k$.
Let $Q$ be a finite quiver, possibly with loops and 2-cycles. Denote by $Q_0$ and $Q_1$ the set of nodes and arrows of $Q$ respectively. Denote by $kQ$ the path algebra and by $\wh{kQ}$ the complete path algebra with respect to the two-sided ideal generated by arrows. For each $a\in Q_1$, we define the cyclic derivative $D_a $ with respect to $a$ as the unique linear map 
\[
D_a: kQ/[kQ,kQ]\to kQ
\] which takes the class of a path $p$ to the sum $\sum_{p=uav} vu$ taken over all decompositions of the path $p$. The definition can be extended to $\wh{kQ}_\cy:=\wh{kQ}/[\wh{kQ},\wh{kQ}]^c$ where the superscript $c$ stands for the completion with respect to the adic topology defined above. An element $w$ in $\wh{kQ}/[\wh{kQ},\wh{kQ}]^c$ is called a \emph{potential} on $Q$. It is given by a (possibly infinite) linear combination of cycles in $Q$.

\begin{definition}(Ginzburg)
Let $Q$ be a finite quiver and $w$ a potential on $Q$. Let $\ol{Q}$ be the graded quiver with the same vertices as $Q$ and whose arrows are 
\begin{enumerate}
\item[$\bullet$] the arrows of $Q$ (of degree 0);
\item[$\bullet$] an arrow $a^*:j\to i$ of degree $-1$ for each arrow $a:i\to j$ of $Q$;
\item[$\bullet$] a loop $t_i:i\to i$ of degree $-2$ for each vertex $i$ of $Q$. 
\end{enumerate}
The (complete) Ginzburg (dg)-algebra $\fD(Q,w)$ is the dg $k$-algebra whose underlying graded algebra is the completion (in the category of graded vector spaces) of the graded path algebra $k\ol{Q}$ with respect to the two-sided ideal generated by the arrows of $\ol{Q}$. Its differential is the unique linear endomorphism homogeneous of degree 1 satisfying the Leibniz rule, and which takes the following values on the arrows of $\ol{Q}$:
\begin{enumerate}
\item[$\bullet$] $da=0$ for $a\in Q_1$;
\item[$\bullet$] $d(a^*)=D_aw$ for $a\in Q_1$;
\item[$\bullet$] $d(t_i)=e_i(\sum_{a\in Q_1}[a,a^*])e_i$ for $i\in Q_0$ where $e_i$ is the idempotent associated to $i$. 
\end{enumerate}
Denote by $l$ the product $\prod_{i\in Q_0} ke_i$. 
Then $\wh{kQ}$ is isomorphic to the complete tensor algebra $\wh{T}_l V$ with $V$ being the vector space spanned by arrows of $Q$. 
\end{definition}

\begin{remark}
In most references, the above definition corresponds to the \emph{complete Ginzburg algebra} while the algebra without taking the graded completion is called \emph{Ginzburg algebra}. The complete Ginzburg algebra $\fD(Q,w)$ is considered as an object of $\PCAlgc(l)$. Because we only consider complete Ginzburg algebra in this paper, we will call it the Ginzburg algebra for simplicity. 
\end{remark}

\begin{definition}
Let $Q$ be a finite quiver and $w$ a potential on $Q$. The \emph{Jacobi algebra} $\Lm(Q,w)$ is the zeroth homology of $\fD(Q,w)$, which is the quotient algebra
\[
\wh{kQ}/((D_aw, a\in Q_1))^c
\] where $((D_aw, a\in Q_1))^c$ is the closed two-sided ideal generated by $D_aw$. A Ginzburg algebra $\fD(Q,w)$ is called \emph{Jacobi-finite} if $\dim_k \Lm(Q,w)<\infty$.
\end{definition}

Van den Bergh showed the following result.
\begin{theorem}(Van den Bergh)\cite[Appendix]{KV09}
Let $Q$ be finite quiver and $w$ be a potential. Then $\fD(Q,w)$ is topologically homologically smooth and bimodule 3CY. 
\end{theorem}
The above theorem was first proved by Van den Bergh 
in the algebraic setting in \cite{KV09}. 
But the same proof  can be adapted to the pseudo-compact case
(cf. \cite{VdB15}).
Given a Jacobi-finite Ginzburg algebra $\Ga:=\fD(Q,w)$, there is an associated cluster category $\cC_\Ga:=\per( \Ga)/\D_{fd}(\Ga)$. 

\begin{remark}\label{canonical CY on Ginz}
There exists a canonical exact CY structure on $\Ga=\fD(Q,w)$. We follow the notation of \cite{VdB15} to write $M_l:=M/[l,M]$ for a $l$-bimodule $M$.
Because the reduced cyclic homology of $\Ga$ is equal to the homology of $(\Ga/l+[\Ga,\Ga])_l$ (see Proof of Theorem 11.2.1 of \cite{VdB08}), a class of $\HC_2(\Ga,\Ga)$ is represented by a degree $-2$ element $\chi$ of  $\Ga$ such that $d\chi \in l+[\Ga,\Ga]$. By the definition of $d$ of $\Ga$, $\chi:=\sum_{i\in Q_0} t_i$ represents a class in $\HC^{\rm{red}}_2(\Ga,\Ga)$. Because $\Ga$ is cofibrant, by Proposition 7.2.1 of \cite{VdB08} the Hochschild chain complex of $\Ga$ is quasi-isomorphic to the mapping cone of
\[
\xymatrix{
\Omega^1_l \Ga/[\Ga,\Omega^1_l \Ga]\ar[r]^-{\partial_1} & \Ga/[l,\Ga]
}
\] with differential defined by $\partial_1 (aDb)=[a,b]$, where $Db=1\ot b-b\ot 1$. In other words, a class in $\HoH_3(\Ga,\Ga)$ is represented by a pair of elements $(\omega,a)$ of degree $(-2,-3)$ satisfying 
$\partial_1(\omega)=da$ and $d\omega=0$. Because $d$ and $D$ commute, 
$(D\chi,0)$ represents a class in $\HoH_3(\Ga,\Ga)$, which is the 
image of $\chi$ under the Connes map. In \cite{VdB08}, a class in $\HoH_d(\Ga,\Ga)$ is called {\em{non-degenerate}} if the corresponding mormphism $\Theta_\Ga\to \Sigma^{-d}\Ga$ is an isomorphism.
By Lemma 11.1.2 of \cite{VdB08}, $(D\chi,0)$ is non-degenerate. 

\end{remark}

For a Ginzburg algebra $\Ga=\fD(Q,w)$, denote $\Lm(Q,w)$ by $\Lm$ for short. 
The image of $w$ under the canonical map $\HoH_0(\wh{kQ},\wh{kQ})=\wh{kQ}_\cy\to \HoH_0(\Lm,\Lm)=\Lm_\cy$, denoted by $[w]$, is a canonical class associated to the Ginzburg algebra $\fD(Q,w)$. Therefore, we see that starting from a Ginzburg algebra $\Ga=\fD(Q,w)$ we get not only a triangulated category $\cC_\Ga$ but an additional piece of information that is a canonical class $[w]$ in the 0-th Hochschild homology of the CY tilted algebra. We will show in the next section that this class is determined by the CY structure up to right equivalences.

\subsection{Existence and uniqueness of potential}
The definition of Ginzburg algebra is not homotopically invariant.
It is important to know when a bimodule 3CY dg algebra admits a model given by a Ginzburg algebra.
The following theorem is due to Van den Bergh.

\begin{theorem}\cite[Theorem 10.2.2]{VdB15}\label{VdBpotential}
Let $k$ be a field and $l$ be a finite dimensional commutative separable $k$-algebra.
Assume that $\Ga$ is a pseudo-compact dg $l$-algebra in $\PCAlgc(l)$
concentrated in nonpositive degrees. Then the following are equivalent
\begin{enumerate}
\item[$(1)$] $\Ga$ is exact 3CY.
\item[$(2)$] $\Ga$ is weakly equivalent to a Ginzburg algebra $\fD(Q,w)$ for some finite quiver $Q$ with $w$ contains only cubic terms and higher.
\end{enumerate}
\end{theorem}
The following result of Van den Bergh provides a lot of examples of dg algebras whose 3CY structures can be lifted to exact ones.

\begin{theorem}\cite[Corollary 9.3]{VdB15}\label{VdBcomplete=>exact}
Assume that $k$ has characteristic zero and let $\Ga$ be a pseudo-compact dg algebra in $\PCAlgc(l)$ concentrated in degree zero. 
Then $\Ga$ is bimodule dCY if and only if it is exact dCY.
\end{theorem}
By putting the above two theorems together, we see that if $\Ga$ is a pseudo-compact dg $l$-algebra in $\PCAlgc(l)$ concentrated in degree zero that is bimodule 3CY then it is 
quasi-isomorphic to a Ginzburg algebra $\fD(Q,w)$ for some finite quiver $Q$ 
and potential $w$. 

 Given a pseudo-compact  $l$-algebra $A$, two classes $[w]$ and $[w^\p]$ in $A/[A,A]^c$ are called \emph{right equivalent} if there exists an automorphism $\gamma$ of $A$ such that $\gamma_*[w]=[w^\p]$.
Now we assume that a bimodule 3CY dg algebra $\Ga$ in $\PCAlgc(l)$ 
is exact. So it admits a model given by $\fD(Q,w)$. Note that being bimodule CY and exact CY are homotopically invariant properties.
The next proposition shows that the right equivalence class of $[w]$ in $\HoH_0(\Lm(Q,w))$ for such a dg algebra is indeed a homotopy invariant. The proof is implicitly contained in Van den Bergh's proof of Theorem \ref{VdBpotential} (cf. proof of Theorem 11.2.1 of \cite{VdB15}). We recall it for completeness. See Remark \ref{Darboux} for a conceptional explanation of Van den Bergh's result.

\begin{prop}(Van den Bergh)\label{canonicalclass}
Let $k$ be a field of characteristic zero and $l=ke_1\times\ldots \times ke_n$. 
Let $\Ga$ be a pseudo-compact 3CY dg $l$-algebra in $\PCAlgc(l)$ with a  fixed exact 3CY structure  $\tau\in\HoH_3(\Ga,\Ga)$.  Suppose there are two pairs $(Q,w)$ and $(Q^\p,w^\p)$ such that $\Gamma$ is quasi-isomorphic to $\fD(Q,w)$ and $\fD(Q^\p,w^\p)$ respectively. Assume that under these quasi-isomorphisms we have $B\chi$ and $B\chi^\p$ identified with $\tau$ where $\chi, \chi^\p$ are the canonical classes in $\HC_2$ defined in Remark \ref{canonical CY on Ginz}. Then $Q=Q^\p$ and $w$ is right equivalent to $w^\p$. In particular, the classes $[w],[w^\p]\in H^0\Ga/[H^0\Ga,H^0\Ga]^c$ are right equivalent.
\end{prop}
\begin{proof}
In Theorem 11.2.1 of \cite{VdB08}, it is proved that there is a weak equivalence $(\wh{T}_lV,d)\to \Ga$ where $V=\Sigma^{-1}(D \Ext^*_\Ga(l,l))_{\leq 0}$ such that
\begin{enumerate}
\item[(1)] $V=V_c+l z$ with $V_c=\Sigma^{-1}(D\Ext^1_\Ga(l,l)\oplus D\Ext^2_\Ga(l,l))$ and $z$ being an $l$-central element of degree $-2$. 
\item[(2)] $dz=\sigma^\p\eta\sigma^\pp$ with $\eta \in (V_c\ot_l V_c)_l$ being a non-degenerate and anti-symmetric element. Here $\sigma^\p\ot \sigma^\pp\in l\ot_k l$ is the Casimir element corresponds to the trace $l\to k$ (see Section 5 of \cite{VdB08}).
\end{enumerate}
Since $Q$ (resp. $Q^\p$) depends only on $l$-bimodule structure on $\Ext^1_\Ga(l,l)$ (see \cite[Section 10.3]{VdB04}, we have $Q=Q^\p$.
Note that the perfect pairing on $V_c$ is determined by the  bimodule 3CY structure $\tau$ but does not depend on the exact lifting. Using the perfect pairing, any non-degenerate and anti-symmetric element in $(V_c\ot_lV_c)_l$ can be reduced to a canonical form by choosing appropriate basis on $\Ext^1_\Ga(l,l)$. The element $\eta$ defines a bisymplectic form $\omega_\eta$ of degree $-1$ on $\wh{T}_l V$ (see definition in Section 10.1 of \cite{VdB08}). By Lemma 11.3.1 of \cite{VdB08} there exists a potential $w\in \wh{T}_lV/[\wh{T}_lV,\wh{T}_l V]$ of degree zero such that for any $f\in \wh{T}_lV$
\[
df=\{w, f\}_{\omega_\eta}
\] where $\{-,-\}_{\omega_\eta}$ is the Poisson bracket associate to the bisymplectic form $\omega_\eta$. 
Since $w$ does not have constant terms, it is uniquely determined by the derivation $d=\{w,-\}_{\omega_\eta}$. In other words, any two potentials $w$ and $w^\p$ without constant terms satisfying the above equation differ by an automorphism $\wh{T}_l V\to \wh{T}_l V$. Moreover, since $w$ and $w^\p$ are of degree zero, therefore does not involve variables in $\Ext^{\geq 2}_\Ga(l,l)$, this automorphism is precisely a formal change of variables on $\Ext^1_\Ga(l,l)$. Such a formal change of variable induces an isomorphism from 
$\fD(Q,w)$ to $\fD(Q,w^\p)$ (see \cite[Theorem 4.3]{HuaZhou}), therefore an automorphism of the complete path algebra of $Q$ and an automorphism of $H^0\Gamma$.
\end{proof}

\begin{remark}\label{Darboux}
Theorem 11.2.1 of \cite{VdB08}, can be viewed as a Darboux-Weinstein theorem in noncommutative formal symplectic geometry. On $\Ext_\Ga^*(l,l)$, the cyclic $A_\infty$-structure can be interpreted as a symplectic structure. The symplectic structure restricts to the truncation $\Ext_\Ga^1(l,l)\op \Ext_\Ga^2(l,l)$ so that $\Ext_\Ga^1(l,l)$ is a (graded) Lagrangian. Then  Theorem 11.2.1 of \cite{VdB08} says that there exists a coordinate on $\Ext_\Ga^1(l,l)$ under which the symplectic form can be normalized so that it has constant coefficients, which is in particular exact. The differential $d$ of $\Ga$ can be interpreted as a homological vector field of degree 1. Then the contraction of the normalized symplectic form by $d$ is the exterior derivation of a potential  $w$ of degree 0. Note that a different choice of Darboux coordinates can only differ by a change of variables on $\Ext^1_\Ga(l,l)$, which leads to the above proposition.
\end{remark}

\begin{corollary}
Let $(\hY,\hf,R)$ be a 3-dimensional formal flopping contraction with reduced exceptional fiber $\Ex(f)=\bigcup_{i=1}^t C_i$, and let $A=\End_R(\bigoplus_{i=1}^t N_i\op R)$ be the NCCR associated to it. Fix a 3CY structure $\eta$ on $\hY$, and therefore on $A$. Denote the CY tilted algebra of 
$\cC(Y,\{C_i\}_{i=1}^t,\eta)\simeq \cC(R,\{ N_i\}_{i=1}^t,\eta)$ by $\Lm$. Then there exists a canonical class, defined up to right equivalence, on $\HoH_0(\Lm)=\Lm/[\Lm,\Lm]^c$ represented by a potential. 
\end{corollary}

The canonical class $[w]$ in the 0-th Hochschild homology of $H^0\Ga$ is part of the ``classical shadow'' of the CY structure. The class plays a crucial role in the geometric applications.  When $\Ga$ is weakly equivalent to a Jacobi-finite Ginzburg algebra $\fD(\wh{F},w)$ for a complete free algebra $\wh{F}=k\lgg x_1,\ldots,x_n\rgg$, then this class vanishes if and only if $w$ is right equivalent to a weighted homogeneous noncommutative polynomial (see Theorem \ref{ncSaito}). Therefore, the quasi-homogeneity of a potential is indeed a homotopy invariant of the CY algebra.
This motivates the following definition.

\begin{definition}
Let $k$ be a field and $l$ be a finite dimensional commutative separable $k$-algebra, and let $\Ga$ be a pseudo-compact dg $l$-algebra in $\PCAlgc(l)$ concentrated in nonpositive degrees. Assume that $\Ga$ is exact 3CY. Then $\Ga$ is called \emph{quasi-homogeneous} if the canonical class $[w]$ is (right equivalent to) zero.
\end{definition}

The notion of quasi-homogeneity is expected to be independent of choices of CY structure. In the case of simple flopping contractions, the first author and Gui-song Zhou have conjectured that this notion of quasi-homogeneity is indeed equivalent to the quasi-homogeneity of the underlying hypersurface singularity $R$ (see Conjecture 4.18 \cite{HuaZhou}).

\subsection{Properties of Jacobi-finite Ginzburg algebras}
In this section, we collect several results about Jacobi-finite Ginzburg algebras. 
We take $k$ to be the field of complex numbers,  though some of the results are valid more generally.

\begin{theorem}\cite[Theorem 3.16]{HuaZhou}
Let $Q$ be a finite quiver and $w$ be a potential in $\wh{\CC Q}_\cy$. Assume that the Jacobi algebra $\Lm(Q,w)$ is finite dimensional. Then $w$ is right equivalent to a formal series with only finitely many nonzero terms.
\end{theorem}
As a consequence, we may assume the potential is a noncommutative polynomial to begin with if the Jacobi algebra is known to be finite dimensional.

\begin{theorem}(Noncommutative Mather-Yau theorem)\cite[Theorem 3.5]{HuaZhou}\label{ncMY}
Let $Q$ be a finite quiver and let $w, w^\p \in \wh{\CC Q}_\cy$ be two potentials with only cubic terms and higher. Suppose that the Jacobi algebras $\Lm(Q , w)$ and $\Lm(Q ,w^\p)$ are both finite dimensional. Then the following two statements are equivalent:
\begin{enumerate}
\item[$(1)$] There is an algebra isomorphism $\gamma : \Lm(Q , w) \cong \Lm(Q , w^\p)$ so that $\gamma_*([w]) = [w^\p]$ in $\Lm(Q,w^\p)_\cy$.
\item[$(2)$] $w$ and $w^\p$ are right equivalent in $\wh{\CC Q}_\cy$.
\end{enumerate}
\end{theorem}

The noncommutative Mather--Yau theorem has an immediate application
to Ginzburg algebras.
\begin{corollary}\cite[Theorem 4.3]{HuaZhou}\label{iso-Ginz}
Fix a finite quiver $Q$. Let $w,w^\p\in \wh{\CC Q}_\cy$ be two potentials with only cubic terms and higher, such that the Jacobi algebras $\Lm(Q , w)$ and 
$\Lm(Q, w^\p)$ are both finite dimensional.  Assume there is an algebra isomorphism 
$\gamma : \Lm(Q ,w) \to \Lm(Q ,w^\p)$ so that $\gamma_*([w]) = [w^\p]$. Then there exists a dg algebra isomorphism
\[\xymatrix{
\beta : \fD(Q , w)\ar[r]^\cong & \fD(Q , w^\p)}
\]
such that $\beta(t_i) = t_i$ for any $i\in Q_0$.
\end{corollary}

\begin{definition}\label{def:wh}
Fix $\wh{F}$ to be the complete free associative algebra $\CC\lgg x_1,\ldots, x_n\rgg$.
Let $(r_1,\ldots, r_n)$ be a tuple of rational numbers. A potential $w\in \wh{F}_{\cy}:=\wh{F}/[\wh{F},\wh{F}]^c$ is said to be \emph{weighted-homogeneous of type $(r_1,\ldots,r_n)$} if it has a representative which is a  linear combination of monomials $x_{i_1}x_{i_2}\cdots x_{i_p}$ such that $r_{i_1}+r_{i_2}+\ldots r_{i_p}=1$.
\end{definition}

\begin{theorem}(Noncommutative Saito theorem)\label{ncSaito}\cite[Theorem 1.2]{HuaZhou2}
Let $w\in \wh{F}_{\cy}$ be a  potential with only cubic terms and higher such that the Jacobi algebra associated to $w$ is  finite dimensional. Then $[w]=0$ if and only if  $w$ is right equivalent to a weighted-homogenous potential of type $(r_1,\ldots,r_n)$ for some rational numbers $r_1,\ldots,r_n$ which lie strictly between $0$ and $1/2$. 
Moreover, in this case, all such types $(r_1, \ldots, r_n)$ agree with each other up to permutations on the indices $1,\ldots, n$.
\end{theorem}

Recall that $\cC_\Ga$ is constructed as the Verdier quotient of the perfect
derived category of $\Ga$ by its full subcategory of dg modules whose
homology is of finite total dimension. The category of perfect dg $\Ga$-modules
has a canonical dg enhancement and we obtain a natural dg enhancement
$(\cC_\Ga)_{dg}$ for $\cC_\Ga$ using the Drinfeld quotient of 
the dg category of strictly perfect dg $\Ga$-modules by its
full subcategory on the dg modules whose homology is of finite
total dimension.
\begin{theorem}\label{Gamma=truncation}
Let $Q$ be a finite quiver and $w$ a Jacobi-finite potential on $Q$.
Let $\Ga$ be the complete Ginzburg algebra associated with $(Q,w)$.
Denote by $\Lm_{dg}$  the dg endomorphism algebra of $\Ga$ in 
the canonical dg enhancement of the cluster category $\cC_\Ga$. 
Then there is a canonical isomorphism in the homotopy category of dg algebras
\[
\xymatrix{
\Ga\ar[r]^-{\sim} & \tau_{\leq 0} \Lm_{dg}.
}
\]
\end{theorem}
\begin{proof} 
There is a canonical morphism
\[
\Ga=\RHom_{\Ga}(\Ga,\Ga)\to(\cC_\Ga)_{dg}(\Ga,\Ga)=\Lm_{dg}
\]
where the right hand side denotes the dg endomorphism algebra
of $\Ga$ in the canonical dg enhancement $(\cC_\Ga)_{dg}$ of $\cC_\Ga$.
It suffices to show that the canonical map
\[
H^{-p}(\Gamma) = \Hom_{\per(\Gamma)}(\Ga,\Sigma^{-p}\Ga) \to
\Hom_{\cC_\Ga}(\Ga, \Sigma^{-p}\Ga)
\]
is invertible for $p\geq 0$. By Proposition~2.8 of \cite{Am}, we have
\[
\Hom_{\cC_\Ga}(\Ga, \Sigma^{-p}\Ga) = 
\colim_n \Hom_{\per(\Ga)}(\tau_{\leq n} \Ga, \tau_{\leq n} (\Sigma^{-p} \Ga)).
\]
We have
\[
\Hom_{\per(\Ga)}(\tau_{\leq n} \Ga, \tau_{\leq n} (\Sigma^{-p} \Ga)) =
\Hom_{\per(\Ga)}(\tau_{\leq n} \Ga, \Sigma^{-p} \Ga).
\]
Consider the canonical triangle
\[
\xymatrix{
\tau_{\leq n} \Ga \ar[r] & \Ga \ar[r] & \tau_{>n} \Ga \ar[r] & \Sigma (\tau_{\leq n} \Ga).}
\]
Recall that by Lemma 2.5 of \cite{Am}, the spaces $H^p(\Ga)$ are finite-dimensional
for all $p\in\ZZ$. Therefore, 
the object $\tau_{>n}\Ga$ belongs to $\D_{fd}(\Ga)$ and $\Sigma^{-p} \Ga$
belongs to $\per(\Ga)$. 

By the $3$-Calabi-Yau property, we have
\[
\Hom_{\per(\Ga)}(\Sigma^{-1} \tau_{>n} \Ga, \Sigma^{-p}\Ga) = 
D\Hom_{\per(\Ga)}(\Ga, \Sigma^{p+2} \tau_{>n}\Ga))
\]
which vanishes because $\tau_{>n}\Ga$ has no homology in degrees $>0$.
Similarly, we have
\[
\Hom_{\per(\Ga)}(\tau_{>n}\Ga, \Sigma^{-p}\Ga) = D\Hom_{\per(\Ga)}(\Ga, \Sigma^{p+3}(\tau_{>n} \Ga))
\]
which vanishes for the same reason. Thus we have
\[
\Hom_{\per(\Ga)}(\tau_{\leq n} \Ga, \Sigma^{-p} \Ga) =\Hom_{\per(\Ga)}(\Ga, \Sigma^{-p} \Ga).
\]
\end{proof}

\begin{corollary}\label{selfinj} 
Let $Q$ be a  quiver with one node and arbitrary number of loops and $\Ga=\fD(Q,w)$ a Jacobi-finite Ginzburg algebra. Then $H^0\Ga$ is self-injective and
there is an isomorphism 
\[\xymatrix{
\Sg^2\Ga\ar[r]^-\sim & \tau_{\leq -1}\Ga
}\]
in the derived category of dg $\Ga$-modules. In particular, we have
$H^i(\Ga)=0$ for odd $i$ and $H^i(\Ga)\cong H^0(\Ga)$ for even $i\leq 0$.
\end{corollary}
\begin{proof} By \cite{Am}, the cluster category $\cC_\Gamma$ is a $\Hom$-finite
$2$-Calabi-Yau category and the image $T$ of
$\Ga$ in $\cC_\Ga$ is a cluster-tilting object in $\cC_\Ga$.
By Theorem~4.1 of \cite{AdachiIyamaReiten14}, the cluster-tilting objects
of $\cC_\Ga$ are in bijection with the support $\tau$-tilting modules over
$\End(T)$. Since $\End(T)=H^0(\Ga)$ is local, the only support $\tau$-tilting
modules over $\End(T)$ are $0$ and $\End(T)$, by example~6.1 of [loc. cit.].
Thus, the only cluster-tilting objects of $\cC_\Ga$ are $T$ and $\Sigma T$.
In particular, $\Sigma^2 T$ has to be isomorphic to $T$ (since
$\Hom(\Sigma T, \Sigma^2 T)=0$ and $\Sg^2T$ must be a cluster-tilting object). This implies that $H^0\Ga=\End(T)$ is
self-injective, since, by the $2$-Calabi-Yau property, 
we have an isomorphism of right $\End(T)$-modules
\[
D \Hom(T,T) = \Hom(T,\Sigma^2 T) = \Hom(T,T).
\]
Let $\phi: \Sigma^2 \Ga \to \Ga$ be a lift of an isomorphism
$\Sigma^2 T \to T$ in $\cC_\Ga$. Let $p\geq 2$. In the commutative
square
\[
\xymatrix{
\Hom_{\per(\Ga)}(\Sigma^p\Ga, \Sigma^2\Ga) \ar[d] \ar[r] &
\Hom_{\cC_\Ga}(\Sigma^p T, \Sigma^2 T) \ar[d]^{\phi_*} \\
\Hom_{\per(\Ga)}(\Sigma^p\Ga, \Ga) \ar[r] &
\Hom_{\cC_\Ga}(\Sigma^p T, T)}
\]
the horizontal arrows are isomorphisms by Theorem~\ref{Gamma=truncation}
and the right vertical
arrow $\phi_*$ is an isomorphism. Thus, the morphism $\phi: \Sigma^2 \Ga \to \Ga$
induces isomorphisms in $H^i$ for $i\leq -2$. Moreover, we have
$H^{-1}(\Ga)=\Hom(T,\Sigma^{-1}T)=\Hom(T,\Sigma T)=0$ since $\Sigma^{-1}T$
is isomorphic to $\Sigma T$.  It follows that $\phi$ induces an isomorphism
\[
\Sigma^2 \Ga \iso \tau_{\leq -1} \Ga.
\]
\end{proof}

\begin{remark}
For the pair $(Q,w)$ associated to 3-dimensional flopping contractions, one can show that $H^0\Ga$ is indeed symmetric  (see Theorem \ref{Z2clustercat}). In the context of general contraction with one dimensional fiber, Kawamata has proved that the classical (multi-pointed) deformation algebra of the reduced exceptional fiber is always self-injective (see Proposition 6.3 \cite{KaComp}). So in particular it is Gorenstein. This result overlaps with the above corollary in the case of simple flopping contractions. For a general finite quiver $Q$, the 0-th homology of a Jacobi-finite Ginzburg algebra $\fD(Q,w)$ is not self-injective. Moreover, Kawamata proves that the deformation algebra is always isomorphic to its opposite algebra (see Corollary 6.3 of \cite{Ka18}). 
\end{remark}

\begin{corollary} Let $\Ga$ be the Ginzburg algebra of a Jacobi-finite quiver 
with potential. Let $T$ be the image of $\Ga$ in the cluster category $\cC_\Ga$.
\begin{itemize}
\item[a)] $H^0\Ga$ is selfinjective if and only if $H^{-1}\Ga$ vanishes
if and only if $T$ is isomorphic to $\Sigma^2 T$ in $\cC_\Ga$.
\item[b)] If the identity functor of $\cC_\Ga$ is isomorphic to $\Sigma^2$, then
$H^0\Ga$ is symmetric and there is an isomorphism
of graded algebras $H^0(\Ga)\ten k[u^{-1}] \iso H^*(\Ga)$, where $u$ is
of degree $2$.
\end{itemize}
\end{corollary}

\begin{proof} a)
By Theorem~\ref{Gamma=truncation}, the space $H^{-1}\Ga$ is isomorphic
to $\cC_\Ga(T,\Sigma^{-1}T)$ and $H^0\Ga$ is isomorphic to the endomorphism
algebra of $T$.
By Proposition~3.6 of \cite{IyamaOppermann13}, the endomorphism
algebra is selfinjective if and only if $\cC_\Ga(T,\Sigma^{-1}T)$ vanishes
if and only if $T$ is isomorphic to $\Sigma^2 T$ in $\cC_\Ga$.

b) By combining the functorial isomorphism from $T$ to $\Sigma^2 T$
with the Calabi-Yau property we get an isomorphism of bimodules over
the endomorphism algebra of $T$
\[
\cC_\Ga(T,T) \iso \cC_\Ga(T,\Sigma^2 T) \iso D \cC_\Ga(T,T).
\]
Since $H^0\Ga$ is in particular selfinjective, the space $H^{-1}\Ga$
vanishes by a). We get an isomorphism of graded algebras
\[
k[u,u^{-1}] \ten_k \cC_\Ga(T,T) \iso \bigoplus_{p\in\ZZ} \cC_\Ga(T,\Sigma^p T)
\]
where $u$ is of degree $2$ by sending $u$ to the functorial isomorphism
$T \iso \Sigma^2 T$. Thanks to Theorem~\ref{Gamma=truncation}, by
truncation, we get an isomorphism of graded algebras
\[
k[u^{-1}]\ten_k H^0\Ga \iso H^*\Ga.
\]
\end{proof}

\subsection{Silting theory for a non positive dg algebra and its
zeroth homology} \label{subsec:silting}
Let $\cT$ be a 
triangulated category. Recall that a {\em tilting object} for $\cT$ is a classical
generator $T$ of $\cT$ such that $\cT(T,\Sigma^p T)$ vanishes for all $p\neq 0$.
A {\em silting object} \cite{KellerVossieck88}
for $\cT$ is a classical generator $T$ of $\cT$ such that
$\cT(T,\Sigma^p T)$ vanishes for all $p>0$. The advantage of silting objects
over tilting objects is that (under suitable finiteness assumptions) they are
stable under mutation \cite{AiharaIyama12}.

We recall fundamental definitions and results from \cite{AiharaIyama12}.
Assume from now on that $\cT$ is $k$-linear, $\Hom$-finite and has split idempotents.
In particular, it is a Krull--Schmidt category, i.e.  indecomposables have
local endomorphism rings and each object is a finite direct sum of
indecomposables (which are then unique up to isomorphism and permutation). 
An object of $\cT$ is {\em basic}
if it is a direct sum of pairwise non isomorphic indecomposables. If
$X$ is an object of $\cT$ and $\cU$ a full additive subcategory stable
under retracts, a {\em left
$\cU$-approximation of $X$} is a morphism $f: X \to U$ to an object of
$\cU$ such that each morphism $X \to V$ to an object of $\cU$ factors
through $f: X \to U$. It is {\em minimal} if each endomorphism $g: U \to U$
such that $g\circ f=f$ is an isomorphism. Notice that the morphism
$f: X\to U$ is a minimal left $\cU$-approximation iff the morphism
$f^* : \cU(U,?) \to \cT(X,?)|_\cU$ is a projective cover in the category of
left $\cU$-modules. In particular, minimal left approximations are
unique up to non unique isomorphism when they exist. Existence
is automatic if $\cU$ has finitely many indecomposables $U_1$, 
\ldots, $U_n$ (which is
the case in our applications) because then the functor $\cT(X,?)|_\cU$
corresponds to a finite-dimensional left module over the finite-dimensional
endomorphism algebra of the sum of the $U_i$.
A {\em (minimal) right $\cU$-approximation} is defined dually.
For an object $X$ of $\cT$, we denote by $\add X$ the full subcategory
formed by all direct factors of finite direct sums of copies of $X$.

Let $M$ be a basic silting object of $\cT$ and $X$ an indecomposable direct
summand of $M$. Denote by $M/X$ the object such that $M\cong X\oplus M/X$.
By definition, the {\em left mutation $\mu_X(M)$ of $M$ at $X$} is the silting object
$M/X\oplus Y$, where $Y$ is defined by a triangle
\[
\xymatrix{
X \ar[r] & E \ar[r] & Y \ar[r] & \Sigma X}
\]
and $X \to E$ is a minimal left $\add(M/X)$-approximation. It is not hard to show
that then $E \to Y$ is a minimal right $\add(M/X)$-approximation,
which implies that $Y$ is indecomposable. Indeed, let us recall the argument: Let $M'=M/X$. 
Since $M$ is silting, we have an exact sequence of $\End(M')$-modules
\[
\Hom(M',E) \to \Hom(M',Y) \to 0 \ko
\]
where $\Hom(M',E)$ is projective over $\End(M')$. Saying that
$E \to Y$ is a minimal right $\add(M/X)$-approximation is
equivalent to saying that $\Hom(M',E) \to \Hom(M',Y)$ is
a projective cover. If $Y$ was decomposable, the
morphism $\Hom(M',E) \to \Hom(M',Y)$ would therefore
split into a direct sum of two surjective morphisms and
this splitting would be induced by a splitting of the
morphism $E \to Y$ as a direct sum of two non
trivial morphisms $E' \to Y'$ and $E''\to Y''$.
But then $X$ would be decomposable, a contradiction.
The right mutation $\mu_X^-(M)$ is defined dually. The right mutation of $\mu_X(M)$ at $Y$ is
isomorphic to $M$.

\begin{example} \label{ex:siltingGinzburg}
Suppose that $\Ga$ is a Jacobi-finite Ginzburg algebra associated with
a finite quiver and a potential not containing cycles of length $\leq 2$.
Then $A=\Ga$ satisfies our assumptions and $\Ga$ is a basic silting
object in $\per(\Ga)$. Let $M$ be a silting object in
$\per(\Ga)$ and $\Ga'$ the derived endomorphism algebra of $M$.
Then the homologies $H^p(\Ga')$ are finite-dimensional and vanish in
degrees $p>0$.
Since $M$ generates $\per(\Ga)$, the $\Ga'$--$\Ga$--bimodule $M$ yields
an algebraic triangle equivalence $\D(\Ga') \iso \D(\Ga)$. Conversely,
if we start from a dg algebra $\Ga'$ whose homologies are finite-dimensional
and vanish in degrees $>0$ and from an algebraic triangle equivalence 
$\D(\Ga') \iso \D(\Ga)$, then the image $M$ of $\Ga'$
in $\per(\Ga)$ is a silting object. In any case, the
dg algebra $\Ga'$ is an exact bimodule $3$-Calabi-Yau and has its
homology concentrated in degrees $\leq 0$. By Van den Bergh's
theorem \cite{VdB15}, the dg algebra $\Ga'$ is again a Jacobi-finite
Ginzburg algebra (up to weak equivalence). In particular, for $M$
we can take the mutation $M'=\mu_X \Ga$, where $X = e_i \Ga$ for
a vertex $i$ of the quiver of $\Ga$. We define the associated
Ginzburg algebra $\Ga'$ to be the {\em left mutation of $\Ga$ at $i$}.
Notice that by construction, we have a canonical
derived equivalence from $\Ga'$ to $\Ga$.
In the same way, we can define the right mutation $\Ga''$ of $\Ga$ at $i$
using the right mutation $M''=\mu_X^-(\Ga)$ of $\Ga$ at $X$.
The right mutation $\Ga''$ turns out to be quasi-isomorphic to the
left mutation $\Ga'$. 
Indeed, by Theorem~\ref{Gamma=truncation}, these algebras are the
$\tau_{\leq 0}$--truncations of the derived endomorphism algebras of
the images $\pi(M')$ and $\pi(M'')$ in the cluster category $\cC_\Ga$.
Now we have $\pi(M')\cong \pi(M'')$ because they are the left resp.
right mutation in the sense of Iyama-Yoshino \cite{IyamaYoshino08}
of the cluster-tilting object $\pi(\Ga)$ at $\pi(X)$ and for cluster-tilting
objects in $2$-Calabi-Yau triangulated categories, right and left mutation 
coincide up to isomorphism.
\end{example}

Now let $A$ be a dg $k$-algebra whose homologies $H^p A$ 
vanish in all degrees $p>0$. An object of $\per(A)$ is
called {\em $2$-term} if it is the cone over a morphism
of $\add(A)$. We write $\twoper(A)$ for the full
subcategory of $\per(A)$ formed by the $2$-term objects.
We write $\2silt(A)$ for the set of isomorphism classes of
{\em 2-term silting objects}.
Our aim is to compare $\2silt(A)$ with $\2silt(H^0A)$.
Note that by our assumption on $A$, we have a canonical morphism
$A \to H^0A$ in the homotopy category of dg algebras. We write 
$F: \per(A) \to \per(H^0 A)$ for the derived tensor product over $A$ with
$H^0A$. Part a) of the following theorem is due to 
Br\"ustle--Yang \cite{BrustleYang13}.

\begin{prop} \label{prop:2silt}
\begin{itemize} 
\item[a)] The functor $F$ induces a bijection $\2silt(A) \iso \2silt(H^0 A)$.
\item[b)] Suppose that we have $H^{-1}(A)=0$. Then the functor
$F$ restricts to an equivalence
\[
\twoper(A) \iso \twoper(H^0 A).
\]
In particular, for each $2$-term object $T$, the functor $F$ induces
an isomorphism 
\[
\End_A(T) \iso \End_{H^0 A}(FT).
\]
\end{itemize}
\end{prop}

\begin{proof} 
Part a) is Proposition~A.3 of \cite{BrustleYang13}. For part b), using the assumption
and the 5-lemma, we check successively that $F$ induces the following bijections:
\begin{itemize}
\item[1)] For $P, Q\in \add(A)$ and $p\geq -1$
\[
\Hom(P,\Sigma^p Q) \iso \Hom(FP, \Sigma^p FQ).
\]
\item[2)] For $P\in \add(A)$, $M\in \twoper(A)$ and $p\in \{-1,0\}$
\[
\Hom(P, \Sigma^p M) \iso \Hom(FP, \Sigma^p FM).
\]
\item[3)] For $M, M'\in \twoper(A)$
\[
\Hom(M,M') \iso \Hom(FM, FM').
\]
\end{itemize}
\end{proof}

Now let $A$ be a pseudocompact dg algebra in $\PCAlgc(l)$ strictly concentrated in
degrees $\leq 0$. Let $e$ be an idempotent of $H^0 A$ and $A'$ the
derived endomorphism algebra of the image of $A$ in the Verdier
quotient of $\per(A)$ by the thick subcategory generated by $eA$.
Then $A'$ is concentrated in degrees $\leq 0$ and
we have a canonical morphism $A \to A'$ in the homotopy category
of $\PCAlgc(l)$. If $A$ is of the form $(\wh{T}_l(V), d)$ for a pseudocompact
$l$-bimodule $V$ concentrated in degrees $\leq 0$, 
where $\wh{T}_l(V)$ is the completed tensor algebra,
then $A'$ is quasi-isomorphic to the quotient of $A$ by the two-sided 
closed ideal generated by $e$ (cf. \cite{BraunChuangLazarev18}).
Put $A_0=H^0 A'$ so that we have a canonical morphism $p: A \to A_0$. 
Let $B$ and $B_0$ be pseudocompact dg algebras in $\PCAlgc(l)$,
$X\in\D(A^{op}\ten B)$ such that $X_B$ is perfect 
and $Q\in \D(B^{op}\ten B_0)$ such that $Q_{B_0}$ is perfect.
\begin{prop} \label{prop:induced-bimodule} 
Suppose that $eX\lten_B Q$ vanishes and that the object 
$X\lten_B Q$ of $\D(B_0)$ has
no selfextensions in degrees $p<0$. Then there is an object
$Y$ of $\D(A_0^{op}\ten B_0)$, unique up to isomorphism, such that
we have an isomorphism
\[
X \lten_B Q \iso A_0 \lten_{A_0} Y
\]
in $\D(A^{op} \ten B_0)$. Thus, the square
\[
\xymatrix{\D(A) \ar[d]_{A_0} \ar[r]^{X} & 
								\D(B) \ar[d]^{Q}\\
\D(A_0) \ar[r]_{Y} & \D(B_0)}
\]
is commutative up to isomorphism, where we write dg bimodules instead of derived tensor
products by dg bimodules.
\end{prop}

\begin{remark} In our applications in this article, the idempotent $e$ will be $0$.
We state and prove the proposition in the general case because it provides
an alternative approach to the problem of relating the tilting theory of maximal 
modification algebras \cite{Wemyss18} to that of the associated contraction algebras as
treated by August in \cite{August18b}. Let $R$ be a complete local cDV
singularity and $M$ a maximal basic rigid object in the category of Cohen-Macaulay
modules over $R$ containing $R$ as a direct summand.
We can take $A=\End_R(M)$ and $e$ the idempotent corresponding to
the projection on $R$. Then $A_0=H^0 A'$ is isomorphic
to the stable endomorphism algebra of $M$, i.e. the contraction
algebra associated with $M$. Let $N$ be another maximal basic rigid object containing
$R$ as a direct summand, $B$ its endomorphism algebra and $B_0$ the
associated contraction algebra. Then $X=\Hom_R(N,M)$ yields a derived
equivalence $?\lten_A X: \D(A) \iso \D(B)$ taking $eA$ to $eB$. Moreover,
the complex $X\lten_B B_0$ is a 2-term silting object of $\per(B_0)$ 
(as it follows from silting reduction \cite{AiharaIyama12}
combined with part a) of Proposition~\ref{prop:2silt}) and 
hence a tilting object since $B_0$ is symmetric. Thus, the hypotheses
of the proposition hold and there is a canonical two-sided tilting complex
$Y$ in $\D(A_0^{op}\ten B_0)$. Clearly, the construction is compatible
with compositions via derived tensor products. 
\end{remark}

\begin{proof}[Proof of the Proposition] Put $U=X\lten_B Q$ viewed as an object in $\D(A^{op}\ten B_0)$.
The morphism $A \to A'$ induces the Verdier quotient 
\[
\per(A) \to  \per(A')=\per(A)/\thick(eA) 
\]
and is therefore a dg quotient. By the universal property of the dg quotient,
there is an object $Z$ in $\D({A'}^{op}\ten B_0)$ unique up to isomorphism
such that the restriction of $Z$ along $A^{op}\ten B_0 \to {A'}^{op}\ten B_0$
is isomorphic to $U$ in $\D(A^{op}\ten B_0)$.
We have a canonical morphism in the homotopy category of dg algebras
\[
\xymatrix{
A' \ar[r] & \RHom_{B_0}(Z,Z).}
\]
Since by assumption $\RHom_{B_0}(Z,Z)$ is concentrated in degrees $\geq 0$
and $A'$ in degrees $\leq 0$, this morphism factors uniquely through a morphism
\[
A_0=H^0 A' \to \RHom_{B_0}(Z,Z)
\]
in the homotopy category of dg algebras. Let us show how to refine this
argument so as to obtain an object $Y$ of $\D(A_0^{op} \ten B_0)$
which restricts to $Z \in \D(A'^{op}\ten B_0)$. We may and will assume that
$A'$ is cofibrant and strictly concentrated in degrees $\leq 0$. We
may and will also assume that $Z$ is cofibrant as a dg $A'$-$B_0$-bimodule.
The left $A'$-module structure on $Z$ then yields a morphism of dg algebras
\[
A' \longrightarrow \Hom_{B_0}(Z,Z).
\]
Since $A'$ is strictly concentrated in  non positive degrees, it factors uniquely
through a morphism of dg algebras
\[
A' \longrightarrow \tau_{\leq 0} \Hom_{B_0}(Z,Z).
\]
Since $Z_{B_0}$ has no negative selfextensions, we have a surjective quasi-isomorphism
of dg algebras
\[
\tau_{\leq 0} \Hom_{B_0}(Z,Z) \longrightarrow H^0 \Hom_{B_0}(Z,Z).
\]
The composition
\[
A' \longrightarrow \tau_{\leq 0} \Hom_{B_0}(Z,Z)  \longrightarrow H^0 \Hom_{B_0}(Z,Z)
\]
uniquely factors through an algebra morphism $A_0=H^0(A') \to H^0\Hom_{B_0}(Z,Z)$. We
thus obtain a commutative square of dg algebra morphisms
\[
\xymatrix{
A' \ar[dd] \ar[r] \ar[dr]& \Hom_{B_0}(Z,Z) \\
                        & \tau_{\leq 0} \Hom_{B_0}(Z,Z) \ar[u] \ar[d] \\
A_0=H^0(A') \ar[r] & H^0 \Hom_{B_0}(Z,Z).
}
\]
We factor the morphism $A' \to A_0$ as the composition $A' \to \wt{A}_0 \to A_0$ 
of an acyclic fibration with a cofibration. 
We consider the diagram
\[
\xymatrix{
A' \ar[d] \ar[r] \ar[dr]& \Hom_{B_0}(Z,Z) \\
 \wt{A}_0 \ar[d] \ar@{.>}[r]     & \tau_{\leq 0} \Hom_{B_0}(Z,Z) \ar[u] \ar[d] \\
A_0=H^0(A') \ar[r] & H^0 \Hom_{B_0}(Z,Z).
}
\]
Here the morphism represented by a dotted arrow exists so that the
diagram becomes commutative because $A' \to \wt{A}_0$ is a cofibration
and $\tau_{\leq 0} \Hom_{B_0}(Z,Z) \to H^0 \Hom_{B_0}(Z,Z)$ is
an acyclic fibration. Thus, we obtain a structure of dg 
$\wt{A}_0$-$B_0$-bimodule on $Z$ wich restricts to the given
structure of dg $A'$-$B_0$-bimodule. Since we have the quasi-isomorphism
$\wt{A}_0 \to A_0$, we can find a bimodule $Y$ in $\D(A_0^{op}\ten B_0)$
unique up to isomorphism and which restricts (up to isomorphism)
to $Z$ in $\D({A'}^{op}\ten B_0)$. 

Now suppose we have a second object $Y'$ in $\D(A_0^{op}\ten B_0)$
which becomes isomorphic to $Z$ in $\D({A'}^{op}\ten H^0 B)$. We have a chain
of isomorphisms
\begin{align*}
\Hom_{\D(A'^{op}\ten B_0)}(Y,Y') 
&=  \Hom_{\D(A'^e)}(A', \RHom_{B_0}(Y,Y')) \\
&= \Hom_{\D(A'^e)}(H^0 A', H^0 \RHom_{B_0}(Y,Y')) \\
&= \Hom_{H^0(A')^e}(H^0 A', H^0 \RHom_{B_0}(Y,Y'))\\
&= \Hom_{\D(A_0^e)}(A_0, H^0 \RHom_{B_0}(Y,Y'))\\
&=\Hom_{\D(A_0^e)}(A_0, \RHom_{B_0}(Y,Y'))\\
&=\Hom_{\D(A_0^{op}\ten B_0)}(Y,Y').
\end{align*}
Clearly, the composition
\[
\Hom_{\D(A_0^{op}\ten B_0)}(Y,Y') \iso \Hom_{\D(A'^{op}\ten B_0)}(Y,Y') 
\]
of these isomorphisms is given by the restriction of scalars
functor
\[
\D(A^0 \ten H^0 B) \to \D(A'\ten H^0 B).
\]
Now any restriction of scalars functor reflects isomorphisms because it
is compatible with the forgetful
functors to the derived category of vector spaces. Thus, isomorphisms
are preserved. This shows the uniqueness. 
\end{proof}

Let $B$ be a dg $k$-algebra whose homologies are
finite-dimensional and vanish in degrees $>0$. Let $C$ be a finite-dimensional
basic $k$-algebra (i.e. $C$ is basic as a right module over itself)
and let $Z$ be an object of $\D(C^{op}\ten H^0 B)$ such that
\[
?\lten_C Z: \D(C) \to \D(H^0 B)
\]
is an equivalence. Notice that $Z_{H^0 B}$ is
a tilting object in $\per(H^0 B)$ and in particular a silting object, which
is basic by our assumption on $C$.

\begin{theorem} \label{thm:liftingEquivalences}
Assume that $Z_{H^0 B}$ is a 2-term silting object and we have $H^{-1}(B)=0$.
Then there is a dg algebra $A$ whose homologies $H^p A$ 
are finite-dimensional and whose components vanish in degrees $p>0$, a derived equivalence 
$?\lten_A X: \D(A) \iso \D(B)$, an isomorphism of algebras $\phi: H^0 A \iso C$ 
and an isomorphism
\[
\mbox{}_\phi Z \iso X\lten_B H^0 B
\]
in $\D(A^{op}\ten H^0 B)$, where the left $A$-module structure on $\mbox{}_\phi Z$
is defined via the composition $A\to H^0 A \to C$. In particular, we have a diagram, 
commutative up to isomorphism 
\[
\xymatrix{
 & \D(A) \ar[dl]_{H^0 A} \ar[r]^X & \D(B) \ar[d]^{H^0 B} \\
\D(H^0 A) \ar[r]_{\mbox{}_\phi C } & \D(C) \ar[r]_-{Z} & \D(H^0 B)\ko}
\]
where we write dg bimodules instead of derived tensor products by dg bimodules.
\end{theorem}

\begin{proof} By part a) of Proposition~\ref{prop:2silt}, there is a
2-term silting object $M$ of $\per(B)$ such that $M\lten_B H^0 B$ is isomorphic
to $Z_{H^0 B}$. Since $M$ is silting, its derived endomorphism algebra
has its homology concentrated in non positive degrees and we define
\[
A=\tau_{\leq 0}\RHom(M,M).
\]
We let $X\in \D(A^{op}\ten B)$ be the dg bimodule given by 
$M$ with its canonical left $A$-action. Since $X\lten_B H^0 B$ is
isomorphic to the tilting object $Z_{H^0 B}$, it has no self-extensions
in degree $<0$. Therefore, Proposition~\ref{prop:induced-bimodule}
yields an object $Y$ of $\D(H^0(A)^{op}\ten H^0(B))$ and an isomorphism
\[
\psi: X \lten_B H^0 B \iso Y|_{A^{op}\ten H^0 B}
\]
in $\D(A^{op} \ten H^0 B)$, where the left $A$-module structure on $Y$ comes
from the canonical morphism $A \to H^0(A)$. By part b) of Proposition~\ref{prop:2silt},
we have an isomorphism $\End(M) \iso \End(M\lten_B H^0 B)$. By
construction, we have an isomorphism $H^0 A \iso \End(M)$ or equivalently
$H^0 A \iso \End(X_B)$. Thus, the composition of $?\lten_B H^0 B$ with
$?\lten_A X$ induces an isomorphism
\[
H^0 A \iso \End_{\D(H^0 B)}(X \lten_B H^0 B).
\]
Via $\psi$, we get an isomorphism
\[
H^0 A \iso \End_{\D(H^0 B)}(Y)
\]
given by the left action of $H^0 A$ on $Y$. We choose an isomorphism
$Y_{H^0 B} \iso Z_{H^0 B}$ in $\D(H^0 B)$ and define $\phi: H^0 A \iso C$ so
as to make the following square commutative
\[
\xymatrix{H^0 A \ar[r] \ar[d]_\phi & \End(Y_{H^0 B}) \ar[d] \\
C \ar[r] & \End(Z_{H^0 B}).}
\]
By Lemma~\ref{lemma:Toda-condition}, the chosen isomorphism
$Y_{H^0 B} \iso Z_{H^0 B}$ lifts to an isomorphism $Y \iso \mbox{}_\phi Z$
in $\D(H^0(A)^{op}\ten H^0(B))$. Whence a composed isomorphism
\[
X\lten_B H^0 B \iso Y|_{A^{op}\ten H^0 B} \iso \mbox{}_\phi Z|_{A^{op}\ten H^0 B}
\]
\end{proof}

\subsection{Cyclic homology and preservation of the canonical class}
\label{subsec:cyclic}
Let $k$ be a field of characteristic $0$ and $l$ a finite product of copies of $k$.
Let $V$ be a pseudocompact $l$-bimodule and $d$ a continous differential
on the completed tensor algebra $\wh{T}_l(V)$. Put $A=(\wh{T}_l(V),d)$.
We define $\Omega_l A$ by the short exact sequence
\[
\xymatrix{
0 \ar[r] & \Omega^1 A \ar[r] & A\ten_l A \ar[r]^-{\mu} & A \ar[r] & 0},
\]
where $\mu$ is the multiplication of $A$. Then the morphism
\[
A \ten_l V \ten_l A \to \Omega^1 A
\]
taking $a\ten v\ten b$ to $av\ten b - a \ten vb$ is an isomorphism of graded
$l$-bimodules, cf. Example 3.10 of \cite{Quillen89}. We can describe the
induced differential on $A \ten_l V \ten_l A$ as follows (cf. Proposition 3.7
of \cite{KV09}): Let $D: A \to A\ten_l V\ten_l A$
be the unique continuous bimodule derivation which restricts to the map
$v\mapsto 1\ten v\ten 1$ on $V$. We have
\[
D(v_1 \ldots v_n) = 1\ten v_1 \ten (v_2 \ldots v_n) + 
\sum_{i=2}^{n-1} v_1 \ldots v_{i-1}\ten v_i \ten v_{i+1} \ldots v_n +
(v_1 \ldots v_{n-1}) \ten v_n \ten 1.
\]
Then the induced differential on $A\ten_l V \ten_l A$ sends $a\ten v \ten b$ to
\[
(-1)^{|a|} a D(dv) b + (da)\ten v \ten b + (-1)^{(|v|+|a|)} a \ten v \ten (db).
\]

For an $l$-bimodule $M$, we write $M_l$ for the coinvariant module $M/[l,M]$.
For an $A$-bimodule $M$, we let $M_\natural$ be the coinvariant module $M/[A,M]$.
We have an isomorphism of graded $A$-modules
\[
(A\ten_l V \ten_l A)_\natural \iso (V\ten_l A)_l
\]
taking $a\ten v\ten b$ to $(-1)^{|a|(|v|+|b|)} v\ten ba$. The induced differential
on the right hand side is given as follows: If 
$D(dv)=\sum_i a_i \ten v_i \ten b_i$, then
\[
d(v\ten a) = (-1)^{|v|} v\ten (da) + \sum_i (-1)^{|a_i|(|v_i|+|b_i|+|a|)} v_i \ten b_i a a_i.
\]

Following section~3 of \cite{Quillen89}
we define morphisms of complexes 
\[
\del_1 : (V\ten_l A)_l \to A_l \quad\mbox{and}\quad
\del_0 : A_l \to (V\ten_l A)_l
\]
as follows: $\del_1$ sends
$v\ten a$ to $va - (-1)^{|v||a|} av$ and $\del_0$ sends $v_1 \ldots v_n$
to 
\[
\sum_i \pm v_i \ten v_{i+1} \ldots v_n v_1 \ldots v_{i-1},
\]
where the sign is determined by the Koszul sign rule. We then have
$\del_0 \del_1 =0 =\del_1 \del_0$.

The (continuous) Hochschild homology of $A$ is computed by the total complex of
\[
\xymatrix{
(V\ten_l A)_l \ar[r]^-{\del_1} & A_l}
\]
and the (continuous) cyclic homology of $A$ is computed by the product total complex of
\[
\xymatrix{
\ldots \ar[r] & (V\ten_l A)_l \ar[r]^-{\del_1} & A_l \ar[r]^-{\del_0} & (V\ten_l A)_l \ar[r]^-{\del_1} & A_l}.
\]
Since $k$ is of characteristic $0$, the morphism $A_l \to A/([A,A]+l)$ induces
a quasi-isomorphism from this complex to $A/([A,A]+l)$. The ISB-sequence 
\[
\xymatrix{
\ldots \ar[r] & \HoH_n \ar[r]^I & \HC_n \ar[r]^S & \HC_{n-2} \ar[r]^B & \HoH_{n-1} \ar[r]
& \ldots}
\]
is induced by the following sequence
\[
\xymatrix{
\ldots \ar[r] & 0 \ar[d] \ar[r] & 0 \ar[d] \ar[r] & (V\ten_l A)_l \ar@{=}[d] \ar[r]^-{\del_1} & 
	A_l \ar@{=}[d]\\
\ldots \ar[r] & (V\ten_l A)_l \ar@{=}[d] \ar[r]^-{\del_1} & A_l \ar@{=}[d] \ar[r]^-{\del_0} 
& (V\ten_l A)_l \ar[d] \ar[r]^-{\del_1} & A_l \ar[d] \\
\ldots \ar[r] & (V\ten_l A)_l \ar[d]\ar[r]^-{\del_1} & A_l \ar[d]^{\del_0} \ar[r] & 0\ar[d] \ar[r] & 
0 \ar[d] \\
\ldots \ar[r] & 0 \ar[r] & (V\ten_l A)_l \ar[r]^-{\del_1} & A_l \ar[r] & 0
}
\]
Notice that the first three rows form a short exact sequence and that
the composition of the last vertical morphism with the second last
vertical morphism is only homotopic to zero.

Now let $Q$ be a finite quiver, $l$ the product over the $ke_i$, where $i$ runs
through the vertices of $Q$, $w$ a potential on $Q$ and $A=\Ga$ the
associated complete Ginzburg algebra with generators the arrows $\alpha$ of
$Q$ in degree $0$, the reversed arrows $\alpha^*$ in degree $-1$ and
the loops $t_i$ in degree $-2$. Then $A$ is of the form $(\wh{T}_l(V), d)$, where
$V$ is the $l$-bimodule with basis given by the arrows $\alpha$, $\alpha^*$ and $t_i$.
Let $t$ be the sum of the $t_i$. By definition,
we have
\[
d(t)=\sum_\alpha [\alpha, \alpha^*].
\]
Thus, $t$ defines an element in $\HC_{2}(A)$. 

\begin{lemma} The image of the class of $t$ under $S: \HC_2(A) \to \HC_0(A)$
is the canonical class $[w]$, i.e. the image of $w$ under the projection
$\HC_0(T_l V) \to \HC_0(A)$.
\end{lemma}

\begin{proof}
We compute $S(t)$ using the
above description of $S$. We need to lift $t$ to an element of the
total complex computing cyclic homology. We have
\[
d(t) = \del_1(\sum_\alpha \alpha\ten \alpha^*).
\]
We have
\[
d(\sum_\alpha \alpha\ten \alpha^*) = \sum_\alpha \alpha\ten d(\alpha^*)
=\sum_\alpha \alpha\ten D_\alpha(w).
\]
Thus, we have
\[
d(\sum_\alpha \alpha\ten \alpha^*) = \del_0(w)
\]
and $S(t)$ is the image of $w$ in $H^0(A/([A,A]+l)=\HC_0(\Ga)$.
Notice that $BS(t)=B(w)$ is indeed a boundary in the Hochschild
complex: It is the differential of
\[
\sum_\alpha \alpha\ten \alpha^* -t.
\]
\end{proof}

\begin{corollary} \label{cor:adapted-QP} Let $\Ga'=\Ga(Q',w')$ be a Ginzburg algebra
and $A$ a pseudo-compact dg algebra in $\PCAlgc(l)$
concentrated in degrees $\leq 0$.
Let $X$ be a dg $A$-$\Ga'$-bimodule such that $?\lten_A X : \D(A) \to \D(\Ga')$ 
is an equivalence. Then there is a quiver with potential $(Q'',w'')$
and a weak equivalence $s: \Ga(Q'',w'') \to A$ such that for the
restriction $_sX$ along $s$, the isomorphism $\HC_0(_sX)$ takes
the class $[w'']$ to $[w']$.
\end{corollary}

\begin{proof} We know that the class $[t'] \in \HC_2(\Ga')$ is non degenerate in the
sense that $B[t']\in \HoH_3(\Ga')$ defines an isomorphism 
$\Sigma^3 \Theta_{\Ga'} \iso \Ga'$
in $\D(\Ga'^e)$, where $\Theta_{\Ga'}$ is the inverse dualizing complex. 
Thus the image $\tau$
of $[t']$ under $\HC_2(X)^{-1}$ is a non degenerate element of $\HC_2(A)$. 
The proof of Theorem 10.2.2 in \cite{VdB15} then shows that there is a quiver
$Q''$, a potential $w''$ and a weak equivalence $s: \Ga(Q'',w'') \to A$
which takes $[t'']$ to $\tau$. Thus, the composition 
$\HC_2(X)\circ \HC_2(s)=\HC_2(_sX)$ takes $[t'']$ to $[t']$ and the 
isomorphism $\HC_0(_sX)$ takes $[w'']=S[t']$ to $[w']=S[t']$.
\end{proof}

\section{CY tilted algebras and singularities}\label{sec:sing}
\subsection{Basics on Hochschild cohomology}
Let $k$ be a commutative ring and $A$ be a unital  $k$-algebra projective
over $k$. Denote by $\ol{A}$ the 
quotient $A/k\cdot 1$.  Define the \emph{normalized bar complex} associated to $A$ to be the complex $B_kA:= A\ot_k T\Sg \ol{A}\ot_k A$ with differential $\sum_{i= 0}^{n-1} (-1)^i b_i: A\ot \ol{A}^{\ot n-1}\ot A\to A\ot \ol{A}^{\ot n-2}\ot A$
\[
b_i(a_0,\ldots,a_{n})=(a_0,\ldots,a_ia_{i+1},\ldots,a_n).
\]
It is a projective bimodule resolution of $A$. 
Let $M$ be an $A$-bimodule. The \emph{Hochschild cochain complex} with coefficients in the bimodule $M$ is defined to be the complex 
$C^*(A,M):=\Hom_{A^e}(B_k(A),M)$ with differential
\[
\delta(f)= -(-1)^n f\circ b
\] for $f: A\ot \ol{A}^{\ot n} \ot A\to M$. The $i$-th \emph{Hochschild cohomology of the algebra $A$ with coefficients in the bimodule $M$} is defined to be 
$\HoH^i(A,M):=H^i(C^*(A,M), \delta)$.

Let $A$ be an augmented dg $k$-algebra. Denote by $\ol{A}$ the kernel of the augmentation.  Then the bar complex $B_kA$ is equipped with a second differential induced from the differential $d_A$ on $A$. Given an $A$-bimodule $M$, the Hochschild cochain complex  $C^*(A,M)$ is equipped with a second differential $d$ induced by $d_A$ and the internal differential $d_M$ on M. The $i$th 
\emph{Hochschild cohomology of the dg algebra $A$ with coefficients in 
the bimodule $M$} is defined to be $\HoH^i(A,M):=H^i(C^*(A,M), d+\delta)$.
 
It is well-known   that 
$\HoH^i(A,M)$ is isomorphic to $\Ext_{A^e}^i(A,M)$. When $A$ is a smooth commutative $k$-algebra, $\HoH^*(A,A)$ is isomorphic to the polyvector fields on 
$\Spec A$ by the Hochschild-Kostant-Rosenberg theorem. 
For non-smooth algebras, there exist different variants of Hochschild cohomology.

Let $A$ be an associative $k$-algebra projective over $k$.
Define the module of \emph{K\"ahler differentials} $\Omega_A$ to be the kernel of the multiplication map $\mu: A\ot A\to A$. Clearly, $\Omega_A$ inherits a bimodule structure from $A\ot A$. It is easy to show that $\Omega_A$ is generated as a bimodule by the elements of the form $xdy:= xy\ot 1- x\ot y$. The left and right module structure are given by
\[
a(xdy)=(ax)dy,~~~~(xdy)a=xd(ya)-xyda.
\]

Define the module of $n$-forms to be the $n$-fold tensor product
\[
\Omega^n_A:=\Omega_A\ot_A\Omega_A\ot\ldots\ot_A\Omega_A.
\] 
Using the above identities, one can check that $\Omega^n_A$ is generated as a bimodule by the elements of the form $a_0 da_1 da_2\ldots da_n$. There is an isomorphism of bimodules $\Omega^n_A\cong A\ot_k \ol{A}^{\ot n}$ defined by 
\[
a_0 da_1 da_2\ldots da_n\mapsto a_0\ot a_1\ot\ldots\ot a_n.
\]
Set $\Omega^0_A=A$ and $\Omega^1_A=\Omega_A$.
Write $\Omega^*_A$ for $\bigoplus_{n\geq 0}\Omega^n_A$. 
The bimodule structure on $\Omega_A$ naturally extends to an associative algebra structure on $\Omega^*_A$. The obvious differential
\[
D: a_0 da_1 da_2\ldots da_n\mapsto da_0 da_1 da_2\ldots da_n
\] 
makes $\Omega^*_A$ into a differential graded algebra.

For $m\in \ZZ$, 
\[
C^m(A,\Sg^n\Omega^n_A)=\Hom_k(\ol{A}^{\ot(n+m)}, A\ot_k \ol{A}^{\ot n}).
\]
Consider the chain maps 
$\theta_n: C^*(A,\Sigma^n\Omega^n_A)\to C^*(A,\Sigma^{n+1}\Omega_A^{n+1})$ between the Hochschild cochain complexes defined by $f\mapsto f\ot \id_{\Sg \ol{A}}$.

\begin{definition} \label{def-sing-HH}
Let $A$ be an associative $k$-algebra. Then the \emph{singular Hochschild cochain
complex} of $A$, denoted by $C_{sg}^* (A, A)$, is defined as the colimit of the inductive system in 
     the category of cochain complexes of $k$-modules,
     \[\xymatrix{
0 \ar[r] & C^* (A,A)  \ar[r]^-{\theta_0} & C^*(A,\Sigma\Omega^1_{A})  \ar[r]^-{\theta_1} &\ldots  \ar[r] & C^*(A,\Sigma^n\Omega_A^n) \ar[r]^-{\theta_n} & \ldots .}
\]
Namely, $C^*_{sg}(A, A) := \colim_n C^*(A, \Sigma^n\Omega^n_{A} )$. 
Its cohomology groups are denoted by $\HoH^*_{sg}(A, A)$.
\end{definition}

By construction, we have a natural chain morphism from $C^*(A,A)$ to $C^*_{sg}(A,A)$, which induces a natural morphism from $\HoH^*(A,A)$ to $\HoH_{sg}^*(A,A)$.

Let $A$ be a Noetherian $k$-algebra. Define $\D_{sg}(A)$ to be the Verdier quotient of $\D^b(A)$ by the subcategory $\per(A)$. We denote the extension group in $\D_{sg}(A)$ by $\ul{\Ext}^i_A(?,?)$. The singular Hochschild cohomology groups are related to the extension groups in $\D_{sg}(A^e)$. 
\begin{prop}(Theorem 3.6 \cite{Wang18})\label{Ext=singHH}
Let $A$ be a Noetherian $k$-algebra. Then there exists a natural isomorphism 
\[
\HoH^*_{sg}(A,A)\iso \ul{\Ext}^*_{A^e}(A,A),
\] such that the following diagram commutes
\[
\xymatrix{
\Ext^*_{A^e}(A,A)\ar[r] & \ul{\Ext}^*_{A^e}(A,A)\\
\HoH^*(A,A)\ar[u]^\cong\ar[r] & \HoH_{sg}^*(A,A)\ar[u]^\cong
}
\]
\end{prop}
From Wang's result, we see that the singular Hochschild cohomology admits a structure of graded commutative algebra.
Notice that this is not immediate from generalities about monoidal triangulated categories because
the singularities category of $A^e$ does not have any obvious monoidal structure.

\subsection{Hochschild cohomology of Gorenstein algebras}
A (not necessarily commutative) Noetherian ring $A$ is called \emph{Gorenstein} if it has finite injective dimension both as a left and right $A$-module. As in the 
commutative case, we denote by $\CM_A$ the category of maximal Cohen-Macaulay (left) $A$-modules and denote by $\ul{\CM}_A$ its stable category. Buchweitz proved 
that if $A$ is Gorenstein, then
$\ul{\CM}_A$ is equipped with a structure of triangulated category and 
$\ul{\CM}_A\cong \D_{sg}(A)$. 

We recall a fundamental result on extension groups in the stable category 
over Gorenstein rings due to Buchweitz.
\begin{prop}(Corollary 6.3.4 of \cite{Bu87})\label{stableExt}
Let $A$ be a Gorenstein ring and  let $X$ and $Y$ be objects in $\D^b(A)$. There exists a positive integer $m$ depending on $A$, $X$ and $Y$ such that the natural morphism $\Ext^i_A(X,Y)\to\ul{\Ext}^i_A(X,Y)$ is surjective for $i=m$ and is an isomorphism for $i>m$. 
\end{prop}


Combining Proposition \ref{Ext=singHH} and Proposition \ref{stableExt}, we obtain the following result.
\begin{corollary}\label{asym-HH}
Let $R$ be a commutative noetherian Gorenstein $k$-algebra.
If $R\ten R$ is noetherian, there exists a positive integer $m$ such that for $i>m$, the natural morphism
\[
\HoH^i(R,R)\to \HoH^i_{sg}(R,R)
\] 
is an isomorphism.
\end{corollary}
\begin{proof}
By Definition \ref{def-sing-HH}, we have a morphism $\HoH^i(R,R)\to \HoH^i_{sg}(R,R)$ for all $i$. 
In order to apply Proposition \ref{stableExt} to show that it is an isomorphism, we need
to check that $R\ten R$ is Noetherian.
This follows from Theorem~1.6 of \cite{TousiYassemi03}.
\end{proof}

A commutative local complete Gorenstein $k$-algebra $\hat{R}$ is called a hypersurface algebra if $\hat{R}\cong k[[x_1,\ldots,x_n]]/(g)$. 
We say that $\hat{R}$ is a hypersurface algebra with isolated singularities if 
$g$ has an isolated critical point.
\begin{theorem}(Theorem 3.2.7 \cite{GGRV})\label{HHvsHHsg}
Let $R=k[x_1,\ldots,x_n]/(g)$ be a hypersurface algebra with isolated singularities.
Denote by $M_g$ the Milnor algebra $k[x_1,\ldots,x_n]/(\frac{\partial g}{\partial x_1},\ldots,\frac{\partial g}{\partial x_n})$, and by $K_g$ and $T_g$ the kernel and cokernel of the endomorphism of $M_g$ defined by multiplication with $g$.
Then for $r\geq n$, there is an isomorphism of $R$-modules 
\[
\begin{cases}
\HoH^{r}(R,R)\cong T_g & r~is~even,\\
\HoH^{r}(R,R)\cong K_g & r~is~odd.
\end{cases}
\]
\end{theorem}
\begin{proof} The proof in \cite{GGRV} shows that in degrees
$r\geq n$, the Hochschild cohomology $\HoH^r(R,R)$ is isomorphic to the
homology in degree $r$ of the complex
\[
k[u]\ten K(R, \frac{\partial g}{\partial x_1},\ldots,\frac{\partial g}{\partial x_n}) ,
\]
where $u$ is of degree $2$ and
$K$ denotes the Koszul complex. Put $P=k[x_1, \ldots, x_n]$.
Since $R$ is quasi-isomorphic
to $K(P,g)$ and the $(\partial g)/(\partial x_i)$ form a regular sequence in $P$, 
the Koszul complex is quasi-isomorphic to $K(M_g,g)$.
\end{proof}
Note that $T_g$ is the \emph{Tyurina algebra} $k[x_1,\ldots,x_n]/(g,\frac{\partial g}{\partial x_1},\ldots,\frac{\partial g}{\partial x_n})$. Since  $K_g$ is the kernel of the multiplication map $g: M_g\to M_g$ it is a module over $T_g$. 

\begin{lemma} Let $A$ be a commutative $k$-algebra such that $A$ and $A^e$
are noetherian. Let $S\subset A$ be a multiplicative subset. If $M$ is a finitely
generated $A$-module and $L$ an $A$-module, we have a canonical isomorphism
\[
\RHom_{A^e}(L,M)_S \iso \RHom_{A_S^e}(L_S,M_S).
\]
\end{lemma}
\begin{proof} Since $L$ is finitely generated over $A$, it is finitely generated over
$A^e$. Since $A^e$ is noetherian, we have a projective resolution $P \to L$
with finitely generated components. This implies that we have isomorphisms
\[
\RHom_{A^e}(L,M)_S = \Hom_{A^e}(P,M)_S = \Hom_{A^e}(P,M_S).
\]
Since $A_S \ten_A P \ten_A A_S \to L_S$ is a projective resolution over $A_S^e$, 
we find
\[
\Hom_{A^e}(P,M_S) = \Hom_{A^e_S}(A_S\ten_A P \ten_A A_S, M_S) =
\RHom_{A^e_S}(L_S, M_S).
\]
\end{proof}

\begin{remark} \label{rk:HH}
In the setting of Theorem~\ref{HHvsHHsg},
assume that $g$ has isolated singularities and that the origin is a singular point
of the vanishing locus of $g$. If we denote by  $\fm$ 
the maximal ideal $(x_1,\ldots,x_n)\subset k[x_1,\ldots,x_n]$, then $g\in \fm$.
Denote by $M_{g,\fm}$, $T_{g,\fm}$ and $K_{g,\fm}$ the localizations of $M_g$, $T_g$ and $K_g$. It follows from the lemma that Theorem~\ref{HHvsHHsg} still
holds if one replaces $R$ by $R_\fm:=k[x_1,\ldots,x_n]_\fm/(g)$
and replaces $T_g$ and $K_g$ by $T_{g,\fm}$ and $K_{g,\fm}$.
\end{remark}

For a Noetherian $k$-algebra $A$, the derived category of singularities $\D_{sg}(A)$ is equipped with a canonical dg enhancement, obtained from its construction as
a Verdier quotient of two canonically enhanced triangulated categories
\cite{Keller99,Drinfeld04}. Instead of $\HoH_{sg}^*(A,A)$, one may also consider 
the Hochschild cohomology of the dg category $\D_{sg}(A)$, 
which we will denote by $\HoH^*(\D_{sg}(A))$. 

\begin{theorem}(Keller \cite{Keller18})\label{HHsg=HHDsg}
There is a canonical morphism of graded algebras
\[
\HoH^*_{sg}(A,A) \to \HoH^*(\D_{sg}(A)).
\]
It is invertible if $\D^b_{dg}(\mod A)$ is smooth.
\end{theorem}

Now we establish the main result of the subsection.
\begin{theorem}\label{HHDsg=Tyurina}
Let $\hat{R}=k[[x_1,\ldots,x_n]]/(g)$ be a hypersurface algebra with isolated singularity. Denote by $\hat{T}_g$ the Tyurina algebra of $g$.
Then there is an isomorphism of $k$-algebras
\[
\HoH^0(\D_{sg}(\hat{R}))\cong \hat{T}_g .
\]
Moreover, 
if $\hat{R}^\p=k[[x_1,\ldots,x_n]]/(g^\p)$ is another hypersurface algebra with isolated singularity such that 
$\D_{sg}(\hat{R}^\p)$ is quasi-equivalent with $\D_{sg}(\hat{R})$ as dg-categories, then $\hat{R}^\p$ is isomorphic to $\hat{R}$.
\end{theorem}
\begin{proof}
Because $g$ has an isolated critical point, we may assume that $g$ is a polynomial without loss of generality.
 Denote by $R$ the algebra $k[x_1,\ldots,x_n]_\fm/(g)$. Notice that $R$ has
an isolated singularity at the origin and that its completion identifies with $\hat{R}$.
By Theorem 5.7 of \cite{Dyc09}, the triangulated category
$\D_{sg}(\hat{R})$ is the Karoubi envelope of 
$\D_{sg}(R)$. Therefore, the two dg categories have equivalent derived
categories and there is a natural isomorphism 
$\HoH^*(\D_{sg}(R))\iso\HoH^*(\D_{sg}(\hat{R}))$. 
Orlov proved in \cite{Orlov04} that $\D_{sg}(R)$ is triangle equivalent with the homotopy category of matrix factorizations $\MF(k[x_1,\ldots,x_n]_\fm, g)$. The triangle equivalence is lifted to an equivalence of dg categories by the work of Blanc, Robalo, To\"en and Vezzosi \cite{BlancRobaloToenVezzosi16}. Therefore,
the dg category $\D_{sg}(R)$ is $2$-periodic and so is its Hochschild
cohomology. So there exists a natural isomorphism of $\hat{R}$-modules
\[
\HoH^0(\D_{sg}(R))\iso \HoH^{2r}(\D_{sg}(R))
\]
for all $r\in \ZZ$. 
By Theorem~B of
 \cite{ElaginLuntsSchnuerer18}, the bounded dg derived category
 $\D^b_{dg}(\mod R)$ is smooth. Thus, by  Theorem \ref{HHsg=HHDsg}, we have
\[
\HoH^{2r}(\D_{sg}(R)) \iso \HoH^{2r}_{sg}(R,R).
\]
By Corollary \ref{asym-HH}, for $r\gg 0$, we have
\[
\HoH^{2r}(R,R) \iso \HoH^{2r}_{sg}(R,R) 
\]
and by Theorem \ref{HHvsHHsg} and Remark \ref{rk:HH}, we have
\[
\HoH^{2r}(R,R) \iso T_g.
\]
Because $g$ has an isolated critical point, there is an isomorphism
\[
\hat{T}_g\cong T_g.
\]
Then the first claim follows. The second claim folllows from the formal version of the Mather-Yau theorem (see Theorem 1.1 \cite{GP}).
\end{proof}

\subsection{Classification of 3-dimensional smooth flops}
Let $(\hY,\hf,R)$ be a 3-dimensional formal flopping contraction with $\Ex(f)=\bigcup_{i=1}^t C_i$, and let $A=\End_R(\bigoplus_{i=1}^tN_i\op R)$ be the NCCR associated to it. We have associated to it an exact 3CY algebra: the derived deformation algebra $\Ga$ of the semi-simple collection $\cO_{C_1},\ldots,\cO_{C_t}$, and the cluster category $\cC_\Ga$. There are two relaxations of the above context.

First we take $Y$ a CY 3-fold with a semi-simple collection of rational curves $\cO_{C_1},\ldots,\cO_{C_t}$. We may still define the derived deformation algebra $\Ga$ and the cluster category $\cC_\Ga$. However, in general $\Ga$ is only  bimodule CY and $\cC_\Ga$ may not be Hom-finite. If we assume that the $Y$ is projective then $\Ga$ will be exact. However, most $\wh{Y}$ are not expected to have CY compactifications.
The second relaxation is to take a 3-dimensional hypersurface ring $R$ with isolated singularities. Associate to it is the derived category of singularity $\D_{sg}(R)$. This is a Hom-finite CY category. One may ask when is it possible to express it as a cluster category of a CY algebra.

In the case of flopping contractions, these two relaxations are related by Theorem \ref{VdBTh}. By \cite{VdB04}, $\wh{Y}$ admits a tilting bundle, $\Ga$ is exact and has finite dimensional cohomology. 
By Corollary \ref{nccr-vs-contraction}, Theorem \ref{VdBpotential} and \ref{VdBcomplete=>exact}, $\Ga$ is weakly equivalent to a Ginzburg algebra $\fD(Q,w)$ with $t$ nodes.  On the other hand  if $R$ admits a NCCR $A$ then $\D_{sg}(R)$ is equivalent to $\cC_\Ga$ for the CY algebra $\Ga$. By Theorem \ref{VdBTh}, $A$ admits a minimal model $(\wh{T}_{\wh{l}}V,d)$ with the dual space of $\Sigma^{-1}\Ext^{\geq 1}_A(\bigoplus_{i=0}^t S_i, \bigoplus_{i=0}^t S_i)$ (see section 4 of \cite{ThV}). Since $A$ is Calabi-Yau, $d$ can be derived from a potential, i.e. $(\wh{T}_{\wh{l}}V,d)$ is a Ginzburg algebra (see  \cite[Proposition 1.2]{ThV}).  By the derived equivalence theorem \ref{VdBeq} of Van den Bergh,  the derived deformation algebra of the semi-simple colleciton $\cO_{C_1},\ldots,\cO_{C_t}$ is isomorphic to the quotient $\wh{T}_{\ol{l}}V/\wh{T}_{\ol{l}}V e_0 \wh{T}_{\ol{l}} V$ in Theorem \ref{VdBTh}.  It is natural to expect that the deformation theory of the exceptional curves and the singularity theory of $R$ should determine each other   since both are governed by the CY algebra $\Ga$.

Recall that the CY tilted algebra 
$\Lm:=H^0\Ga$ is isomorphic to $\End_{\cC_\Ga}(\pi(\Ga))$ (c.f. Theorem \ref{Amiot}).  Donovan and Wemyss conjectured that $\Lm$ alone can already determine the analytic type of $R$.

\begin{conjecture}\label{DWconj}
Let $(\hY,\hf,R)$ and $(\hY^\p,\hf^\p,R^\p)$ be two 3-dimensional simple formal flopping contractions with associated CY tilted algebras $\Lm$ and $\Lm^\p$. 
Then the following are equivalent
\begin{enumerate}
\item[$(1)$] $R$ is isomorphic to $R^\p$.
\item[$(2)$] $\Lm$ is isomorphic to $\Lm^\p$.
\end{enumerate}
\end{conjecture}

Donovan and Wemyss have extended this conjecture to the case of
not necessarily simple formal flopping contractions by replacing (2) with
\begin{quote}
(2') $\Lm$ is derived equivalent to $\Lm^\p$,
\end{quote}
cf. Conjecture~1.3 in \cite{August18a}. In this situation, the implication from
(1) to (2) is known to be true by iterating a construction of Dugas \cite{Dugas15}.
The implication from (2') to (1) is one of the main open problems in the homological minimal model program for 3-folds. 
In this section, we will prove a slightly weaker version of this implication.


The exactness of $\Ga$ poses a strong constraint on $\Lm$, i.e. the relations of $\Lm$ can be written as cyclic derivative of a potential $w$ by Theorem \ref{VdBpotential}. If we fix the exact CY structure then $w$ is uniquely determined up to right equivalences (Proposition \ref{canonicalclass}). 

\begin{theorem}\label{weakDW}
Let $(\hY, \hf, R)$ and $(\hY^\p, \hf^\p,R^\p)$ be formal flopping contractions. Given exact CY structures $\eta$ and $\eta^\p$ on $\hY$ and $\hY^\p$ respectively, denote by $(\Lm,[w])$ and $(\Lm^\p,[w^\p])$ the associated CY tilted algebras and the canonical classes of potentials. If there exists a derived equivalence
from $\Lm$ to $\Lm^\p$ given by a bimodule complex $Z$ such that the
induced map $\HoH_0(Z)$ takes $[w]$ to $[w^\p]$, then $R$ is isomorphic to $R^\p$.
\end{theorem}

The proof of this theorem will take up the rest of the section. Here we highlight major components of the proof.
First we prove that the cluster category $\cC_\Ga$ is dg equivalent to the $\D_{sg}(R)$ with its canonical $\ZZ$-graded dg structure (Lemma \ref{eqdgenhancement}). The second step is to establish the result that analytic type of isolated hypersurface singularity (with fixed embedded dimension) is determined by its $\ZZ$-graded dg category of singularities (Theorem \ref{HHDsg=Tyurina}). In the last step, we prove that $\Ga$ can be reconstructed from the CY tilted algebra $\Lm$ together with the class $[w]\in \HoH_0(\Lm)$ represented by the potential.

Let $(R,\fm)$ be a complete commutative Noetherian local Gorenstein $k$-algebra of Krull dimension $n$ with isolated singularity and with residue field $k$. Suppose that $R$ admits a NCCR, then $\D_{sg}(R)$ has another dg model via 
the triangle equivalence $\D_{sg}(R)\simeq\cC_\Ga$ (c.f. Theorem \ref{VdBTh}). We first prove that these two models are dg quasi-equivalent. 

\begin{lemma}\label{eqdgenhancement} In the homotopy category of
dg categories, there is an isomorphism between 
$\cC_\Gamma=\per(\Ga)/\D_{fd}(\Ga)$ and the category of singularities 
$\D_{sg}(R)=\D^b(R)/K^b(\proj (R))$, both equipped with their canonical dg enhancements.
\end{lemma}

\begin{proof} Let $\cA$ and $\cB$ be two pretriangulated dg categories.
We call a triangle functor $F: H^0(\cA) \to H^0(\cB)$ algebraic if there
is a dg $\cA$-$\cB$-bimodule $X$ such that we have a square of triangle functors,
commutative up to isomorphism
\[
\xymatrix{
H^0(\cA) \ar[r]^F \ar[d] & H^0(\cB) \ar[d] \\
\D(\cA) \ar[r]_{?\lten_\cA X} & \D(\cB),}
\]
where the vertical arrows are induced by the Yoneda functors. We know
from \cite{Toen07} that morphisms
$\cA \to \cB$ in the homotopy category of dg categories are in bijection with
isomorphism classes of right quasi-representable $\cA$-$\cB$-bimodules
in the derived category of bimodules. Thus, it suffices to show
that the triangle equivalence $\cC_\Ga \iso \D_{sg}(R)$ is algebraic. 
We use the notation of subsection \ref{ss:NCCR} and put
$N=N_1\oplus \ldots \oplus N_t$. Let $\cF$ denote the
thick triangulated subcategory of $\per(A)$ generated by the simples
$S_1$, \ldots, $S_t$. Let us recall from Proposition~3 of
\cite{Palu09} that we have a diagram of triangle functors, 
commutative up to isomorphism and whose rows and columns
are exact sequences of triangulated categories
\[
\xymatrix{
            &                                 &  0                             & 0                         &   \\
0 \ar[r] & \cF \ar[r] \ar@{=}[d] & \per(A)/\per(R)\ar[u] \ar[r] & \ul{\CM}_R \ar[u] \ar[r] & 0 \\
0 \ar[r] & \cF \ar[r] & \per(A) \ar[u] \ar[r] & \D^b(\CM_R) \ar[u] \ar[r] & 0 \\
           &                & \per(R) \ar[u] \ar@{=}[r] & \per(R) \ar[u]                    &   \\
           &                &  0\ar[u]                          & 0 \ar[u]}
\]
Here the category $\D^b(\CM_R)$ is the bounded derived category of
the exact category $\CM_R$, the functor $\per(R) \to \D^b(\CM_R)$ is
induced by the inclusion $\proj(R) \to \CM_R$, 
the functor $\per(A) \to \D^b(\CM_R)$ is induced by $?\lten_A(R\oplus N)$
and $\per(R) \to \per(A)$ is induced by the inclusion 
\[
\add(R) \to \add(R\oplus N) = \proj(A) \ko
\]
where the categories $\add(R)$ and $\add(R\oplus N)$ are 
full subcategories of $\CM_R$
and the last equality denotes the equivalence given by the
functor $\Hom(R\oplus N,?)$. We endow $\per(A)/\per(R)$ with
the dg enhancement given by the dg quotient \cite{Keller99, Drinfeld04}.
It is then clear that the triangle functors of the middle row and of
the middle column are algebraic. Let us show that the functor
$\D^b(\CM_R) \to \ul{\CM}_R$ is algebraic. The canonical dg enhancement
of $\ul{\CM}_R$ is given by the triangle equivalence from the homotopy
category $H^0(\cA)$ of the dg category $\cA$ of 
acyclic complexes over $\proj (R)$ to $\ul{\CM}_R$
taking an acyclic complex $P$ to its zero cycles $Z^0(P)$. Let $\cB$
be the dg enhanced derived category $\D^b_{dg}(\CM_R)$. We define a
$\cB$-$\cA$-bimodule $X$ by putting
\[
X(P,M)=\Hom(P,\Sigma M)
\]
where $P$ is an acyclic complex of finitely generated projective $R$-modules
and $M$ a bounded complex over $\CM_R$. For $M\in\cB$, denote by
$M^\wedge$ the representable dg functor $\cB(?,M)$.  We have
$M^\wedge \ten_{\cB} X = X(?,M)$. If $M$ is concentrated in degree $0$,
we have canonical isomorphisms
\[
H^0 X(P,M) = \Hom_{\ul{\CM}_R}(Z^0(P), M) \quad\mbox{and}\quad
H^p X(P,M)= \Hom_{\ul{\CM}_R}(Z^0(P), \Sigma^p M) \ko
\]
which shows that $X(?,M)$ is quasi-representable by a complete resolution
of $M$. By d\'evissage, it follows that $X(?,M)$ is quasi-representable
for any bounded complex $M$ and one checks easily that the
(derived=non derived) tensor product with $X$ induces the canonical
triangle functor $\D^b(\CM_R) \to \ul{\CM}_R$. It follows that, at the level of
dg categories, $\ul{\CM}_R$ identifies with the dg quotient of $\D^b(\CM_R)$
by $\per(R)$. In other words, the canonical equivalence
\[
\D^b(\CM_R)/\per(R) \iso \ul{\CM}_R
\]
is algebraic. 
Therefore, the induced functor $\per(A)/\per(R) \to \ul{\CM}_R$
is algebraic. Thus, the whole diagram is made up of algebraic functors.
Now notice that the inclusion $\CM_R \subset \mod R$
induces  algebraic equivalences $\D^b(\CM_R) \iso \D^b(\mod R)$
so that we get algebraic equivalences
\[
\ul{\CM}_R \liso \D^b(\CM_R)/\per(R) \iso \D^b(\mod R)/\per(R)
\]
Thus, it will suffice to prove that the equivalence between 
the cluster category of $\Gamma$ and the stable category
$\ul{\CM}_R$ is algebraic. 

Now using the notations of Theorem~\ref{VdBTh} put 
$\wt{\Ga}=(\wh{T}_lV,d)$ so that we have a quasi-isomorphism
$\wt{\Ga} \iso A$. It induces an algebraic equivalence $\per(\wt{\Ga}) \iso \per(A)$.
This equivalence induces an algebraic equivalence $\tria(e_0\wt{\Ga}) \to \per(R)$,
where $\tria(e_0\tilde{\Ga})$ is the triangulated subcategory generated by
the $\wt{\Ga}$-module $e_0 \wt{\Ga}$.
The quotient map $\wt{\Ga} \to \Ga$ induces an algebraic triangle
functor $\per(\wt{\Ga}) \to \per(\Ga)$ and we know from Lemma~7.2
of \cite{KY18} that it is a localization with kernel $\tria(e_0\wt{\Ga})$.
We obtain a diagram of triangle functors, commutative up to isomorphism,
whose vertical arrows are equivalences and whose rows are exact.
\[
\xymatrix{
0 \ar[r] & \tria(e_0\wt{\Ga}) \ar[r] \ar[d] & \per(\wt{\Ga}) \ar[d] \ar[r] & 
					\per(\Ga) \ar[d] \ar[r] & 0 \\
0 \ar[r] & \per(R) \ar[r] & \per(A) \ar[r] & \per(A)/\per(R) \ar[r] & 0.}
\]
By passage to the dg quotient, the rightmost vertical arrow is an algebraic
triangle equivalence. By composing the algebraic inclusion
$\cF \to \per(A)/\per(R)$ with a quasi-inverse of the algebraic
equivalence $\per(\Ga) \iso \per(A)/\per(R)$, we obtain an
algebraic inclusion $\cF \to \per(\Ga)$ whose image identifies
with the thick subcategory of $\per(\Ga)$ generated by
the simple $H^0(\Ga)$-modules. This subcategory equals
$\D_{fg}(\Ga)$. Indeed, clearly it is contained in $\D_{fd}(\Ga)$ and
conversely, it contains the category $\mod H^0(\Ga)$ of finite-dimensional
$H^0(\Ga)$-modules since $H^0(\Ga)$ is finite-dimensional so that
every finite-dimensional $H^0(\Ga)$-module is a finite iterated
extension of simples. Since $\mod H^0(\Ga)$ is the heart of
a bounded $t$-structure on $\D_{fd}(\Ga)$, the image of $\cF$
equals $\D_{fd}(\Ga)$. This yields the first
row of the following diagram with exact rows
whose vertical arrows are equivalences.
\[
\xymatrix{
0 \ar[r] & \cF \ar[r] \ar@{=}[d] & \per(\Ga) \ar[d] \ar[r] & \cC_\Ga \ar[r] \ar[d] & 0 \\
0 \ar[r] & \cF \ar[r] & \per(A)/\per(R) \ar[r] & \ul{\CM}_R \ar[r] & 0.}
\]
Again by passage to the dg quotient, the rightmost vertical arrow is
an algebraic triangle equivalence.
\end{proof}

\begin{proof}[{\bf{Proof of Theorem \ref{weakDW}}}] Let $(Q,w)$ and
$(Q',w')$ be the quivers with potential constructed from the formal
flopping contractions and $\Ga$ and $\Ga'$ the associated Ginzburg dg algebras.
Let $\La=H^0(\Ga)$ and $\La'=H^0(\Ga')$ be the associated contraction algebras.
Recall that these algebras are symmetric so that tilting objects coincide
with silting objects in their derived categories. We will construct
a quiver with potential $(Q'', w'')$ with associated Ginzburg algebra
$\Ga''$, a dg $\Ga''$-$\Ga'$-bimodule $\tilde{Z}$ and an isomorphism
$\psi$ from $\La''= H^0(\Ga'')$ to  $C$ such that the square
\[
\xymatrix{
\Ga'' \ar[r]^{\tilde{Z}} \ar[d]_{\Lambda''} & \Ga' \ar[d]^{\Lambda'} \\
\Lambda'' \ar[r]_{\mbox{}_\psi Z} & \Lambda'}
\]
is commutative and the isomorphism $HC_0(\tilde{Z})$ takes
the class $[w'']$ to $[w']$. Here we write dg bimodules instead of derived tensor products, algebras
instead of their derived categories and the top and bottom arrows are equivalences. 
Notice that by \cite{Keller98}, Hochschild homology is functorial with respect to right perfect
dg bimodules so that the notation $HC_0(\tilde{Z})$ does make sense.

To construct the above square, we may assume that $Z$ is a $2$-term silting object 
since by part (3) of Theorem~7.2 of \cite{August18b}, the standard derived
equivalence given by $Z$ is the composition of equivalences
given by $2$-term tilting complexes and their inverses.
Let $A$ be the dg algebra obtained by applying Theorem~\ref{thm:liftingEquivalences}
to $B=\Ga'$, $H^0(B)=\Lambda'$, $C=\Lambda$ and $Z=Z$.
Its homologies $H^p(A)$ are finite-dimensional and vanish for $p>0$.
Thus, it is quasi-isomorphic to a dg algebra of the form $(\wh{T_l} V, d)$, where
$V$ is a graded bimodule whose components vanish in degrees $>0$ and
are finite-dimensional in all degrees $\leq 0$. So we may assume that
$A$ is in $\PCAlgc(l)$. By Theorem~\ref{thm:liftingEquivalences}, there is moreover
a dg $A$-$\Ga'$-bimodule $X$ yielding an
equivalence $?\lten_A X: \D(A) \iso \D(\Ga')$ and an isomorphism
$\phi: H^0(A) \to \Lambda$ such that we have an isomorphism
\[
\mbox{}_\phi Z \iso X\lten_{\Ga'} \Lambda'
\]
in the derived category $\D(A^{op}\ten \Lambda')$ and in particular a square,
commutative up to isomorphism
\[
\xymatrix{
A \ar[rr]^X \ar[d]_{H^0 A} & & \Ga' \ar[d]^{\Lambda'} \\
H^0(A) \ar[r]^{_\sim}_\phi & \Lambda \ar[r]_Z & \Lambda'.}
\]
By Corollary~\ref{cor:adapted-QP},
there is a quiver with potential $(Q'',w'')$ and a weak equivalence
$s:\Ga'' \to A$ from the associated Ginzburg algebra $\Ga''$ to $A$
such that the isomorphism $HC_0(\mbox{}_s X)$ takes
the class $[w'']$ to $[w']$. We define $\tilde{Z}=\mbox{}_{s}X$
and $\psi=\phi \circ H^0(s)$ to obtain the diagram
\[
\xymatrix{
\Ga'' \ar[d]  \ar[r]^s & A \ar[rr]^X \ar[d]_{H^0 A} & & \Ga' \ar[d]^{\Lambda'} \\
H^0(\Ga'') \ar[r]^\sim_{H^0s} & H^0(A) \ar[r]^-{_\sim}_-\phi & \Lambda \ar[r]_Z & \Lambda'.}
\]
So by construction, the isomorphism $HH_0(\mbox{}_{s} X)=HH_0(\mbox{}_\psi Z)$
takes $[w'']$ to $[w']$.
We may assume all potentials to contain no cycles of
length $\leq 2$ and then it follows that $\psi=\phi\circ H^0(s)$ induces an
isomorphism of quivers $Q'' \iso Q$. Indeed, it induces an isomorphism
in the Jacobian algebras and the vertices $i$
of the quiver $Q$ are in bijection with the isomorphism classes of simple modules $S_i$
of the pseudocompact algebra $\Lambda=H^0(\Ga)$  and the number of
arrows from $i$ to $j$ equals the dimension of the space
of extensions $\Ext^1_{\Lambda}(S_j,S_i)$.
By Corollary~\ref{iso-Ginz},
there is an isomorphism $\beta: \Ga'' \to \Ga$.
The dg bimodule $\mbox{}_{s\beta^{-1}}X$ yields
an algebraic triangle equivalence $\per(\Ga) \iso \per(\Ga')$. Such an
equivalence induces an equivalence between the subcategories
of dg modules with finite-dimensional homologies because
their objects $M$ are characterized as those for which
$\Hom(P,M)$ is finite-dimensional for any object $P$. 
Thus, the algebraic triangle equivalence $\per(\Ga) \iso \per(\Ga')$
induces an algebraic triangle equivalence in the cluster categories. 
Hence, by Lemma~\ref{eqdgenhancement}, it induces an algebraic
triangle equivalence $\D_{sg}(R) \iso \D_{sg}(R')$ and therefore an algebra 
isomorphism $\HoH^0(\D_{sg}(R))\cong \HoH^0(\D_{sg}(R^\p))$. By Theorem \ref{HHDsg=Tyurina}, we get an isomorphism $R\cong R^\p$.
\end{proof}

\section{Contractibility of rational curve}\label{sec:contractibility}
\subsection{dg $\kiu$-algebras}
In this section, we define dg-$\kiu$-algebras and study their properties. All the definitions and results can be adapted to the pseudo-compact case,
with notations appropriately replaced by their pseudo-compact counterparts.
\begin{definition}
Let $S$ be a commutative dg algebra. Denote by $\sC(S)$ the category of 
complexes of $S$-modules with the monoidal structure given by the tensor product 
over $S$. A dg $S$-algebra is an algebra in $\sC(S)$, i.e. a dg $S$-module with an $S$-bilinear multiplication and a unit.
\end{definition}

Let $k$ be a field, $S$ a commutative dg $k$-algebra and $\Ga$ a dg $S$-algebra.
By restriction, each dg $\Ga$-module becomes a dg $S$-module and the
morphism complexes between dg $\Ga$-modules are naturally dg $S$-modules.
Thus, the derived category $\D(\Ga)$ is naturally enriched over $\D(S)$.

\begin{definition}
Let $S=k[u^{-1}]$ be the commutative dg algebra with $\deg(u)=2, \deg(u^{-1})=-2$ and zero differential. We call a dg $k$-algebra  $A$  \emph{$\kiu$-enhanced} if $A$ is isomorphic to a dg $\kiu$-algebra in the homotopy category of dg $k$-algebras.
\end{definition}

Let us put $S=\kiu$ and $K=\kuiu$. For a dg $S$-module $M$, we have
a canonical isomorphism
\[
H^*(M\ten_S K)= H^*(M)\ten_S K.
\]
We call $M$ a {\em torsion module} if $M\ten_S K$ is acyclic. This happens if
and only if $H^*(M)$ is a torsion module, i.e. for each $m$ in $H^*(M)$,
there exists a $p\gg 0$ such that $m u^{-p}=0$.

Let $A$ be a dg $S$-algebra. The functor taking a dg $A$-module $M$ to
the dg $A\ten_S K$-module $M\ten_S K$ preserves quasi-isomorphisms.
Thus, it induces a functor $?\ten_S K: \D(A) \to \D(A\ten_S K)$. The kernel
of this functor consists of the dg $A$-modules which are torsion as dg
$S$-modules. We write $\D(A)_\ut$ for the kernel and $\per(A)_\ut$ for
its intersection with the perfect derived category $\per(A)$.

\begin{lemma} \label{lemma:exact-seq} 
We have exact sequences of triangulated categories
\[
\xymatrix{
0\ar[r] & \D(A)_\ut \ar[r] & \D(A) \ar[r] & \D(A\ten_S K) \ar[r] & 0}
\]
and
\[
\xymatrix{
0\ar[r] & \per(A)_\ut \ar[r] & \per(A) \ar[r] & \per(A\ten_S K) \ar[r] & 0.}
\]
\end{lemma}

\begin{proof} The restriction along $A \to A\ten_S K$ induces a fully faithful
right adjoint to $?\ten_S K: \D(A) \to \D(A\ten_S K)$. Thus, the latter functor
is a localization functor. By definition, its kernel is $\D(A)_\ut$ so that we
obtain the first sequence. To deduce the second one, it suffices to show
that the kernel of $\D(A) \to \D(A\ten_S K)$ is compactly generated.
Indeed, let $P$ be the cone over the morphism
\[
A \to \Sigma^{-2} A
\]
given by the multiplication by $u^{-1}$. Clearly $P$ is compact and we claim
that it generates the kernel. For this, it suffices to show that the right orthogonal
of $P$ in the kernel vanishes. Indeed, let $M$ be in the kernel. If
$\RHom_A(P,M)$ vanishes, then the morphism $\Sigma^2 M \to M$
given by the multiplication by $u^{-1}$ is a quasi-isomorphism. Thus,
$u^{-1}$ acts in $H^*(M)$ by an isomorphism. But on the other hand,
$H^*(M)$ is torsion. So $H^*(M)$ vanishes and $M$ is acyclic as was
to be shown.
\end{proof}

A dg $S$-algebra $A$ concentrated in non positive
degrees is {\em non degenerate} if the morphism
$A \to \Sigma^{-2} A$ given by the multiplication by $u^{-1}$ induces
isomorphisms $H^n(A) \iso H^{n-2}(A)$ for all $n\leq 0$.

\begin{lemma}
Let $A$ be a dg $k$-algebra concentrated in nonpositive degrees. Assume that $A$ is homologically smooth, that $H^n(A)$ is finite-dimensional for each $n\in\ZZ$ 
and that $A$ admits a 
non degenerate $\kiu$-enhancement. Then the subcategory 
$\D_{fd}(A)\subset \per(A)$ coincides with $\per(A)_\ut$.
\end{lemma}
\begin{proof}
Since $A$ is homologically smooth, $\D_{fd}(A)$ is contained in $\per(A)$ and
clearly it consists of torsion modules. Conversely, we know from the
proof of Lemma \ref{lemma:exact-seq} that $\per(A)_\ut$ is the thick
subcategory of $\D(A)$ generated by the cone $P$ over the morphism
$A \to \Sigma^{-2} A$ given by the multiplication by $u^{-1}$. Since
$A$ is non degenerate, the object $P$ lies in $\D_{fd}(A)$.
\end{proof}

\begin{prop}\label{Z2clustercat}
 Let $Q$ be a finite quiver with potential $w$ such that the associated Ginzburg algebra $\Ga = \fD(Q,w)$ has finite-dimensional Jacobi algebra. Assume that $\Ga$ is equipped with a nondegenerate 
 $\kiu$-enhancement. Then $\cC_\Ga$ is a $\ZZ/2$-graded 0CY triangulated 
 category, equivalent with the category of perfect modules over $\Ga\ot_{\kiu}\kuiu$ as $k[u,u^{-1}]$-enhanced triangulated categories. In particular, the Jacobi algebra $H^0(\Ga)$ is a symmetric Frobenius algebra.
\end{prop}
\begin{proof}
 We may assume that $\Gamma$ itself is a differential graded
$k[u^{-1}]$-algebra. The multiplication with $u^{-1}$ yields a functorial morphism
$\id \to \Sigma^{-2}$ of triangle functors $\D(\Gamma)\to \D(\Gamma)$ and
$\per(\Gamma)\to \per(\Gamma)$. This induces a functorial morphism
of triangle functors $\id\to\Sigma^{-2}$ in the cluster category $\cC_\Gamma$.
To check that it is invertible, it is enough to check that its action on
the cluster-tilting object $\Gamma\in \cC_\Gamma$ is invertible (since $\cC_\Gamma$
equals its thick subcategory generated by $\Gamma$). Now by our
assumption, the morphism $u^{-1}: \Sigma^2\Gamma \to \Gamma$ induces
isomorphisms in $H^n$ for $n\leq -2$ and the $0$-map for $n\geq -1$. Thus,
it induces a quasi-isomorphism $\Sigma^2\Gamma \to \tau_{\leq -2}\Gamma$.
We claim that the canonical morphism $\tau_{\leq -2}\Gamma \to \Gamma$ becomes
invertible in the cluster category. Indeed, the homology of its cone is
of finite total dimension since $H^p \Ga$ is of finite dimension for all 
integers $p$ by Lemma 2.5 of \cite{Am}.
It follows that $u$ induces an isomorphism
$\Sigma^2 \Gamma \iso \Gamma$ in $\cC_\Gamma$. The rest follows because
we have isomorphisms of $H^0(\Gamma)$-bimodules
\[
H^0(\Gamma) \iso \cC_\Gamma(\Gamma,\Gamma) \iso
\cC_\Gamma(\Gamma, \Sigma^2\Gamma) \iso D\cC_\Gamma(\Gamma,\Gamma).\]
\end{proof}

\subsection{$\kiu$-structures on Ginzburg algebras associated to contractible curves}
Recall that for a Jacobi-finite Ginzburg algebra $\Ga:=\fD(Q,w)$ associate to quiver $Q$ with only one vertex (and multiple loops), $H^0(\Ga)$ is self-injective. On the other hand, if $\Ga$ is non-degenerately $\kiu$-enhanced then $H^0(\Ga)$ is symmetric.

The simplest case is when $\Ga=k[t]$ with zero differential. It is $\kiu$-enhanced by setting $u^{-1}=t$. This is the derived deformation algebra for $\cO_C$ of a $(-1,-1)$-curve, which is always contractible.

\begin{prop}
Let $F=k\lgg x\rgg$ be the complete free algebra of one generator and $w\in F$ be an element with no constant term. Then $\Ga:=\fD(F,w)$ is $\kiu$-enhanced.
\end{prop}
\begin{proof}
A general element $w\in F$ is of the form:
\[
w=x^{n+1}+\text{higher order terms}.
\] When $n=1$ we are in the case of $(-1,-1)$ curve. We assume that $n\geq 2$.
The Jacobi algebra of $\fD(F,w)$ is isomorphic to $k[[x]]/(x^{n})$. It is always finite dimensional. Because $w=x^{n+1} \cdot u$ for some unit $u\in k[[x]]$, $[w]=0$ in $k[[x]]/(x^{n})$.  By Theorem \ref{ncSaito}, $w$ is right equivalent to $x^{n+1}$. Without loss of generality, we may assume $w=\frac{x^{n+1}}{n+1}$ to begin with. Then the Ginzburg algebra $\fD(F,w)$ is isomorphic to $k\lgg x,\theta,t\rgg$ with $dt=[x,\theta]$ and $d\theta=x^{n}$. It is easy to check that the two-sided differential ideal $(t,[x,\theta])$ is acyclic. As a consequence, the quotient morphism
\[
\Ga=(k\lgg x,\theta,t\rgg, d)\to \Ga^\p:=(k\lgg x,\theta,t\rgg/(t,[x,\theta]), d)
\] is a quasi-isomorphism of dg algebras. Note that $\Ga^\p$ is isomorphic to the complex
\[
\xymatrix{
\ldots k[[x]]\theta^3\ar[r]^d & k[[x]]\theta^2\ar[r]^d & k[[x]]\theta\ar[r]^d & k[[x]]\ar[r] & 0
}
\] where
\[
d(\theta^{2k})=0,~~~~~d(\theta^{2k+1})=x^n\theta^{2k}.
\]
Define the action of $u^{-1}$ on $\Ga^\p$ by multiplication by $\theta^2$. It is easy to check it makes $\Ga^\p$ a dg $\kiu$-algebra.
\end{proof}
From the above proposition, we see that the Ginzburg algebras associated to the ``one-loop quiver'' are essentially classified by the dimension of their Jacobi algebras. Moreover, they all admit $\kiu$-enhancements. 
If $\dim_k(H^0\Ga)=n$ for $n>1$, then $\Ga$ is equivalent to the derived deformation algebra of a floppable $(0,-2)$ curve of width $n$ (see \cite{Reid} for the geometric definition of width). The following corollary can be viewed as a  noncommutative counter part of the classification theorem of Reid \cite{Reid}.
\begin{corollary}
Let $C$ be a rational curve in a quasi-projective smooth CY 3-fold $Y$ of normal bundle $\cO_C\op \cO_C(-2)$. Denote its derived deformation algebra by $\Ga$. Then
\begin{enumerate}
\item[$(1)$] $C$ is movable if and only if $\Ga$ has infinite dimensional Jacobi algebra;
\item[$(2)$] If $C$ is rigid then it is contractible. The dimension of $H^0\Ga$ is equal to $n$ for $n>1$ if and only if the underlying singularity is isomorphic to the germ of hypersurface $x^2+y^2+u^2+v^{2n}=0$ at the origin.
\end{enumerate}
\end{corollary}

It is proved by Laufer (c.f. \cite{Pink83}) that a contractible rational curve in a CY 3-fold must have normal bundle of types $(-1,-1)$, $(0,-2)$ or $(1,-3)$. Donovan and Wemyss give an example of a rigid rational curve of type $(1,-3)$ that is not nc rigid (see Example 6.4 in \cite{DW15}). In their example, there exists a birational morphism that contracts a divisor containing the $(1,-3)$ curve. Kawamata asked  whether it is true that $C$ is contractible if it is nc rigid (see Question 6.6 of \cite{Ka18}). We formulate a conjecture in terms of the derived deformation algebra.
\begin{conjecture}\label{conj-contractibility}
Let $C\subset Y$ be a nc rigid rational curve in a smooth quasi-projective CY 3-fold. Denote its associated derived deformation algebra by $\Ga_C^Y$. Then $C$ is contractible if and only if $\Ga_C^Y$ is $\kiu$-enhanced.
\end{conjecture}
Note that one direction of the conjecture follows from our Theorem \ref{Gamma=truncation}. 
\begin{prop}\label{contractible=>kiu}
Let $C\subset Y$ be  a contractible rational curve  in a smooth quasi-projective CY 3-fold. Then $\Ga_C^Y$ is $\kiu$-enhanced.
\end{prop}
\begin{proof}
Denote $R$ for the ring of formal functions on the singularity underlying the contraction.
For simplicity, we denote the derived deformation algebra $\Ga_C^Y$ by $\Ga$. By Proposition \ref{eqdgenhancement}, $\cC_\Ga$ is quasi-equivalent to 
$\D_{sg}(R)$ as dg categories. Under the equivalence, the projection image of $\Ga$ is identified with the Cohen-Macaulay module $N\in \D_{sg}(R)\cong \ul{\CM}_R$. By Theorem \ref{Gamma=truncation}, $\Ga$ is isomorphic to $\tau_{\leq 0}\Lm_{dg}$ where $\Lm_{dg}$ is the dg endomorphism algebra of $N$ in $\ul{\CM}_R$. Because $R$ is a hypersurface ring, the dg category $\ul{\CM}_R$ carries a canonical $\ZZ/2$-graded structure (equivalently $\kuiu$-structure) by Eisenbud's theorem \cite{Eis80}. 
Therefore, $\Ga$ is $\kiu$-enhanced.
\end{proof}
We have already seen that in $(-1,-1)$ and $(0,-2)$ cases, the $\kiu$-structure on $\Ga$ can be computed explicitly. However, we don't have any explicit construction for the $(1,-3)$ case even though we know it must exist. We do have an explicit formula for the symmetric Frobenius structure on the CY tilted algebra $H^0\Ga$ in term of the residue map of matrix factorizations (see \cite{HT16}).

\section{Appendix: Serre duality for sheaves and modules}
In this Appendix, we give two proofs of the link the between the
inverse dualizing sheaf on a smooth quasi-projective variety $Y$ and 
the inverse dualizing bimodule for the derived endomorphism
algebra of any perfect generator of $\D(\Qcoh(Y))$.

The proof in section~\ref{ss:Gaitsgory} is based on the work of
Gaitsgory \cite{Gaitspaper} and Gaitsgory--Rosenblyum \cite{Gaits11}.
The proof in section~\ref{ss:BFN} is essentially taken from Example~2.7 of 
\cite{KinjoMasuda21}. It combines Grothendieck duality with the results
of Ben-Zvi--Francis--Nadler's article \cite{BenZviFrancisNadler10}.
Both proofs rely on the foundational work of To\"en, Joyal, Lurie
and many others.

\subsection{Dg Serre duality, after Gaitsgory} \label{ss:Gaitsgory}
We recall a result of Gaitsgory and Gaitsgory--Rosen\-blyum which could be viewed as Serre duality for the 
dg category of coherent sheaves. They use the category of ind coherent sheaves which behaves 
better on singular spaces. Though we only need the smooth case, we recall the basics of ind coherent sheaves and their properties. We will follow \cite{Gaits11} and \cite{Gaitspaper}.

We consider (quasi-)coherent sheaves on quasi-compact separated Noetherian schemes.
We write $\Qcoh$ for the dg category of (fibrant replacements of) unbounded complexes of quasi-coherent sheaves and $\coh$ for its full dg subcategory of complexes with coherent cohomology and bounded cohomological amplitude.
Let $\IndQ$ be the ind-completion of $\coh$. Thus, the dg category $\IndQ$ is quasi-equivalent to the dg enhanced
derived category of the (essentially) small dg category $\coh$. There is a natural functor
$\Psi: \IndQ\to \Qcoh$ which commutes with coproducts and, restricted to $\coh$, becomes the inclusion
 $\coh\to \Qcoh$. If $X$ is smooth then $\Psi_X: \IndQ(X)\to \Qcoh(X)$ is an equivalence \cite[Chapter 4, Lemma 1.1.3]{Gaits11}. Let $f: X\to Y$ be a proper morphism. There exists a continuous pushforward functor 
 $f^{\IndQ}_*:\IndQ(X)\to \IndQ(Y)$ with a commutative diagram
\[
\xymatrix{
\IndQ(X)\ar[r]^{f^\IndQ_*} \ar[d]^{\Psi_X} & \IndQ(Y)\ar[d]^{\Psi_Y}\\
\Qcoh(X)\ar[r]^{f_*} & \Qcoh(Y)
}
\] See \cite[Ch.~4, Prop.~2.1.2]{Gaits11}. Moreover, since $f$ is proper $f^{\IndQ}_*$ sends $\coh(X)$ to $\coh(Y)$. The above commutative diagram is compatible with tensor products by
quasi-coherent complexes (\cite[Ch.~4, Prop.~2.1.4]{Gaits11}). The pushforward functor $f^{\IndQ}_*$ admits a \emph{continuous} right adjoint functor $f^!:\IndQ(Y)\to \IndQ(X)$ (\cite[Ch.~4, 5.1.5]{Gaits11}) and it is compatible with 
the tensor product with quasi-coherent complexes (\cite[Ch.~4, 5.1.7]{Gaits11}). To distinguish $f^!$ from the right adjoint of $f_*$ in $\Qcoh$, 
we denote the latter by $f^{\Qcoh, !}$. We have the following comparison theorem between these two 
functors (\cite[Ch.~4, Lemma 5.1.9]{Gaits11}): There is a commutative diagram 
\begin{align}\label{CD:!pullback}
\xymatrix{
\IndQ(X)^+\ar[r]^{\Psi_X}  & \Qcoh(X)^+\\
\IndQ(Y)^+\ar[u]^{f^!}\ar[r]^{\Psi_Y} & \Qcoh(Y)^+\ar[u]^{f^{\Qcoh,!}}
}
\end{align}  where the superscript + refers to the subcategory consisting of objects whose cohomological amplitude is bounded below. Note that the similar diagram with $+$ removed is not commutative (see \cite[Ch.~4, Remark 5.1.10]{Gaits11}). This result shows that we can compute $f^!\cE$ via $f^{\Qcoh,!}\cE$ for $\cE\in\coh(Y)$. The left adjoint of $f_*$, denoted by $f^{\IndQ, *}$ exists. It is compatible with tensor products by quasi-coherent complexes and satisfies a similar 
commutative diagram 
as diagram \ref{CD:!pullback} with $+$ removed (\cite[Ch.~4, Prop.~3.1.6]{Gaits11}).

For a scheme $X$ of finite type, by \cite[Ch.~5, Th.~4.2.5]{Gaits11}
there is a canonical equivalence 
\[
\DD_X: \IndQ(X)^\vee\simeq \IndQ(X)
\] 
such that 
\[
\IndQ(X)^{\vee\vee}\simeq \IndQ(X)^\vee\simeq \IndQ(X).
\] The dual category $\cC^\vee$ of a dualizable dg category $\cC$ can be identified with the category of continuous dg functors from $\cC$ to $dg_k$ where $dg_k$ is the category of dg $k$-modules.

\begin{theorem}\label{Serredg}
Let $X$ be a smooth and proper $k$-scheme of dimension $d$ for $k=\CC$. Let $M,N$ be objects of $\coh(X)$. There is a bifunctorial quasi-isomorphism
\[
D\Hom^{dg}_{\coh(X)}(M,N)\simeq \Hom^{dg}_{\coh(X)}(N,M\ot \omega_X[d]).
\]
\end{theorem}
\begin{proof}
We will interpret 
\[
M^\vee\bt N\mapsto D\Hom^{dg}_{\coh(X)}(M,N)
\] and 
\[
M^\vee\bt N\mapsto \Hom^{dg}_{\coh(X)}(N,M\ot \omega_Y[d])
\] as dg functors from $\coh(X\times X)$ to $dg_k$ and show there is a natural isomorphism between them.
By \cite[Chapter 3, Proposition 3.1.7]{Gaits11}, we have an equivalence
$\Qcoh(X\times X)\simeq \Qcoh(X)\ot \Qcoh(X)$. 
Therefore we obtain a $\coh(X)^{op}\ot \coh(X)$-module structure by letting $M$ and $N$ vary in $\coh(X)$.

Let $f:X\to X\times X$ be the diagonal map and $p:X\to \Spec k$ be the counit map. By \cite[Chapter 5, Theorem 4.2.5]{Gaits11} there is a commutative diagram of dg functors (both the upper and the lower square commute, c.f. \cite[9.2.3]{Gaitspaper}.)
\[
\xymatrix{
\IndQ(k)^\vee\ar[r]^{\DD_k} & \IndQ(k) \\
\IndQ(X)^\vee \ar[r]^{\DD_X}\ar[u]^{(p^!)^\vee} & \IndQ(X)\ar[u]^{p_*} \\
\IndQ(X\times X)^\vee\ar[r]^{\DD_{X\times X}}\ar[u]^{f_*^\vee} & \IndQ(X\times X)\ar[u]^{f^!}
}
\]
By \cite[Corollary 9.5.9]{Gaitspaper}, if we restrict to $\coh(X)$ (resp. $\coh(k)$ and 
$\coh(X\times X)$) we get a commutative diagram
\[
\xymatrix{
\coh(k)^{op}\ar[r]^{\DD_k} & \coh(k) \\
\coh(X)^{op} \ar[r]^{\DD_X}\ar[u]^{(p_*)^{op}} & \coh(X)\ar[u]^{p_*} \\
\coh(X\times X)^{op}\ar[r]^{\DD_{X\times X}}\ar[u]^{(f^*)^{op}} & \coh(X\times X)\ar[u]^{f^!}
} 
\]
Then we have a natural isomorphism
\[
p_*\circ f^!\circ \DD_{X\times X}\simeq \DD_k\circ (p_*)^{op}\circ (f^*)^{op}.
\]
Since $X$ is smooth coherent sheaves are dualizable. In this case 
\[
\DD_X(E)=E^\vee\ot \omega_X[d]
\] for $E\in \coh(X)$ (\cite[Lemma 9.5.5]{Gaitspaper}).
\footnote{In \cite{Gaitspaper}, Gaitsgory denotes by $\omega_X$ the $!$-pullback $p^!k$ for $p:X\to k$, which differs from the the standard notion of dualizing complex by $[d]$.}

Given an object $M^\vee\bt N$ in $\coh(X\times X)$, we may compute
\begin{align*}
p_*\circ f^!\circ \DD_{X\times X}(M^\vee\bt N)&\cong p_*(f^!((M\ot\omega_X[d])\bt (N^\vee\ot\omega_X[d])))\\
&\cong p_*(f^*(M\bt N^\vee)\ot \omega_X[d])\\
&\cong \Hom^{dg}_{\coh(X)}(N,M\ot \omega_X[d])
\end{align*} 
In the last step, we use the condition that $X$ is smooth and proper. On the other hand
\begin{align*}
\DD_k\circ (p_*)^{op}\circ (f^*)^{op}(M^\vee\bt N)&\cong 
\DD_k\Hom^{dg}_{\coh(X)}(M,N)\\
&\cong D\Hom^{dg}_{\coh(X)}(M,N)
\end{align*} 
\end{proof}

\subsection{Inverse dualizing sheaves and bimodules, after Ben-Zvi--Francis--Nadler} \label{ss:BFN}
The following proposition is a consequence of Theorem~4.7 and Corollary~4.8 of Ben-Zvi--Francis-Nadler's paper 
\cite{BenZviFrancisNadler10} combined with Grothendieck duality. The proof we give is an elaboration on Example~2.7 
of \cite{KinjoMasuda21}.

\begin{prop}[Kinjo--Masuda \cite{KinjoMasuda21}] \label{prop:generators}
Let $Y$ be a smooth quasi-projective variety of dimension~$d$.
Let $G$ be a perfect generator of $\D(\Qcoh Y)$ and $B=\RHom(G,G)$.
Then there is a canonical equivalence
\[
\D(\Qcoh (Y \times Y)) \to \D(B \ten B^{op})
\]
taking $\Delta_*(\Sg^{-d}\omega_Y^{-1})$ to the inverse dualizing complex
$\Theta_B=\Hom_{B^e}(B,B^e)$ of $B$ and $\Delta_*(\mathcal{O}_Y)$ to the identity bimodule $B$. 
In particular, if $Y$ is $d$-Calabi-Yau, then $B$ is bimodule $d$-Calabi-Yau. 
\end{prop}

\begin{remark} One could make further use of the results of \cite{BenZviFrancisNadler10}
to show more precisely that the derived endomorphism algebra of $G\boxtimes G^\vee$
in $\D(Y\times Y)$ is quasi-isomorphic to $B^e$ and that the canonical
equivalence of the theorem is given by
\[
\RHom(G \boxtimes G^\vee): \D(Y \times Y) \iso \D(B^e).
\]
\end{remark}

\begin{proof} Let $k=\CC$. We mostly work in the $\infty$-category $\St$ of $k$-linear
stable presentable $\infty$-categories whose $1$-morphisms are cocontinuous
$k$-linear exact $\infty$-functors (equivalently: left adjoints of $k$-linear 
exact $\infty$-functors). We recommend section~2 of \cite{KinjoMasuda21} for
a concise but readable introduction to this setting. Each (large) pretriangulated dg $k$-category 
$\cA$ such that $H^0(\cA)$ has arbitrary (set-indexed) coproducts and is compactly generated 
gives rise to an object of $\St$; each dg functor $F: \cA \to \cB$ between two such categories such that $H^0(F)$
commutes with arbitrary coproducts gives rise to a $1$-morphism in 
$\St$ and this $1$-morphism is an equivalence iff $H^0(F)$ is an equivalence.
To make these statements more precise, let us denote by $\dgcat_k^{pretr,\oplus}$ the category
\begin{itemize}
\item[-] whose objects are the pretriangulated dg categories $\cA$ such that $H^0(\cA)$ has arbitrary coproducts 
and is compactly generated and 
\item[-] whose morphisms are the dg functors $F$ such that $H^0(F)$ commutes with arbitrary
coproducts. 
\end{itemize}
We then have a canonical $\infty$-functor (where on the left, we write
the ordinary category instead of its nerve)
\[
\mathrm{can}: \dgcat_k^{pretr,\oplus} \longrightarrow \St.
\]
Each object $\cX$ of $\St$ has an underlying $\infty$-category $\cX_\infty$ (obtained
by forgetting the $k$-linear structure). It is a stable presentable $\infty$-category
and its $1$-categorical truncation $\tau(\cX_\infty)$ is naturally a triangulated category
with arbitrary coproducts, cf.~section~1.4.4 of \cite{LurieHA}. 
Similarly, if $f: \cX \to \cY$ is a $1$-morphism of $\St$, then
its underlying $\infty$-functor $f_\infty: \cX_\infty \to \cY_\infty$ is exact and cocontinuous and
its $1$-categorical truncation $\tau(f_\infty): \tau(\cX_\infty) \to \tau(\cY_\infty)$ is naturally a triangulated 
functor and commutes with arbitrary coproducts. To state these facts more precisely, let us
denote by $\Tria^{\oplus}$ the category whose objects are the triangulated categories
with arbitrary coproducts and whose morphisms are the triangle functors which commute
with arbitary coproducts;  let us denote by $\Pr_{st}^{L}$ Lurie's $\infty$-category of presentable 
stable $\infty$-categories whose $1$-morphisms are the exact left adjoint $\infty$-functors. 
Then we have $\infty$-functors
\[
\xymatrix{
\dgcat_k^{pretr,\oplus} \ar[r]^-{\mathrm{can}} &
\St \ar[r]^{?_\infty} &
\Pr_{st}^L \ar[r]^\tau & 
\Tria^\oplus
}.
\]
The $\infty$-category $\Pr_{st}^L$ is by definition a subcategory of the
$\infty$-category $\Pr_{st}$ of presentable stable $\infty$-categories whose
$1$-morphisms are all exact $\infty$-functors.  We will use the following facts: 
\begin{itemize}
\item[a)] If $f: \cX \to \cY$ is a $1$-morphism of $\Pr_{st}$, then the values of
$f$ and of $\tau(f)$ on objects are equal in $\tau(\cY)$. 
\item[b)] A $1$-morphism $f: \cX \to \cY$ of $\Pr_{st}$ belongs to $\Pr_{st}^L$ 
iff its truncation $\tau(f): \tau(\cX) \to \tau(\cY)$ is a (triangle) functor which 
commutes with arbitrary coproducts, cf.~Prop.~1.4.4.1 (2) of \cite{LurieHA}.
\item[c)] A $1$-morphism $f: \cX \to \cY$ of $\Pr_{st}$ admits a left adjoint 
$f_\lambda: \cY \to \cX$ in $\Pr_{st}$ iff the (triangle) functor $\tau(f)$ admits a left adjoint 
$\tau(f)_\lambda$. By definition, in this case, the $1$-morphism $f_\lambda$ belongs to $\Pr_{st}^L$.
Moreover, the functor $\tau(f_\lambda)$ is isomorphic to $\tau(f)_\lambda$ since truncation
preserves adjunctions.
\end{itemize}

For a dg algebra $A$, we denote by $\iD(A)$  the object of $\St$ corresponding to the
dg category of $K$-projective dg $A$-modules (i.e. dg $A$-modules which are
homotopy equivalent to cofibrant dg $A$-modules).  Notice that $\iD(A)$
denotes an object of $\St$ whereas $\D(A)$ denotes the ordinary derived
category (with its triangulated structure). The ordinary derived category
$\D(A)$ is naturally equivalent to the truncation $\tau(\iD(A)_\infty)$.

Let $X$ be a quasi-projective variety (or more generally, a quasi-projective
separated scheme). Let $\iD(X)$ denote the object of $\St$ corresponding to the
dg category of fibrant complexes of quasi-coherent sheaves on $X$.
Then the truncation $\tau(\iD(X)_\infty)$ identifies with
the derived category $\D(\Qcoh(X))$. Let $T$ be a compact generator
of the unbounded derived category of quasi-coherent sheaves on $X$.
Suppose that $T$ is a fibrant complex of quasi-coherent sheaves and
let $A=\Hom(T,T)$ be its dg endomorphism algebra. The homotopy category
of fibrant complexes of quasi-coherent sheaves on $X$ is compactly
generated by $T$ and so the dg $A$-module $\Hom(T,I)$ is $K$-projective
for all fibrant complexes $I$ of quasi-coherent sheaves on $X$. Moreover,
the dg functor $\Hom(T,?)$ induces a quasi-equivalence between the
homotopy categories of fibrant complexes on $X$ and $K$-projective
dg $A$-modules. Thus, the dg functor $\Hom(T,?)$ induces a
$1$-morphism $\iD(X) \to \iD(A)$ in $\St$. Its truncation $\tau(\Hom(T,?)_\infty)$
identifies with $\RHom(T,?): \D(\Qcoh(X)) \to \D(\Qcoh(Y))$, which is an
equivalence by our assumption on $T$. Thus, we have an
equivalence $\iD(X) \to \iD(A)$ in $\St$.

As detailed in section~4.8 of \cite{LurieHA} and section~2.2 of \cite{KinjoMasuda21}, 
the $\infty$-category $\St$ is endowed with a symmetric monoidal structure. 
The unit of the monoidal structure on $\St$ is the $k$-linear (symmetric monoidal) $\infty$-category $\iD(k)$.
If $A_1$ and $A_2$ are two dg algebras, then, by part (2) of Prop.~4.1 of \cite{BenZviFrancisNadler10}, 
we have a canonical equivalence
\[
\xymatrix{
\iD(A_1) \ten \iD(A_2) \ar[r]^-{_\sim} & \iD(A_1\ten A_2)
}
\]
induced by the dg functor taking a pair $(L,M)$ of cofibrant dg modules
to the cofibrant $A_1\ten A_2$-module $L \ten M$.

For any dg algebra $A$, the object $\iD(A)$ is dualizable in the monoidal $\infty$-category $\St$ 
and its dual is equivalent to $\iD(A^{op})$, cf.~Prop.~4.3 (3) of \cite{BenZviFrancisNadler10}. 
More precisely, the duality between  $\iD(A)$ and $\iD(A^{op})$ is given by the evaluation morphism 
\[
\iD(A^{op})\ten \iD(A) \to \iD(k)
\]
and the co-evaluation morphism
\[
\iD(k) \to \iD(A)\ten \iD(A^{op}).
\]
The composition of the co-evaluation morphism with the canonical morphism 
$\iD(A)\ten \iD(A^{op}) \iso \iD(A\ten A^{op})$ is the morphism
\[
\iD(k) \to \iD(A \ten A^{op})
\]
induced by the dg functor taking a complex $V$ to $V \ten A$, where $A$ is the identity bimodule.
The evaluation morphism is  obtained as the composition
\[
\xymatrix{
\iD(A^{op}) \ten \iD(A) \ar[r]^-{_\sim} & \iD(A^{op}\ten A) \ar[r]^-{\ev^A} &  \iD(k), 
}
\]
where the second morphism is induced by the dg functor taking
a cofibrant dg bimodule $M$ to $M \ten_{A^e} A$. Let us 
abbreviate $\ev^A_\infty$ by $\Ev^A$. Then
the truncation $\tau(\Ev^A)$ identifies with the derived functor
\[
?\lten_{A^e} A : \D(A^{op}\ten A) \to \D(k).
\]
When $A$ is smooth, this functor is isomorphic to $\RHom_{A^e}(\Theta,?)$,
where $\Theta=\Theta_A=\RHom_{A^e}(A, A^e)$ is the inverse dualizing bimodule.
Thus, it admits the left adjoint 
\[
?\ten \Theta: \D(k) \to \D(A^{op} \ten A).
\]
Hence, by fact c) above, the $\infty$-functor $\Ev^A$ admits a left
adjoint $\Ev^A_\lambda$ and $\tau(\Ev^A_\lambda)$ is isomorphic
to $?\lten_{A^e} \Theta$. In particular, the left adjoint $\Ev^A_\lambda$ 
takes $k$ to $\Theta$ (fact a) above).

Let $X_1$ and $X_2$ be quasi-projective varieties (or more generally,
quasi-compact, separated schemes).  They are in particular perfect
stacks in the sense of \cite{BenZviFrancisNadler10}.
By Theorem~4.7 of [loc.~cit.], we have a canonical equivalence
\[
\iD(X_1)\ten \iD(X_2) \iso \iD(X_1 \times X_2).
\]
Fix a quasi-projective variety $X$.
By Corollary~4.8 of \cite{BenZviFrancisNadler10} and the first three lines of 
its proof, the object $\iD(X)$ becomes its own dual and the evaluation morphism 
is given by a composition
\[
\xymatrix{
\iD(X) \ten \iD(X) \ar[r]^-{_\sim} &  \iD(X \times X) \ar[r]^-{\ev^X} &  \iD(k)},
\]
where the truncation $\tau(\ev^X_\infty)$ identifies with the composition of derived functors
\[
\xymatrix{
\D(X\times X) \ar[r]^-{\Delta^*} &  \D(X) \ar[r]^{\pi_*} & \D(k).
}
\]
The co-evaluation morphism is the composition
\[
\xymatrix{
\iD(k) \ar[r] & \iD(X\times X) & \iD(X)\ten \iD(X) \ar[l]_-{_\sim}},
\]
where the first morphism is induced by the dg functor taking
a complex $V$ to $V\ten I$, where $I$ is an injective
resolution of $\Delta_*(\mathcal{O}_X)$.

Now consider the smooth quasiprojective variety $Y$ of the claim.
In this case, by Grothendieck duality,
the functor $\Delta^*: \D(Y\times Y) \to \D(Y)$ admits the left adjoint $\Delta_*(?)\ten \Sg^{-d}\omega_Y^{-1}$.
Thus, the composed functor
\[
\xymatrix{
\D(Y\times Y) \ar[r]^-{\Delta^*} &  \D(Y) \ar[r]^{\pi_*} & \D(k).
}
\]
admits the left adjoint
\[
\Delta_*(\pi^*(?)\ten \Sg^{-d}\omega_Y^{-1})
\]
which takes $k$ to $\Delta_*(\Sg^{-d}\omega_Y^{-1})$.
As above, it follows that the $\infty$-functor $\Ev^Y=\ev^Y_\infty$
admits a left adjoint $\Ev^Y_\lambda$ and that this left adjoint takes $k$ to $\Delta_*(\Sg^{-d} \omega_Y^{-1})$.

As we have seen above, under our hypotheses, the dg functor taking a fibrant complex of
quasi-coherent sheaves $C$ to $\Hom(G,C)$ induces an equivalence
\[
\iD(Y) \iso \iD(B).
\]
We deduce an equivalence
\[
\xymatrix{
\iD(Y\times Y)  & \iD(Y) \otimes \iD(Y) \ar[l]_-{_\sim} \ar[r]^-{_\sim} &
\iD(Y)\otimes \iD(Y)^\vee \ar[r]^-{_\sim} &
 \iD(B) \otimes \iD(B^{op}) \ar[r]^-{_\sim} & \iD(B\otimes B^{op})}.
\]
The intrinsic descriptions above show that under the $1$-categorical
truncation of this equivalence, the object $\Delta_*(\Sg^{-d}\omega_Y^{-1})$ 
corresponds to the inverse dualizing complex $\Theta_B=\Hom_{B^e}(B,B^e)$ of $B$.
Moreover, the object $\Delta_*(\mathcal{O}_Y)$ corresponds to the
identity bimodule $B\in \D(B^e)$ as we see by examining the co-evaluation morphisms.
Thus, an isomorphism $\omega_Y \iso \mathcal{O}_Y$ yields an isomorphism
$\Theta \iso \Sigma^{-d} B$ in $\D(B^e)$.
\end{proof}

\end{document}